%% file: template.tex
\documentclass{m2an}

\usepackage{bm}
\usepackage{algorithm}% http://ctan.org/pkg/algorithm
\usepackage{algpseudocode}
\usepackage{color}
\usepackage{systeme}
\usepackage{comment}
\usepackage{scalerel}
\usepackage{cases}
\usepackage{graphicx}
\usepackage{grffile}
\usepackage{multirow}
\usepackage{pdflscape}
\usepackage{amsmath,amsfonts,amssymb}
\usepackage{bm}
\usepackage{booktabs}
\usepackage[plainpages = false, pdfpagelabels,
	     bookmarks,
	     bookmarksopen = true,
	     bookmarksnumbered = true,
	     breaklinks = true,
	     linktocpage,
	     pagebackref,         % shows hyperrefs for usage of refs
	     colorlinks = true,  % was true
         hypertexnames=false,   %destination with the same identifier error fixed
	     linkcolor = green!50!black,
	     urlcolor  = blue!50!black,
	     citecolor = blue!50!black,
	     anchorcolor = blue,
	     hyperindex = true,
	     hyperfigures
	     ]{hyperref}
\usepackage[capitalise,noabbrev]{cleveref}
\usepackage{algorithm}
\usepackage{algorithmicx}
\usepackage{algpseudocode}
\usepackage{mathtools}
\usepackage{comment}
\usepackage{tikz}
\newlength{\figureheight}
\newlength{\figurewidth}
\usepackage{tikzsymbols}
\usepackage{pgfplots}
\pgfplotsset{compat=newest} %use latest pgf version
\usetikzlibrary{plotmarks}
\usetikzlibrary{arrows.meta}
\usetikzlibrary{positioning}
\usepgfplotslibrary{patchplots}
\usepackage{grffile}
\usepackage{multirow}
\usepackage{pdflscape}
 \usepackage{fontawesome}
\pgfplotsset{every axis/.append style={
                    label style={font=\scriptsize},
                    tick label style={font=\scriptsize},
                    legend style={font=\scriptsize}
                    }}
\pgfplotsset{compat=newest}
\pgfplotsset{plot coordinates/math parser=false}
\pgfplotsset{grid style={dotted,gray}}
\usepgfplotslibrary{groupplots}
\pgfdeclarelayer{background}
\pgfdeclarelayer{foreground}
\pgfsetlayers{background,main,foreground}
\usepackage{soul}
\usepackage{stmaryrd}
\usepackage{subcaption}
\usepackage{textcomp}
\usepackage{xcolor}
\usepackage{amsmath}
\usepackage{dsfont}
\usepackage{graphicx}
\usepackage{hyperref}
\usepackage{markdown}

\newcommand\norm[1]{\left\lVert#1\right\rVert}

\newcommand{\matr}[1]{#1}

\newtheorem{example}{Motivational example}

\makeatletter
\patchcmd{\ALG@step}{\addtocounter{ALG@line}{1}}{\refstepcounter{ALG@line}}{}{}
\newcommand{\ALG@lineautorefname}{Line}
\makeatother

\definecolor{darkgreen}{rgb}{0.0,0.5,0.0}

\begin{document}
%%-----------------------------
%%      the top matter
%%-----------------------------
\title{Optimal control for a class of linear transport-dominated systems via the shifted proper orthogonal decomposition}\thanks{We gratefully acknowledge the support of the Deutsche Forschungsgemeinschaft (DFG) as part of GRK2433 DAEDALUS. We also acknowledge the financial support from the SFB TRR154 under the sub-project B03.}% At most 5 thanks
\author{Tobias Breiten}\address{Institute of Mathematics, Technische  Universit\"at Berlin, Stra\ss e des 17. Juni 136, Berlin, 10623, Berlin, Germany;  \email{tobias.breiten@tu-berlin.de}, \email{burela@tnt.tu-berlin.de}, \email{pschulze@math.tu-berlin.de}}
\author{Shubhaditya Burela}\sameaddress{1}
\author{Philipp Schulze}\sameaddress{1}
%
% \date{...}
%
\begin{abstract} 
Solving optimal control problems for transport-dominated partial differential equations (PDEs) can become computationally expensive, especially when dealing with high-dimensional systems. 
To overcome this challenge, we focus on developing and deriving reduced-order models that can replace the full PDE system in solving the optimal control problem.
Specifically, we explore the use of the shifted proper orthogonal decomposition (POD) as a reduced-order model, which is particularly effective for capturing high-fidelity, low-dimensional representations of transport-dominated phenomena.
Furthermore, we propose two distinct frameworks for addressing these problems: one where the reduced-order model is constructed first, followed by optimization of the reduced system, and another where the original PDE system is optimized first, with the reduced-order model subsequently applied to the optimality system.
We consider a 1D linear advection equation problem and compare the computational performance of the shifted POD method against the conventional methods like the standard POD when the reduced-order models are used as surrogates within a backtracking line search.
\end{abstract}

\maketitle
%%-----------------------------
%%      your text
%%-----------------------------

\section{Introduction}
In this work, we discuss the use of reduced-order methods for optimal control problems for transport-dominated phenomena, associated for example with a linear first-order hyperbolic \textit{partial differential equation} (PDE) of the form
\begin{equation}\label{eq_def:PDE_continuous}
\begin{aligned}
    \partial_t y(x, t) &= \mathcal{A} y(x, t) + (\mathcal{B}u(t))(x),  &(x, t) \in \mathcal{Q}:=\Omega \times [0,t_f]\\
    y(x, 0) &= y_0(x), & x \in \Omega:=[0,\ell]\subset \mathbb R,
\end{aligned}
\end{equation}
minimizing a standard quadratic tracking type cost functional 
\begin{equation}\label{eq_def:FOM_costFunc_continuous}
    \mathcal{J}(y, u) = \frac{1}{2} \int^{t_f}_0\int_{\Omega} \bigl(y(x, t) - y_\mathrm{d}(x, t)\bigr)^2 \mathrm{d} x\, \mathrm{d}t + \frac{\mu}{2} \int_0^{t_f}\norm{u(t)}^2_{\mathbb{R}^{n_c}}\,\mathrm{d}t.
\end{equation} 
Our main focus will be on (a spatially discretized version of) the case $\mathcal{A}=-v \tfrac{\partial}{\partial x}$, i.e., the linear transport equation. We further make the following assumptions. The number of controls $n_c$ is finite and for fixed $t$ we search for $u(t) \in U_{ad} := \mathbb{R}^{n_c}$. For simplicity, we assume the control operator $\mathcal{B}$ to be bounded, i.e., it holds that $\mathcal{B}\in \mathcal{L}( \mathbb{R}^{n_c},L^2(\Omega))$. For the desired state, we assume $y_{\mathrm{d}} \in L^2(0,t_f;\Omega)$ and the regularization parameter  $\mu>0$ to be positive.

Optimal control problems of an abstract form as in \eqref{eq_def:PDE_continuous}-\eqref{eq_def:FOM_costFunc_continuous} are well understood and detailed treatises can be found in, e.g., the monographs \cite{hinze_optimization_2009,troltzsch_optimal_2010}. First-order necessary optimality conditions can be derived in a straightforward way by means of the formal Lagrange method and lead to the adjoint equation
\begin{equation}\label{eq_def:PDE_adjoint_continuous}
\begin{aligned}
    -\partial_t \lambda(x, t) &= \mathcal{A}^* \lambda(x, t) +y(x, t) - y_\mathrm{d}(x, t), \\
    \lambda(x, t_f) &= 0,\\
    \mu u(t) + \mathcal{B}^* \lambda(\cdot, t) &= 0,
\end{aligned}
\end{equation}
where $\lambda$ denotes the adjoint or co-state and $\mathcal{B}^*\colon L^2(\Omega) \rightarrow \mathbb{R}^{n_c}$ is the adjoint of $\mathcal{B}$. 

PDE-constrained optimal control problems are often hard to solve numerically as they lead to large-scale optimization problems, in particular for higher spatial dimensions and/or in cases where a fine discretization is required. 
For this reason, one is interested in mathematical methods that reduce the complexity of the underlying dynamical system and to speed up its simulation by using \textit{reduced-order models} (ROMs), see, e.g., \cite{benner_model_2017} for an overview also including linear quadratic optimal control problems. 
While for parabolic problems, the use of ROMs for linear quadratic control problems is rather standard by now, hyperbolic problems still pose significant challenges as they exhibit a slower decay of the Kolmogorov $n$-width \cite{greif_decay_2019}, rendering many conventional reduced order modeling methods ineffective.

The main goal of this work is to construct ROMs for \eqref{eq_def:PDE_continuous} and \eqref{eq_def:PDE_adjoint_continuous} specifically utilizing the \textit{shifted proper orthogonal decomposition} (sPOD) \cite{reiss_shifted_2018} method. This method alleviates the issue of the slow decay of the Kolmogorov $n$-width, and we aim to study its impact within the context of optimal control. 
Over the past three decades, significant progress has been made in developing efficient \textit{model order reduction} (MOR) methods for optimal control problems, with a particular focus on projection-based methods. 
These methods rely on projecting the dynamical system onto subspaces consisting of basis elements that contain characteristics of the expected solution. 
The \textit{reduced basis} (RB) \cite{benner_model_2017} method is one of such kind which constructs a low-dimensional reduced basis space that is spanned by the snapshots of the original solution. The reduced basis approximation is then efficiently obtained by Galerkin projection onto this space. 
Early research on RB methods for optimal control, such as \cite{ito_reduced_2001, quarteroni_reduced_nodate}, focused on flow problems, while later studies \cite{negri_reduced_2013, karcher_certified_2018} extended RB methods to the optimal control of parametrized problems.
Another MOR technique, \textit{balanced truncation} (BT), transforms the state-space system into a balanced form, where the controllability and observability Gramians are diagonal and equal, allowing states that are hard to observe and hard to reach to be truncated. 
Applications of balanced truncation to PDE control problems have been explored in works such as \cite{benner_balancing-related_2013, de_los_reyes_balanced_2011}.

Probably the most widely used MOR method for linear or non-linear optimal control problems is \textit{proper orthogonal decomposition} (POD). 
Early applications of POD in optimal control can be traced back to \cite{ravindran_reduced-order_2000, kunisch_control_1999}, where its use in fluid flow control was explored. 
POD has also been effectively employed to compute reduced-order controllers \cite{ravindran_adaptive_2002} with nonlinear observers \cite{atwell_reduced_nodate} and to design model predictive controllers \cite{alla_asymptotic_2015, ghiglieri_optimal_2014}. 
In many optimal control problems, it is often necessary to compute a feedback control law instead of a static optimal control. 
POD-based reduced order models have been successfully applied to address such problems as well \cite{kunisch_hjb-pod-based_2004}. 
Additionally, there has been considerable research on error analysis for POD-based reduced-order models in optimal control. 
This includes analysis for nonlinear dynamical systems \cite{hinze_proper_2005}, abstract LQ systems \cite{hinze_error_2008}, parabolic problems \cite{gubisch_pod_2014, studinger_numerical_2012}, and elliptic problems \cite{kahlbacher_pod_2012}. 
A comparison of a posteriori error estimators for RB and POD methods in linear-quadratic optimal control problems was conducted in \cite{tonn_comparison_2011}. 
Furthermore, numerous successful applications of POD-based optimal control have been demonstrated in fields such as fluid dynamics and aerodynamic shape optimization \cite{LY2001223, choi_gradient-based_2020}, chemistry \cite{amsallem_design_2015}, microfluidics \cite{antil_reduced_2012}, and finance \cite{sachs_priori_2013}.

One challenge with using the POD approach in optimal control problems is that the basis must be precomputed based on a reference control which may differ significantly from the final optimal control. 
As a result, the suboptimal control derived from the POD model may not provide a good approximation of the full-order optimal control. 
To address this, several adaptive strategies have been developed \cite{Afanasiev_wake_2001, ravindran_adaptive_2002, alla_adaptive_2013} that update the POD basis during the optimization process. 
Some of these methods have established strong mathematical foundations over the years, such as the \textit{optimality system POD} (OSPOD) \cite{kunisch_proper_2008} and \textit{trust region POD} (TRPOD) \cite{arian_trust-region_nodate}. 
OSPOD addresses the issue of unmodeled dynamics by updating the POD basis in the direction that minimizes the cost functional, ensuring that the POD-reduced system is computed from the trajectory associated with the optimal control. 
Convergence results and a-posteriori error estimates for OSPOD have been studied in \cite{volkwein_optimality_nodate, kunisch_uniform_2015}. 
TRPOD, on the other hand, manages model approximation quality within a trust-region framework by comparing the predicted reduction with the actual reduction in the model function value. 
TRPOD has also seen notable applications in recent years \cite{bergmann_optimal_2008, yue_accelerating_2013}.

Until now, it has been assumed that POD approximations of the state are sufficiently accurate when using a small POD basis, and in many PDEs, particularly those with a faster decay of the Kolmogorov $n$-width, a low-dimensional subspace with low approximation error can indeed be effectively captured using POD. 
However, for problems where the Kolmogorov $n$-width decays slowly, the trial subspace constructed with POD may no longer remain low-dimensional.
This issue is especially prevalent when dealing with \textit{transport-dominated fluid systems} (TDFS), such as propagating flame fronts or traveling acoustic and shock waves \cite{krah_front_2022}. 
To address the challenge posed by slowly decaying singular values, new MOR methods for TDFS have been developed in recent years. 
One such method is the shifted POD (sPOD) \cite{reiss_shifted_2018, krahmarminzorawskireissschneider2024}, which forms a nonlinear projection framework to capture the high-dimensional space, followed by constructing the reduced-order approximation via Galerkin projection \cite{black_efficient_2021, black_projection-based_2020}. 
Other MOR techniques for TDFS include \cite{rim_transport_2018}, which uses transport reversal and template fitting; \cite{nonino_overcoming_2019}, which employs transport maps; and \cite{krah_front_2022}, which approximates the field variable using a front shape function and a level set function for efficient model reduction.

From here onwards, and to focus on the key aspects of our contribution, instead of the continuous system \eqref{eq_def:PDE_continuous}, we consider a semi-discrete model corresponding to a discretization of the spatial domain $\Omega$ into $m$ grid points. 
The approximated PDE in terms of the discretized state variable $q(t)\in \mathbb{R}^m$ and the control $u(t)\in \mathbb{R}^{n_c}$ now reads:
\begin{equation}\label{eq_def:FOM_discrete}
\begin{aligned} 
    &\dot{q}(t) =  Aq(t) + Bu(t), \\
    &q(0) = q_0, \\
\end{aligned}
\end{equation}
where $A\in \mathbb{R}^{m\times m}$ is the discrete periodic approximation of the linear advection operator $\mathcal{A}:=-v \frac{\partial }{\partial x}$ with $v$ being a constant advection velocity and $B\in \mathbb{R}^{m \times n_c}$ is the discrete approximation of the control operator $\mathcal{B}$. In particular, we assume the continuous control term $(\mathcal{B}(x)u)(t)$ to be given as $\sum^{n_c}_{k=1} \mathcal{B}_k(x) u_k(t)$ where $\mathcal{B}_k$ are control shape functions and $u_k$ are the control intensities. We then minimize the spatially discrete cost functional 
\begin{equation}\label{eq_def:FOM_costFunc_discrete}
    \mathcal{J}(q, u) = \frac{1}{2} \int^{t_f}_0  \norm{C (q(t) - q_\mathrm{d}(t))}^2_{\mathbb{R}^{p}}  + \mu   \norm{u(t)}^2_{ \mathbb{R}^{n_c}} \mathrm{d}t
\end{equation}
where $q_\mathrm{d}$ is the target state vector and $C=\mathrm{diag}\left(\frac{\sqrt{\mathrm{d}x}}{\sqrt{2}}, \sqrt{\mathrm{d}x}, \ldots, \sqrt{\mathrm{d}x}, \frac{\sqrt{\mathrm{d}x}}{\sqrt{2}}\right) \in \mathbb{R}^{m\times m}$ is the discrete weighting associated with the trapezoidal sum approximation of the $L^2(\Omega)$-norm when using a uniform grid spacing $\mathrm{d}x$. For the linear semi-discretized finite-dimensional optimal control problem, necessary optimality conditions are well-known and read
\begin{subnumcases}{\mathrm{OC}_{\scaleto{\mathrm{FOM}}{4pt}} := }
    \dot{q} =  Aq + Bu\label{eq_def:OC_FOM_1},  \\
    q(0) = q_0\label{eq_def:OC_FOM_2}, \\[1ex]
    -\dot{p} = A^\top\,p  + C^\top C(q - q_\mathrm{d})\label{eq_def:OC_FOM_3},  \\
    p(t_f) = 0\label{eq_def:OC_FOM_4},\\[1ex]
    \mu u + B^\top p = 0. \label{eq_def:OC_FOM_5}
\end{subnumcases} 
The equations \eqref{eq_def:OC_FOM_3}, \eqref{eq_def:OC_FOM_4} and \eqref{eq_def:OC_FOM_5} can also be viewed as the semi-discretized finite-dimensional approximation of \eqref{eq_def:PDE_adjoint_continuous}. 
At this point, it is important to note that the necessary optimality conditions still scale with the dimension $m$ of the \textit{full order model} (FOM). 
As discussed earlier, we can replace the FOM with ROMs to solve the optimal control problem. In doing so, we have the option of employing either the \textit{first optimize then reduce} (FOTR) framework or the \textit{first reduce then optimize} (FRTO) framework when using MOR methods for solving the problem. 
With this context in mind, the following are the contributions of this paper:
\begin{itemize} 
    \item We derive the necessary optimality conditions for the \textit{sPOD-Galerkin} (sPOD-G) method for both the FRTO and the FOTR frameworks which especially in the first case becomes non-trivial since the sPOD-G method itself generically results in a non-linear set of equations. 
    \item We comment on the commutativity of the two presented frameworks for the sPOD-G method both with respect to theoretical but also numerical aspects. While commutativity results for the POD-G method have already been investigated \cite{hinze_optimization_2009}, to the best of our knowledge, no such study is available for the sPOD-G method.
    \item We propose algorithmic simplifications and improvements to the frameworks with the sPOD-G method by exploiting specific structures in the PDE model to speed up the computation required for solving the optimal control problem.
    Although such simplifications are problem-specific, nevertheless, they serve as a first crucial step in making the sPOD-G method computationally comparable to other model reduction methods.
    \item We examine a 1D linear advection test case, comparing the results in terms of reduced order dimension, convergence behavior, and computational time for both the sPOD-G and POD-G methods. Here, the reduced order models are used as computationally cheap surrogates within a line search that generally requires multiple evaluations of the forward model for finding a descent direction.
\end{itemize}
In \Cref{sec:MOR_methods}, we begin by outlining the foundational concepts for the MOR methods used in this work, specifically the POD-G and sPOD-G methods. 
We briefly recall the (well-known) necessary optimality conditions for the POD-G method in \Cref{ssec:POD-G} and subsequently derive the sPOD-G in \Cref{ssec:sPOD-G} in detail.
In \Cref{sec:algorithms} we discuss the algorithmic details, along with specific simplifications and improvements.
Numerical results for the 1D linear advection equation, including a comprehensive comparison and timing analysis, are presented in \Cref{sec:results}. 
Finally, we summarize our findings in \Cref{sec:conclusion}.

\subsection*{Notation}\label{ssec:notations}
Matrices are denoted by upper case letters, vectors are denoted by lower case letters. 
The symbols $\langle \cdot, \cdot \rangle$, and $\norm{\cdot}$ denote the Euclidean inner product and norm, respectively. 
The space of real $m\times n$ matrices is denoted by $\mathbb{R}^{m\times n}$. 
Moreover, the abbreviation $\mathrm{blkdiag}$ refers to a block diagonal matrix constructed as:
\begin{equation*}
    \mathrm{blkdiag}(C_1, \ldots, C_n) :=
\begin{pmatrix}
C_{1} \\
& \ddots \\
& & C_{n}
\end{pmatrix}
\end{equation*}
where $C_1, \ldots, C_n$ are arbitrary matrices of arbitrary sizes. 
For brevity, we often omit the explicit $t$ dependencies of all vectors and matrices. 
For an arbitrary parameter-dependent matrix $C(p)$, its derivative (with respect to the parameter $p$) is denoted as $C'(p)$ or, in short, simply $C'$. For 3-tensors arising from derivatives of parameter-dependent matrices encountered in our work, we simply use the operation ($'$) on the matrix itself, most commonly as $C'$ or $(C^\top)'$.
If clear from the context, all 3-tensors are multiplied from the appropriate side such that the common dimension with a vector is contracted to form a matrix. Given an arbitrary tensor $C' \in \mathbb{R}^{r_1\times r_2\times r_3}$ and an arbitrary vector $a\in \mathbb{R}^{r_2}$, the result will be $C = C' a \in \mathbb{R}^{r_1\times r_3}$.
We also abuse the notation by using the same variables while describing the FOTR and FRTO framework. 
Since these frameworks are described distinctly and hold no similarity unless otherwise stated, the symbols and notations concerning them should be treated distinctly as well.
In the POD-G method an arbitrary quantity $c$ is subscripted as $c_\mathrm{p}$ for the state equation and $c_\mathrm{pa}$ for the adjoint equation. 
Similarly in the sPOD-G method, we use $c_\mathrm{s}$ for the state equation and $c_\mathrm{sa}$ for the adjoint equation.

\section{Model order reduction methods}\label{sec:MOR_methods}

\subsection{POD-Galerkin method}
The POD approximates the quantity $y(x, t)$ from \eqref{eq_def:PDE_continuous} as:
\begin{equation}
    y(x, t) = \sum^{\mathrm{p}}_{i=1} \phi_i(x)\alpha_i(t)
\end{equation}
where $\phi_i(x)$ are the spatial basis functions and $\alpha_i(t)$ are the time coefficients.
Subsequently, after discretizing we extract optimal basis vectors from a collection of snapshots. 
The snapshot matrix $Q = [q(t_1), \cdots, q(t_n)]\in \mathbb{R}^{m \times n}$ comprises of snapshots $q(t_j)$ computed from solving \eqref{eq_def:FOM_discrete} arranged in a column-wise fashion for each time step $t_j$ where $j=1, \ldots, n$. 
The single entries of each column correspond to approximations of the original PDE solution at the spatial grid points. 
The POD method approximates the snapshot matrix $Q$ with the help of the singular value decomposition (SVD) 
\begin{equation}\label{eq_def:POD}
      Q \approx Q_{\mathrm{POD}} = U_\mathrm{p} \Sigma_\mathrm{p} (V_\mathrm{p})^\top\, .
\end{equation}
Here, $\mathrm{p}\ll \min(n, m)$ is the truncation rank, $\Sigma_{\mathrm{p}}=\mathrm{diag}{(\sigma_1,\dots,\sigma_{\mathrm{p}})}$ is a diagonal matrix containing the leading $\mathrm{p}$ singular values $\sigma_1\ge\sigma_2\ge\dots\ge\sigma_{\mathrm{p}}$ of $Q$ and $U_{\mathrm{p}}\in\mathbb{R}^{m\times {\mathrm{p}}}$ and $V_{\mathrm{p}}\in\mathbb{R}^{n\times {\mathrm{p}}}$ satisfy $U_{\mathrm{p}}^\top U_{\mathrm{p}}=I_\mathrm{p}$ and $V_{\mathrm{p}}^\top V_{\mathrm{p}}=I_\mathrm{p}$. 
The POD yields the best rank-$\mathrm{p}$ approximation of the snapshot matrix $Q$ with respect to the Frobenius norm \cite{eckart_approximation_1936}. 
Once the basis is computed we then perform Galerkin projection on \eqref{eq_def:FOM_discrete} and obtain the dynamical reduced order model in terms of the POD amplitudes $a_{\mathrm{p}}$:
\begin{equation}\label{eq_def:POD_Galerkin}
\begin{aligned}
    & \dot{a}_{\mathrm{p}} = U_\mathrm{p}^\top A U_\mathrm{p} a_{\mathrm{p}} + U_\mathrm{p}^\top B u,  \\
    & a_{\mathrm{p}}(0) = U_\mathrm{p}^\top q_0.
\end{aligned}
\end{equation}
Once $a_{\mathrm{p}}$ is obtained by simulating the system \eqref{eq_def:POD_Galerkin} it is lifted with the basis vector $U_\mathrm{p}$ and the approximation obtained is given as
\begin{equation}
    q(t) \approx U_\mathrm{p} a_{\mathrm{p}}(t).
\end{equation}

\subsection{sPOD-Galerkin method}\label{ssec:sPODG}
The sPOD was introduced in \cite{reiss_shifted_2018} and further algorithmic developments were presented in \cite{reiss_optimization-based_2021, krah_non-linear_nodate,krahmarminzorawskireissschneider2024, black_projection-based_2020, black_efficient_2021}.
The sPOD method decomposes $y(x, t)$ using a non-linear decomposition ansatz:
\begin{equation}
    y(x, t) = \sum^{K}_{k=1}\mathcal{T}^k y^k(x, t) \quad\quad \text{with}\quad \quad y^k(x, t) = \sum^{\mathrm{p}^k}_{i=1} \phi^k_i(x) \alpha^k_i(t)
\end{equation}
into multiple co-moving quantities $y^k(x, t)$ where $\mathcal{T}^k$ are the continuous transport operators that transform the co-moving quantities into the reference frame. 
In our work we consider the transport operators to perform shift transformations as $\mathcal{T}^ky^k(x, t) = y^k(x - \Delta^k(t), t)$ where $\Delta^k(t)$ are the continuous time-dependent shifts.
Subsequently, after discretization we construct the snapshot matrix ${Q}\in\mathbb{R}^{m\times n}$ and decompose it into multiple co-moving fields $\{Q^k\in\mathbb{R}^{m\times n}\}_{k=1,\dots,K}$ as:
\begin{align}
    Q &\approx Q_\mathrm{sPOD} :=\sum_{k=1}^K T^{z^k} \matr{Q}^k\,.
    \label{eq_def:sPOD-decomposition-discrete}
\end{align}
Here, $T^{z^k}\colon\mathbb{R}^{m \times n}\rightarrow\mathbb{R}^{m \times n}$ are the discretized transport operators that implement a discrete time-dependent shift $z^k\in \mathbb{R}^n$ on $Q$:
\begin{equation*}
    T^{z^k} Q = \left[ \bar{T}^{z_1^k}q(t_1),\dots,\bar{T}^{z_n^k}q(t_{n})\right]\in \mathbb{R}^{m\times n},
\end{equation*}
where $\bar{T}^{z^k}\colon \mathbb{R}^m \rightarrow \mathbb{R}^m$ are the discrete transport operators acting on the individual columns of $Q$. 

\begin{rmrk}
    In our numerical experiments, $\bar{T}^{z^k}\colon \mathbb{R}^m \rightarrow \mathbb{R}^m$ are constructed by using Lagrange polynomials of order $5$, which introduce an error of the order $\mathcal{O}(h^6)$. This comparably high order of the polynomial interpolants is chosen as we observed small errors in the shifting process to impact the overall convergence of the optimal control loop.
\end{rmrk}

The idea behind sPOD is that for transport-dominated systems the ansatz \eqref{eq_def:sPOD-decomposition-discrete} can decompose the snapshot data more efficiently than POD as the transported quantity remains stationary in the co-moving data frame. For clarity, we refer to the following motivational example:
\begin{example}
\footnote{A similar example has also been shown in \cite{burela_parametric_2023}.} We consider the advection equation shown in \eqref{eq_def:OC_FOM_1} with initial condition ${q}(x,0) = \exp(-(x_i - L/20)^2 / 0.01)$ and solve it on a discrete grid $x_i=ih$ for $i=1, \dots, m$ where $(x, t) \in ]0, L]\times [0, T[$. The spatial domain length is $L=10 \: \mathrm{m}$ and the temporal domain is $T=10 \: \mathrm{s}$. We discretize the spatial and temporal domain into $m=2000$ and $n=2000$ equidistant points respectively. 
\Cref{fig:svd} shows the discrete traveling wave solution of the advection equation for the constant advection velocity $v=0.9 \: \mathrm{m} / \mathrm{s}$ along with the corresponding singular value decay of the solution. 
It also shows the shifted stationary wave profile and its singular value decay. 
We observe that the singular value decay of the traveling wave is significantly slower compared to that of the stationary wave. 
\begin{figure}[h]
\centering
\includegraphics[scale=0.5]{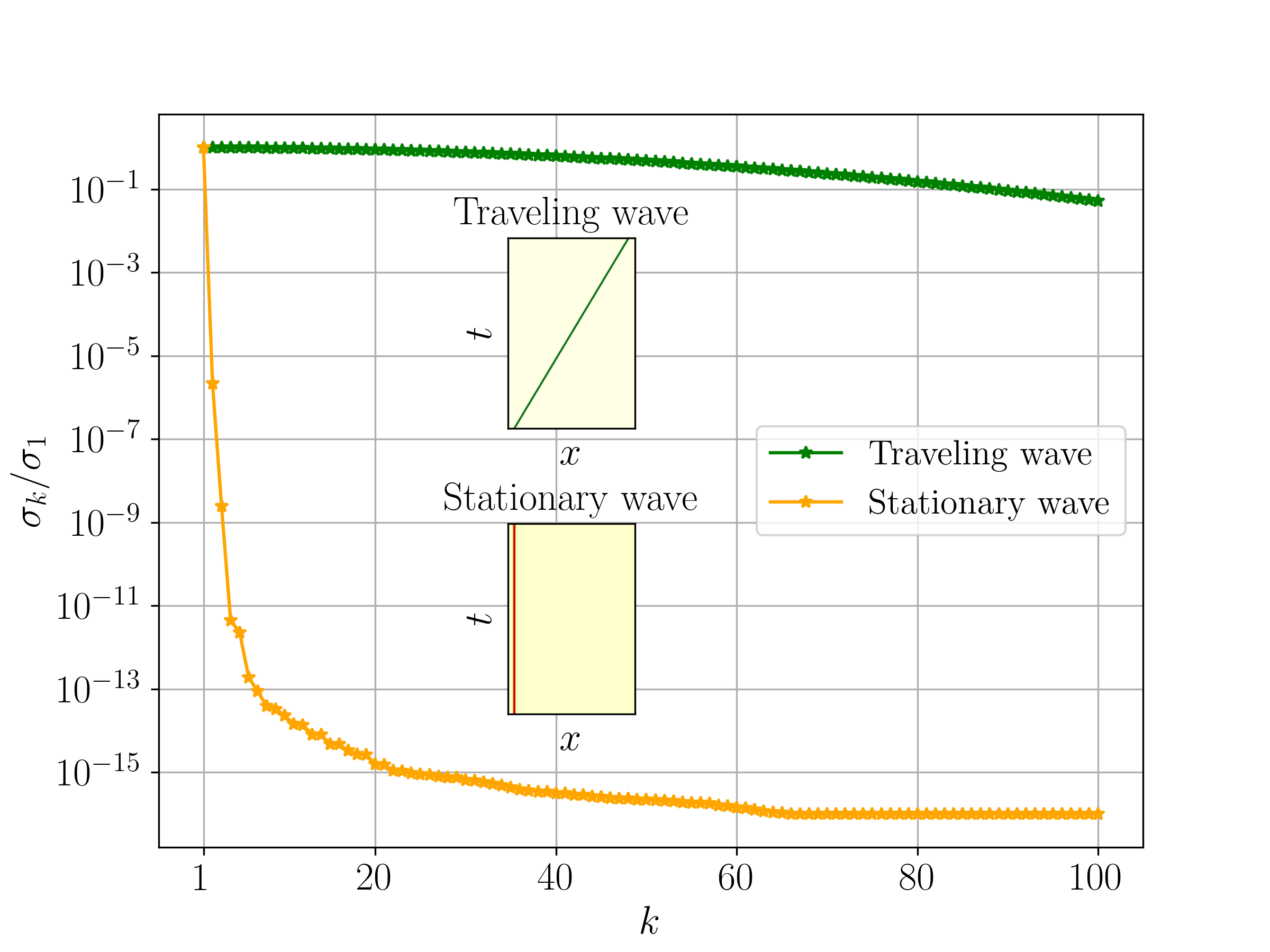}
\caption{Singular value decay of traveling and stationary wave}\label{fig:svd}
\end{figure}
\end{example}
The transported quantity in the co-moving data frame thus could be approximated efficiently with just a few modes with a truncated SVD: 
\begin{align}
         \matr{Q}^k \approx \matr{U}_{\mathrm{p}^k}^k\matr{\Sigma}_{\mathrm{p}^k}^k (\matr{V}_{\mathrm{p}^k}^k)^\top\,\quad k=1,\dots,K\,.
\label{eq-def:sPOD-Qk}
\end{align}
Here,  $\mathrm{p}^k\ll \min(m, n)$ is the truncation rank of each co-moving field, $\Sigma^k_{\mathrm{p}^k}=\mathrm{diag}{(\sigma_1^k,\dots,\sigma_{\mathrm{p}^k}^k)}$ is a diagonal matrix containing the singular values $\sigma_1^k\ge\sigma_2^k\ge\cdots \ge \sigma_{\mathrm{p}^k}^k$ and $\matr{U}_{\mathrm{p}^k}^k\in\mathbb{R}^{m\times {\mathrm{p}^k}}$, $\matr{V}_{\mathrm{p}^k}^k\in\mathbb{R}^{n\times {\mathrm{p}^k}}$ satisfy $(U^k_{\mathrm{p}^k})^\top U^k_{\mathrm{p}^k}=I_{\mathrm{p}^k}$ and $(V^k_{\mathrm{p}^k}) ^\top V^k_{\mathrm{p}^k}=I_{\mathrm{p}^k}$. 
Once the basis is constructed we obtain the reduced dynamical model by Galerkin projection as explained in \cite{black_projection-based_2020, black_efficient_2021} using the sPOD ansatz \eqref{eq_def:sPOD-decomposition-discrete}:
\begin{equation}\label{eq_def:sPOD_galerkin}
   \begin{aligned}
    \begin{bmatrix} M_{1, \mathrm{s}}(z_\mathrm{s}) & N_\mathrm{s}(z_\mathrm{s})D_\mathrm{s}(a_\mathrm{s}) \\ D_\mathrm{s}(a_\mathrm{s})^\top N_\mathrm{s}(z_\mathrm{s})^\top & D_\mathrm{s}(a_\mathrm{s})^\top M_{2, \mathrm{s}}(z_\mathrm{s}) D_\mathrm{s}(a_\mathrm{s})\end{bmatrix} \begin{bmatrix} \dot{a}_{\mathrm{s}} \\ \dot{z_\mathrm{s}} \end{bmatrix} &= \begin{bmatrix} A_{1, \mathrm{s}}(z_\mathrm{s}) & 0 \\ D_\mathrm{s}(a_\mathrm{s})^\top A_{2, \mathrm{s}}(z_\mathrm{s}) & 0\end{bmatrix}
    \begin{bmatrix} a_{\mathrm{s}} \\ z_\mathrm{s} \end{bmatrix}  + \begin{bmatrix}  V_\mathrm{s}(z_\mathrm{s})^\top B u \\ \mathcal{D}_\mathrm{s}(a_\mathrm{s})^\top W_\mathrm{s}(z_\mathrm{s})^\top B u \end{bmatrix}\\
    V_s(z_\mathrm{s}(0))^\top V_s(z_\mathrm{s}(0)) a_{\mathrm{s}}(0) = V_s(z_\mathrm{s}(0))^\top q_0, &\quad \quad z_\mathrm{s}(0) = z_0,
      \end{aligned}
\end{equation}
where the shift-dependent matrices $V_\mathrm{s}(z_\mathrm{s})\colon\mathbb{R}^r \rightarrow \mathbb{R}^m$ and $W_\mathrm{s}(z_\mathrm{s})\colon\mathbb{R}^r \rightarrow \mathbb{R}^m$ are given as:
\begin{equation}\label{eq_def:sPOD_galerkin_basis}
    \begin{aligned}
    V_\mathrm{s}(z_\mathrm{s}) &= \begin{bmatrix} \bar{T}^{z^1}U^1_1,\ldots, \bar{T}^{z^1}U^1_{\mathrm{p}^1}, & \ldots & \bar{T}^{z^K}U^K_1,\ldots, \bar{T}^{z^K}U^K_{\mathrm{p}^K} \end{bmatrix} \in \mathbb{R}^{m\times r} \\
    W_\mathrm{s}(z_\mathrm{s}) &= \begin{bmatrix} \bar{T}'^{z^1}U^1_1,\ldots, \bar{T}'^{z^1}U^1_{\mathrm{p}^1}, & \ldots & \bar{T}'^{z^K}U^K_1,\ldots, \bar{T}'^{z^K}U^K_{\mathrm{p}^K} \end{bmatrix} \in \mathbb{R}^{m\times r}
    \end{aligned} 
\end{equation}
with $r = \sum^{K}_{k=1}\mathrm{p}^k$ and the system matrices in terms of \eqref{eq_def:sPOD_galerkin_basis} are given as:
\begin{equation}\label{eq_def:sPOD_galerkin_sysmat}
    \begin{aligned}
     M_{1, \mathrm{s}}(z_\mathrm{s}) = V_\mathrm{s}(z_\mathrm{s})^\top V_\mathrm{s}(z_\mathrm{s})\in \mathbb{R}^{r\times r} & ,\quad  M_{2, \mathrm{s}}(z_\mathrm{s}) = W_\mathrm{s}(z_\mathrm{s})^\top W_\mathrm{s}(z_\mathrm{s})\in \mathbb{R}^{r\times r},\quad N_\mathrm{s}(z_\mathrm{s}) = V_\mathrm{s}(z_\mathrm{s})^\top W_\mathrm{s}(z_\mathrm{s}) \in \mathbb{R}^{r\times r},\\
    A_{1, \mathrm{s}}(z_\mathrm{s}) = V_\mathrm{s}(z_\mathrm{s})^\top A V_\mathrm{s}(z_\mathrm{s})\in \mathbb{R}^{r\times r} & ,\quad A_{2, \mathrm{s}}(z_\mathrm{s}) = W_\mathrm{s}(z_\mathrm{s})^\top A V_\mathrm{s}(z_\mathrm{s})\in \mathbb{R}^{r\times r}.
    \end{aligned}
\end{equation}
Note that for more than one shift, the matrix $V_s(z_s(0))$ is generally not expected to have orthonormal columns so the reduced initial value $a_s(0)$ is obtained by solving a linear system.
The terms $$a_{\mathrm{s}}(t) = \begin{bmatrix} a^1_1(t) \ldots a^1_{\mathrm{p}^1}(t) & a^2_1(t) \ldots a^2_{\mathrm{p}^2}(t) & \ldots & a^K_1(t) \ldots a^K_{\mathrm{p}^K}(t)\end{bmatrix}^\top\, \in \mathbb{R}^{r}, \ \  z_\mathrm{s}(t) = \begin{bmatrix} z^1(t) & z^2(t) & \ldots & z^K(t) \end{bmatrix}^\top \in \mathbb{R}^{K}$$ are the vectors of sPOD amplitudes corresponding to each co-moving field and the time-dependent shifts associated with them respectively. 
In addition, $D_\mathrm{s}(a_\mathrm{s}) = \mathrm{blkdiag}(a^1(t), \ldots, a^K(t)) \in \mathbb{R}^{r\times K}$. We observe here that \eqref{eq_def:sPOD_galerkin} is a non-linear equation due to the non-linear dependence of $z^k$ on $T^{z^k}$. Once \eqref{eq_def:sPOD_galerkin} is simulated the approximation of the FOM state is obtained as:
\begin{equation}
    q(t) \approx V_s(z_\mathrm{s}(t)) a_\mathrm{s}(t)\ .
\end{equation}

\begin{rmrk}
    For the remainder of the manuscript, we omit the explicit $z_\mathrm{s}$ and $a_\mathrm{s}$ dependency from the shift-dependent matrices $V_\mathrm{s}$, $W_\mathrm{s}$, $M_{1, \mathrm{s}}$, $M_{2, \mathrm{s}}$, $A_{1, \mathrm{s}}$, $A_{2, \mathrm{s}}$, $N_\mathrm{s}$ and the reduced state-dependent matrix $D_\mathrm{s}$ respectively for simplicity. 
\end{rmrk}

\section{Necessary optimality conditions for MOR methods}\label{sec:NOC_MOR_methods}
In this section, we discuss the first-order necessary optimality conditions for both the POD-G and sPOD-G methods. 
While the conditions for POD-G are well-known, we will mainly focus on sPOD-G and provide a detailed derivation for this method.

\subsection{Optimality conditions for the POD-G method}\label{ssec:POD-G}
Let us first look at the FRTO framework, which builds upon the reduced state equation \eqref{eq_def:POD_Galerkin}. The associated optimality system therefore reads
\begin{subnumcases}{\mathrm{OC}^{\scaleto{\mathrm{FRTO}}{4pt}}_{\scaleto{\mathrm{POD-G}}{4pt}}:=} 
    \dot{a}_\mathrm{p} = U_\mathrm{p}^\top A U_\mathrm{p} a_\mathrm{p} + U_\mathrm{p}^\top B u \label{eq_def:OC_FRTO_PODG_1},\\
    a_\mathrm{p}(0) = U_\mathrm{p}^\top q_0 \label{eq_def:OC_FRTO_PODG_2},\\[1em]
    -\dot{a}_\mathrm{pa} =  U^\top_\mathrm{p} A^\top U_\mathrm{p} a_\mathrm{pa} + U^\top_\mathrm{p} C^\top C(U_\mathrm{p} a_\mathrm{p} - q_{\mathrm{d}})  \label{eq_def:OC_FRTO_PODG_3}, \\
    a_\mathrm{pa}(t_f) = 0  \label{eq_def:OC_FRTO_PODG_4},\\[1em]
    \mu u + B^\top U_\mathrm{p} a_\mathrm{pa} = 0 \label{eq_def:OC_FRTO_PODG_5}
\end{subnumcases}
Here, the first two equations are the reduced-order state equations as mentioned before and the next two equations are the adjoint equations derived from the reduced-order state equation. 
The last equation \eqref{eq_def:OC_FRTO_PODG_5} is analogous to \eqref{eq_def:OC_FOM_5} and incorporates the POD projection $U_\mathrm{p}^\top B$ of the control matrix $B$. 
We also note that we now minimize an approximation of the cost functional \eqref{eq_def:FOM_costFunc_discrete} given as:
\begin{equation}\label{eq_def:PODG_costFunc_discrete}
    \mathcal{J}_{\scaleto{\mathrm{POD-G}}{4pt}}(a_\mathrm{p}, u) = \frac{1}{2} \int^{t_f}_0  \norm{C(U_\mathrm{p} a_\mathrm{p}(t) - q_\mathrm{d}(t))}^2_{ \mathbb{R}^{m}}  + \mu  \norm{u(t)}^2_{\mathbb{R}^{n_c}} \mathrm{d}t.
\end{equation}

For the FOTR framework, we build the reduced order approximations for the already mentioned necessary optimality conditions  $\mathrm{OC}_{\scaleto{\mathrm{FOM}}{4pt}}$. 
Performing Galerkin projection we get:
\begin{subnumcases}{\mathrm{OC}^{\scaleto{\mathrm{FOTR}}{4pt}}_{\scaleto{\mathrm{POD-G}}{4pt}}:=}
    \dot{a}_\mathrm{p} = U_{\mathrm{p}}^\top\, A U_{\mathrm{p}} a_\mathrm{p} + U_{\mathrm{p}}^\top\,B u  \label{eq_def:OC_FOTR_PODG_1},  \\
    a_\mathrm{p}(0) = U_{\mathrm{p}}^\top\, q_0,\label{eq_def:OC_FOTR_PODG_2}\\[1em]
    -\dot{a}_\mathrm{pa} = U_{\mathrm{pa}}^\top\, A^\top\,U_{\mathrm{pa}} a_\mathrm{pa} + U_{\mathrm{pa}}^\top C^\top C(U_\mathrm{p} a_\mathrm{p} - q_{\mathrm{d}})\label{eq_def:OC_FOTR_PODG_3},  \\
    a_\mathrm{pa}(t_f) = 0, \label{eq_def:OC_FOTR_PODG_4}\\[1em]
    \mu u + B^\top U_\mathrm{pa} a_\mathrm{pa} = 0 \label{eq_def:OC_FOTR_PODG_5}
\end{subnumcases}
and we minimize the cost functional \eqref{eq_def:PODG_costFunc_discrete}.
A key difference between $\mathrm{OC}^{\scaleto{\mathrm{FRTO}}{4pt}}_{\scaleto{\mathrm{POD-G}}{4pt}}$ and $\mathrm{OC}^{\scaleto{\mathrm{FOTR}}{4pt}}_{\scaleto{\mathrm{POD-G}}{4pt}}$ is the inclusion of a separate basis $U_\mathrm{pa}$ for the adjoint in the FOTR framework. 
It is known \cite{hinze_optimization_2009} that the two frameworks commute if $U_\mathrm{p} = U_\mathrm{pa}$ meaning that if a common basis is used for both the state and adjoint equations, the FOTR framework effectively becomes equivalent to the FRTO framework. 
This is typically achieved by concatenating the state and adjoint snapshots into a single matrix and performing an SVD, though Gram-Schmidt orthogonalization may be needed to adjust for scale differences between the state and adjoint solutions \cite{negri_reduced_2013}.

\subsection{Optimality conditions for the sPOD-G method}\label{ssec:sPOD-G}

Again, let us first focus on the sPOD-G method for which we consider the optimal control problem
\begin{equation}\label{eq_def:sPOD_costFunc_discrete}
    \underset{a_\mathrm{s}, z_\mathrm{s}, u}{\min} \:\:\mathcal{J}_{\scaleto{\mathrm{sPOD-G}}{4pt}}(a_\mathrm{s}, z_\mathrm{s}, u):= \frac{1}{2} \int^{t_f}_0  \norm{C(V_\mathrm{s} a_\mathrm{s} - q_\mathrm{d})}^2_{\mathbb{R}^{m}}  +  \mu \norm{u(t)}^2_{\mathbb{R}^{n_c}}\mathrm{d}t
\end{equation}
subject to the dynamical constraint \eqref{eq_def:sPOD_galerkin}. In the following, we often suppress the time dependency. For the Lagrangian, we obtain 
\begin{multline}\label{eq_def:sPOD_galerkin_linear_Lagr}
    \mathcal{L}(a_\mathrm{s}, z_\mathrm{s}, u, a_\mathrm{sa}, z_\mathrm{sa}) = \mathcal{J}_{\scaleto{\mathrm{sPOD-G}}{4pt}}(a_\mathrm{s}, z_\mathrm{s}, u) -  \int^{t_f}_0 \begin{bmatrix} a_\mathrm{sa} \\ z_\mathrm{sa} \end{bmatrix}^\top \left( \begin{bmatrix} M_{1, \mathrm{s}} & N_\mathrm{s} D_\mathrm{s} \\ D_\mathrm{s}^\top\, N_\mathrm{s}^\top\, & D_\mathrm{s}^\top\, M_{2, \mathrm{s}} D_\mathrm{s}\end{bmatrix} \begin{bmatrix} \dot{a}_\mathrm{s} \\ \dot{z}_\mathrm{s} \end{bmatrix}
    - \begin{bmatrix} A_{1, \mathrm{s}} & 0 \\ D_\mathrm{s}^\top\, A_{2, \mathrm{s}} & 0\end{bmatrix}
    \begin{bmatrix} a_\mathrm{s} \\ z_\mathrm{s} \end{bmatrix} \right.\\ \left.- \begin{bmatrix}  V_\mathrm{s}^\top\, B u \\ D_\mathrm{s}^\top\, W_\mathrm{s}^\top\, B u \end{bmatrix}\right)\mathrm{d}t
\end{multline}
where $a_\mathrm{sa}$ and $z_\mathrm{sa}$ are the adjoint variables. Applying integration by parts we obtain
\begin{multline}
    \mathcal{L}(a_\mathrm{s}, z_\mathrm{s}, u, a_\mathrm{sa}, z_\mathrm{sa}) = \frac{1}{2} \int^{t_f}_0  \norm{C(V_\mathrm{s} a_\mathrm{s} - q_\mathrm{d})}^2_{\mathbb{R}^{m}}  + \mu  \norm{u}^2_{\mathbb{R}^{n_c}} \mathrm{d}t\\ -\left[\langle a_\mathrm{sa},(M_{1, \mathrm{s}} a_\mathrm{s} + N_\mathrm{s} D_\mathrm{s} z_\mathrm{s})\rangle + \langle z_\mathrm{sa}, (D_\mathrm{s}^\top N_\mathrm{s}^\top a_\mathrm{s} + D_\mathrm{s}^\top M_{2, \mathrm{s}} D_\mathrm{s} z_\mathrm{s})\rangle \right]^{t_f}_{0} \\ + \int^{t_f}_{0} \left(\langle\dot{a}_\mathrm{sa}, M_{1, \mathrm{s}} a_\mathrm{s}\rangle + \langle a_\mathrm{sa} ,(M'_{1, \mathrm{s}}\dot{z}_\mathrm{s}) a_\mathrm{s}\rangle + \langle \dot{z}_\mathrm{sa}, D_\mathrm{s}^\top N_\mathrm{s}^\top a_\mathrm{s}\rangle + \langle z_\mathrm{sa}, ((D_\mathrm{s}^\top)'\dot{a}_\mathrm{s}) N_\mathrm{s}^\top a_\mathrm{s}\rangle + \langle z_\mathrm{sa},  D_\mathrm{s}^\top ((N_\mathrm{s}^\top)'\dot{z}_\mathrm{s}) a_\mathrm{s}\rangle \right.\\+ \left.\langle \dot{a}_\mathrm{sa}, N_\mathrm{s} D_\mathrm{s}z_\mathrm{s}\rangle +\langle a_\mathrm{sa},(N_\mathrm{s}'\dot{z}_\mathrm{s}) D_\mathrm{s}z_\mathrm{s}\rangle + \langle a_\mathrm{sa}, N_\mathrm{s}(D_\mathrm{s}'\dot{a}_\mathrm{s}) z_\mathrm{s}\rangle + \langle \dot{z}_\mathrm{sa}, D_\mathrm{s}^\top M_{2, \mathrm{s}} D_\mathrm{s}z_\mathrm{s} \rangle + \langle z_\mathrm{sa}, ((D_\mathrm{s}^\top)'\dot{a}_\mathrm{s}) M_{2, \mathrm{s}} D_\mathrm{s}z_\mathrm{s}\rangle \right.\\+ \left. \langle z_\mathrm{sa}, D_\mathrm{s}^\top (M_{2, \mathrm{s}}'\dot{z}_\mathrm{s}) D_\mathrm{s}z_\mathrm{s}\rangle + \langle z_\mathrm{sa}, D_\mathrm{s}^\top M_{2, \mathrm{s}} (D'_\mathrm{s}\dot{a}_\mathrm{s}) z_\mathrm{s}\rangle + \langle a_\mathrm{sa}, A_{1, \mathrm{s}} a_\mathrm{s}\rangle + \langle z_\mathrm{sa}, D_\mathrm{s}^\top A_{2, \mathrm{s}} a_\mathrm{s} \rangle + \langle a_\mathrm{sa}, V_\mathrm{s}^\top B u\rangle  \right. \\ + \left. \langle z_\mathrm{sa}, D_\mathrm{s}^\top W_\mathrm{s}^\top B u \rangle \right)\mathrm{d}t.
\end{multline}
In the above equation, we encounter 3-tensors which we contract as outlined in the notation section. For example, if the tensors $(M'_{1, \mathrm{s}}\dot{z}_\mathrm{s})$,    $M'_{1, \mathrm{s}} \in \mathbb{R}^{r\times r\times K}$ are contracted by $\dot{z}_\mathrm{s} \in \mathbb{R}^{K}$, the result is a matrix of dimension $r\times r$.
For obtaining the adjoint equation for $a_\mathrm{sa}$, we utilize that $\mathcal{D}_{a_\mathrm{s}}\mathcal{L} h = 0$ for any direction $h$ which results in\footnote{Here, for terms involving $\dot{h}$, we apply integration by parts a second time.}
\begin{multline}
    \int^{t_f}_0 \langle C(V_\mathrm{s}a_\mathrm{s} -  q_{\mathrm{d}}), CV_\mathrm{s} h\rangle - \langle \dot{a}_\mathrm{sa} ,M_{1, \mathrm{s}} h \rangle + \langle \dot{z}_\mathrm{sa}, D_\mathrm{s}^\top N_\mathrm{s}^\top h\rangle + \langle a_\mathrm{sa}, \left((M'_{1, \mathrm{s}}\dot{z}_\mathrm{s} h) - N_\mathrm{s} (D_\mathrm{s}'h\dot{z}_\mathrm{s}) + A_{1, \mathrm{s}} h\right)\rangle \\ + \langle z_\mathrm{sa}, \left(((D_\mathrm{s}^\top)'\dot{a}_\mathrm{s}) N_\mathrm{s}^\top h + D_\mathrm{s}^\top ((N_\mathrm{s}^\top)'\dot{z}_\mathrm{s}) h - (D_\mathrm{s}^\top)'h (N_\mathrm{s}^\top\dot{a}_\mathrm{s}) - (D_\mathrm{s}^\top)'h (M_{2, \mathrm{s}} D_\mathrm{s}\dot{z}_\mathrm{s}) - D_\mathrm{s}^\top M_{2, \mathrm{s}} (D_\mathrm{s}' h \dot{z}_\mathrm{s}) + (D_\mathrm{s}^\top)' h (A_{2, \mathrm{s}} a_\mathrm{s}) \right. \\ \left.+ D_\mathrm{s}^\top A_{2, \mathrm{s}} h + (D_\mathrm{s}^\top)'h (W_\mathrm{s}^\top B u)\right)\rangle \mathrm{d}t - \left[\langle a_\mathrm{sa}, (M_{1, \mathrm{s}} h + N_\mathrm{s} (D_\mathrm{s}'h z_\mathrm{s}))\rangle +  \langle z_\mathrm{sa},((D_\mathrm{s}^\top\,)'h (N_\mathrm{s}^\top\, a_\mathrm{s}) + D_\mathrm{s}^\top\, N_\mathrm{s}^\top\, h \right.\\ + \left. (D_\mathrm{s}^\top\,)'h (M_{2, \mathrm{s}} D_\mathrm{s} z_\mathrm{s}) + D_\mathrm{s}^\top\, M_{2, \mathrm{s}} (D_\mathrm{s}'h z_\mathrm{s}))\rangle \right]^{t_f}_{0} + \left[\langle a_\mathrm{sa}, N_\mathrm{s} (D_\mathrm{s}'hz_\mathrm{s}) \rangle\right]^{t_f}_{0} + \left[\langle z_\mathrm{sa}, (D_\mathrm{s}^\top)'h (N_\mathrm{s}^\top a_\mathrm{s}) \rangle \right]^{t_f}_{0} \\ + \left[\langle z_\mathrm{sa},(D_\mathrm{s}^\top)'h (M_{2, \mathrm{s}} D_\mathrm{s} z_\mathrm{s}) \rangle \right]^{t_f}_{0} + \left[\langle z_\mathrm{sa}, D_\mathrm{s}^\top M_{2, \mathrm{s}} (D_\mathrm{s}' hz_\mathrm{s}) \rangle \right]^{t_f}_{0}  = 0.
\end{multline}
Consequently, the first adjoint equation is given by
\begin{multline}\label{eq_def:sPOD_galerkin_linear_adjoint_1}
    M_{1, \mathrm{s}}^\top \dot{a}_\mathrm{sa}  + N_\mathrm{s} D_\mathrm{s} \dot{z}_\mathrm{sa} + \left(M'_{1, \mathrm{s}}\dot{z}_\mathrm{s} - N_\mathrm{s} (D_\mathrm{s}'\dot{z}_\mathrm{s}) + A_{1, \mathrm{s}}\right)^\top a_\mathrm{sa} + \left(((D_\mathrm{s}^\top)'\dot{a}_\mathrm{s}) N_\mathrm{s}^\top  + D_\mathrm{s}^\top ((N_\mathrm{s}^\top)'\dot{z}_\mathrm{s}) - (D_\mathrm{s}^\top)' (N_\mathrm{s}^\top\dot{a}_\mathrm{s}) \right.\\ - \left. (D_\mathrm{s}^\top)' (M_{2, \mathrm{s}} D_\mathrm{s}\dot{z}_\mathrm{s}) - D_\mathrm{s}^\top M_{2, \mathrm{s}} (D_\mathrm{s}'  \dot{z}_\mathrm{s}) + (D_\mathrm{s}^\top)'  (A_{2, \mathrm{s}} a_\mathrm{s}) + D_\mathrm{s}^\top A_{2, \mathrm{s}}  + (D_\mathrm{s}^\top)' (W_\mathrm{s}^\top B u)\right)^\top z_\mathrm{sa} \\ +  V_\mathrm{s}^\top C^\top C(V_\mathrm{s}a_\mathrm{s} - q_{\mathrm{d}}) = 0\,.
\end{multline}
If $h$ additionally satisfies $h(0) = 0$, we further have that
\begin{equation}\label{eq_def:sPOD_galerkin_linear_adjoint_boundcond_1}
     M_{1, \mathrm{s}}(t_f)^\top a_\mathrm{sa}(t_f)  = - N_\mathrm{s}(t_f) D_\mathrm{s}(t_f) z_\mathrm{sa}(t_f). 
\end{equation}
Similarly, with $\mathcal{D}_{z_\mathrm{s}}\mathcal{L} h = 0$, we obtain 
\begin{multline}
    \int^{t_f}_{0} \langle C(V_\mathrm{s}a_\mathrm{s} - q_{\mathrm{d}}), C (V_\mathrm{s}'a_\mathrm{s}) h \rangle + \langle \dot{a}_\mathrm{sa}, N_\mathrm{s} D_\mathrm{s}h \rangle + \langle \dot{z}_\mathrm{sa}, D_\mathrm{s}^\top M_{2, \mathrm{s}} D_\mathrm{s} h \rangle + \langle a_\mathrm{sa}, \Bigl((N'_\mathrm{s}\dot{z}_\mathrm{s}) D_\mathrm{s}h + N_\mathrm{s} (D'_\mathrm{s}\dot{a}_\mathrm{s}) h -  (M_{1, \mathrm{s}}' h\dot{a}_\mathrm{s}) \\ - N_\mathrm{s}'h (D_\mathrm{s}\dot{z}_\mathrm{s}) + (A_{1, \mathrm{s}}' h a_\mathrm{s}) + (V_\mathrm{s}^\top)' h (Bu)\Bigr)\rangle + \langle z_\mathrm{sa}, \Bigl(((D_\mathrm{s}^\top)'\dot{a}_\mathrm{s}) M_{2, \mathrm{s}} D_\mathrm{s} h + D_\mathrm{s}^\top (M'_{2, \mathrm{s}}\dot{z}_\mathrm{s}) D_\mathrm{s}h + D_\mathrm{s}^\top M_{2, \mathrm{s}} (D'_\mathrm{s}\dot{a}_\mathrm{s}) h - D_\mathrm{s}^\top ((N_\mathrm{s}^\top)'h \dot{a}_\mathrm{s}) \\ - D_\mathrm{s}^\top (M_{2, \mathrm{s}}' h (D_\mathrm{s}\dot{z}_\mathrm{s})) + D_\mathrm{s}^\top (A_{2, \mathrm{s}}' h a_\mathrm{s}) + D_\mathrm{s}^\top ((W_\mathrm{s}^\top)'h (B u))\Bigr)\rangle \mathrm{d}t - \Bigl[\langle a_\mathrm{sa}, \Bigl( M'_{1, \mathrm{s}} h a_\mathrm{s} + N_\mathrm{s}'h (D_\mathrm{s} z_\mathrm{s}) + N_\mathrm{s} D_\mathrm{s} h\Bigr)\rangle + \\ \langle z_\mathrm{sa}, \Bigl(D_\mathrm{s}^\top((N_\mathrm{s}^\top)'h a_\mathrm{s}) + D_\mathrm{s}^\top (M'_{2, \mathrm{s}} h (D_\mathrm{s} z_\mathrm{s})) + D_\mathrm{s}^\top M_{2, \mathrm{s}} D_\mathrm{s} h\Bigr)\rangle \Bigr]^{t_f}_{0}  + \Bigl[\langle a_\mathrm{sa}, (M_{1, \mathrm{s}}' h a_\mathrm{s})\rangle\Bigr]^{t_f}_{0} + \Bigl[\langle a_\mathrm{sa}, (N_\mathrm{s}' h (D_\mathrm{s}z_\mathrm{s}))\rangle\Bigr]^{t_f}_{0} \\ + \Bigl[\langle z_\mathrm{sa}, D_\mathrm{s}^\top ((N_\mathrm{s}^\top)' h a_\mathrm{s})\rangle\Bigr]^{t_f}_{0} + \Bigl[\langle z_\mathrm{sa}, D_\mathrm{s}^\top (M_{2, \mathrm{s}}' h (D_\mathrm{s}z_\mathrm{s}))\rangle\Bigr]^{t_f}_{0} = 0.
\end{multline}
The second adjoint equation therefore reads
\begin{multline}\label{eq_def:sPOD_galerkin_linear_adjoint_2}
    D_\mathrm{s}^\top N_\mathrm{s}^\top \dot{a}_\mathrm{sa} + D_\mathrm{s}^\top M_{2, \mathrm{s}}^\top D_\mathrm{s} \dot{z}_\mathrm{sa} + \Bigl((N'_\mathrm{s}\dot{a}_\mathrm{s}) D_\mathrm{s} + N_\mathrm{s} (D'_\mathrm{s}\dot{a}_\mathrm{s})  -  M_{1, \mathrm{s}}' \dot{a}_\mathrm{s} - N_\mathrm{s}' (D_\mathrm{s}\dot{z}_\mathrm{s}) + A_{1, \mathrm{s}}'  a_\mathrm{s} + (V_\mathrm{s}^\top)'  (B u)\Bigr)^\top a_\mathrm{sa} \\ + \Bigl(((D_\mathrm{s}^\top)'\mathrm{a}_\mathrm{s}) M_{2, \mathrm{s}} D_\mathrm{s}  + D_\mathrm{s}^\top  (M'_{2, \mathrm{s}}\dot{z}_\mathrm{s}) D_\mathrm{s} + D_\mathrm{s}^\top M_{2, \mathrm{s}} (D'_\mathrm{s}\dot{a}_\mathrm{s}) - D_\mathrm{s}^\top ((N_\mathrm{s}^\top)' \dot{a}_\mathrm{s}) - D_\mathrm{s}^\top (M_{2, \mathrm{s}}'  (D_\mathrm{s}\dot{z}_\mathrm{s})) + D_\mathrm{s}^\top (A_{2, \mathrm{s}}'  a_\mathrm{s}) \\ + D_\mathrm{s}^\top ((W_\mathrm{s}^\top)' (B u))\Bigr)^\top z_\mathrm{sa} + (C (V_\mathrm{s}' a_\mathrm{s}))^\top C( V_\mathrm{s}a_\mathrm{s} - q_{\mathrm{d}}) = 0\,.
\end{multline}
Again simplifying the boundary terms results in
\begin{equation}\label{eq_def:sPOD_galerkin_linear_adjoint_boundcond_2}
    D_\mathrm{s}(t_f)^\top N_\mathrm{s}(t_f)^\top a_\mathrm{sa}(t_f) = - D_\mathrm{s}(t_f)^\top M_{2, \mathrm{s}}(t_f)^\top D_\mathrm{s}(t_f) z_\mathrm{sa}(t_f).
\end{equation}
Combining  \eqref{eq_def:sPOD_galerkin_linear_adjoint_1} and \eqref{eq_def:sPOD_galerkin_linear_adjoint_2}, we obtain the following system of adjoint equations
\begin{equation}\label{eq_def:sPOD_galerkin_linear_adjoint}
    \begin{bmatrix} M_{1, \mathrm{s}}^\top & N_{\mathrm{s}} D_{\mathrm{s}} \\ D_{\mathrm{s}}^\top N_{\mathrm{s}}^\top & D_{\mathrm{s}}^\top M_{2, \mathrm{s}}^\top D_{\mathrm{s}}\end{bmatrix} \begin{bmatrix} \dot{a}_{\mathrm{sa}} \\ \dot{z}_{\mathrm{sa}} \end{bmatrix} = - \begin{bmatrix} E^\top_{11} & E^\top_{12} \\ E^\top_{21} & E^\top_{22} \end{bmatrix}
    \begin{bmatrix} a_{\mathrm{sa}} \\ z_{\mathrm{sa}} \end{bmatrix} - \begin{bmatrix}  V_\mathrm{s}^\top C^\top C(V_\mathrm{s}a_\mathrm{s} - q_{\mathrm{d}}) \\ (C (V_\mathrm{s}' a_\mathrm{s}))^\top C( V_\mathrm{s}a_\mathrm{s} - q_{\mathrm{d}})\end{bmatrix}
\end{equation}
where
\begin{align*}
&   E_{11} = M'_{1, \mathrm{s}}\dot{z}_\mathrm{s} - N_\mathrm{s} (D_\mathrm{s}'\dot{z}_\mathrm{s}) + A_{1, \mathrm{s}}      \\[1em]
&   \begin{multlined}[0.7\linewidth]
    E_{12} = ((D_\mathrm{s}^\top)'\dot{a}_\mathrm{s}) N_\mathrm{s}^\top  + D_\mathrm{s}^\top ((N_\mathrm{s}^\top)' \dot{z}_\mathrm{s}) - ( D_\mathrm{s}^\top)' (N_\mathrm{s}^\top\dot{a}_\mathrm{s}) - (D_\mathrm{s}^\top)' (M_{2, \mathrm{s}} D_\mathrm{s}\dot{z_\mathrm{s}}) - D_\mathrm{s}^\top M_{2, \mathrm{s}} (D_\mathrm{s}'  z_\mathrm{s}) \\ + (D_\mathrm{s}^\top)'  (A_{2, \mathrm{s}} a_\mathrm{s}) + D_\mathrm{s}^\top A_{2, \mathrm{s}}  + (D_\mathrm{s}^\top)' (W_\mathrm{s}^\top B u)
    \end{multlined}        \\[1em]
&   E_{21} = (N'_\mathrm{s} \dot{z}_\mathrm{s})D_\mathrm{s} + N_\mathrm{s} (D'_\mathrm{s}\dot{a}_\mathrm{s})  -  M_{1, \mathrm{s}}' \dot{a}_\mathrm{s} - N_\mathrm{s}'(D_\mathrm{s}\dot{z_\mathrm{s}}) + A_{1, \mathrm{s}}'  a_\mathrm{s} + (V_\mathrm{s}^\top)'  (B u)   \\[1em]
&    \begin{multlined}[0.7\linewidth]
        E_{22} = ((D_\mathrm{s}^\top)'\dot{a}_\mathrm{s}) M_{2, \mathrm{s}} D_\mathrm{s}  + D_\mathrm{s}^\top (M'_{2, \mathrm{s}}\dot{z}_\mathrm{s}) D_\mathrm{s} + D_\mathrm{s}^\top M_{2, \mathrm{s}} (D'_\mathrm{s}\dot{a}_\mathrm{s}) - D_\mathrm{s}^\top ((N_\mathrm{s}^\top)' \dot{a}_\mathrm{s}) - D_\mathrm{s}^\top (M_{2, \mathrm{s}}'  (D_\mathrm{s}\dot{z_\mathrm{s}})) \\ + D_\mathrm{s}^\top (A_{2, \mathrm{s}}'  a_\mathrm{s}) + D_\mathrm{s}^\top ((W_\mathrm{s}^\top)' (B u))
    \end{multlined}. 
\end{align*}
Consolidating the two terminal conditions \eqref{eq_def:sPOD_galerkin_linear_adjoint_boundcond_1} and \eqref{eq_def:sPOD_galerkin_linear_adjoint_boundcond_2}, it follows that
\begin{equation}
    \begin{bmatrix} M_{1, \mathrm{s}}^\top (t_f) & N_{\mathrm{s}}(t_f) D_{\mathrm{s}}(t_f) \\ D_{\mathrm{s}}^\top(t_f) N_{\mathrm{s}}^\top(t_f) & D_{\mathrm{s}}^\top(t_f) M_{2, \mathrm{s}}^\top(t_f) D_{\mathrm{s}}(t_f)\end{bmatrix} \begin{bmatrix} a_{\mathrm{sa}}(t_f) \\ z_{\mathrm{sa}}(t_f) \end{bmatrix} =
    \begin{bmatrix} 0 \\ 0 \end{bmatrix}
\end{equation}
determine the terminal condition for the adjoint variables. Finally, for the gradient of the Lagrangian with respect to the control $u$ we have 
\begin{equation}\label{eq_def:sPOD_galerkin_linear_dldu}
    \mathcal{D}_{u}\mathcal{L} h = \int^{t_f}_0 (\mu u + B^\top V_\mathrm{s} a_\mathrm{sa} + B^\top W_\mathrm{s} D_\mathrm{s} z_\mathrm{sa}) h \, \mathrm{d}t.
\end{equation}
The necessary optimality conditions are thus given as
  \begin{subnumcases}{\mathrm{OC}^{\scaleto{\mathrm{FRTO}}{4pt}}_{\scaleto{\mathrm{sPOD-G}}{4pt}}:=} 
    \begin{bmatrix} M_{1, \mathrm{s}} & N_{\mathrm{s}} D_{\mathrm{s}} \\ D_{\mathrm{s}}^\top N_{\mathrm{s}}^\top & D_{\mathrm{s}}^\top M_{2, \mathrm{s}} D_{\mathrm{s}}\end{bmatrix} \begin{bmatrix} \dot{a}_\mathrm{s} \\ \dot{z}_\mathrm{s} \end{bmatrix} = \begin{bmatrix} A_{1, \mathrm{s}} & 0 \\ D_{\mathrm{s}}^\top A_{2, \mathrm{s}} & 0\end{bmatrix}
    \begin{bmatrix} a_\mathrm{s} \\ z_\mathrm{s} \end{bmatrix} + \begin{bmatrix}  V_{\mathrm{s}}^\top B u \\ D_{\mathrm{s}}^\top W_{\mathrm{s}}^\top B u \end{bmatrix}, \label{eq_def:OC_FRTO_sPODG_1} \\
    a_{\mathrm{s}}(0) = V_{\mathrm{s}}(z_\mathrm{s}(0))^\top q_0, \quad \quad z_\mathrm{s}(0) = z_0 ,\label{eq_def:OC_FRTO_sPODG_2}\\[1em]
    \begin{bmatrix} M_{1, \mathrm{s}}^\top & N_{\mathrm{s}} D_{\mathrm{s}} \\ D_{\mathrm{s}}^\top N_{\mathrm{s}}^\top & D_{\mathrm{s}}^\top M_{2, \mathrm{s}}^\top D_{\mathrm{s}}\end{bmatrix} \begin{bmatrix} \dot{a}_{\mathrm{sa}} \\ \dot{z}_{\mathrm{sa}} \end{bmatrix} = - \begin{bmatrix} E^\top_{11} & E^\top_{12} \\ E^\top_{21} & E^\top_{22} \end{bmatrix}
    \begin{bmatrix} a_{\mathrm{sa}} \\ z_{\mathrm{sa}} \end{bmatrix} - \begin{bmatrix}  V_\mathrm{s}^\top C^\top C(V_\mathrm{s}a_\mathrm{s} - q_{\mathrm{d}}) \\ (C (V_\mathrm{s}' a_\mathrm{s}))^\top C( V_\mathrm{s}a_\mathrm{s} - q_{\mathrm{d}})\end{bmatrix}, \label{eq_def:OC_FRTO_sPODG_3} \\[1em]
    \begin{bmatrix} M_{1, \mathrm{s}}^\top (t_f) & N_{\mathrm{s}}(t_f) D_{\mathrm{s}}(t_f) \\ D_{\mathrm{s}}^\top(t_f) N_{\mathrm{s}}^\top(t_f) & D_{\mathrm{s}}^\top(t_f) M_{2, \mathrm{s}}^\top(t_f) D_{\mathrm{s}}(t_f)\end{bmatrix} \begin{bmatrix} a_{\mathrm{sa}}(t_f) \\ z_{\mathrm{sa}}(t_f) \end{bmatrix} =
    \begin{bmatrix} 0 \\ 0 \end{bmatrix} \label{eq_def:OC_FRTO_sPODG_4},\\[1.5em]
    \mu u + B^\top V_\mathrm{s}a_{\mathrm{sa}} + B^\top W_\mathrm{s} D_\mathrm{s} z_{\mathrm{sa}} = 0. \label{eq_def:OC_FRTO_sPODG_5}
 \end{subnumcases}
\begin{rmrk}
    We observe from the necessary optimality conditions that the adjoint equation contains 3-tensors like $D_\mathrm{s}'$, $(D_\mathrm{s}^\top\,)'$, $V'_\mathrm{s}$, $(V_\mathrm{s}^\top)'$, $N'_\mathrm{s}$, $(N^\top_\mathrm{s})'$, $M'_{1, \mathrm{s}}$, $M'_{2, \mathrm{s}}$, $A'_{1, \mathrm{s}}$, $A'_{2, \mathrm{s}}$, $(W_\mathrm{s}^\top)'$.
    All of these tensors, except  for $V'_\mathrm{s} \in \mathbb{R}^{m\times r\times K}$, $(V_\mathrm{s}^\top)' \in \mathbb{R}^{r\times m\times K}$ and $(W_\mathrm{s}^\top)' \in \mathbb{R}^{r\times m\times K}$ scale with the reduced order dimension $r$.
\end{rmrk}
Now for the FOTR framework, the derivation is rather straightforward and we follow the steps shown in \cite{black_projection-based_2020}. In particular, we only need to reduce the individual equations of $\mathrm{OC}_{\scaleto{\mathrm{FOM}}{4pt}}$ so that the necessary optimality conditions are given as
  \begin{subnumcases}{\mathrm{OC}^{\scaleto{\mathrm{FOTR}}{4pt}}_{\scaleto{\mathrm{sPOD-G}}{4pt}}:=}
    \begin{bmatrix} M_{1, \mathrm{s}} & N_\mathrm{s}D_\mathrm{s} \\ D_\mathrm{s}^\top N_\mathrm{s}^\top & D_\mathrm{s}^\top M_{2, \mathrm{s}} D_\mathrm{s}\end{bmatrix} \begin{bmatrix} \dot{a}_\mathrm{s} \\ \dot{z}_\mathrm{s} \end{bmatrix} = \begin{bmatrix} A_{1, \mathrm{s}} & 0 \\ D_\mathrm{s}^\top A_{2, \mathrm{s}} & 0\end{bmatrix}
    \begin{bmatrix} a_\mathrm{s} \\ z_\mathrm{s} \end{bmatrix} + \begin{bmatrix}  V_\mathrm{s}^\top B u \\ D_{\mathrm{s}}^\top W_{\mathrm{s}}^\top B u \end{bmatrix}\label{eq_def:OC_FOTR_sPODG_1},  \\
    a_{\mathrm{s}}(0) = V_\mathrm{s}(z_{\mathrm{s}}(0))^\top q_0, \quad \quad z_{\mathrm{s}}(0) = z_0 \label{eq_def:OC_FOTR_sPODG_2},\\[1em]
    \begin{bmatrix} M_{1, \mathrm{sa}} & N_\mathrm{sa} D_\mathrm{sa} \\ D_\mathrm{sa}^\top N_\mathrm{sa}^\top & D_\mathrm{sa}^\top M_{2, \mathrm{sa}} D_\mathrm{sa}\end{bmatrix} \begin{bmatrix} \dot{a}_\mathrm{sa} \\ \dot{z}_\mathrm{sa} \end{bmatrix} = - \begin{bmatrix} A_{1, \mathrm{sa}} & 0 \\ D_\mathrm{sa}^\top A_{2, \mathrm{sa}} & 0\end{bmatrix}\begin{bmatrix} a_\mathrm{\mathrm{sa}} \\ z_\mathrm{sa} \end{bmatrix} - \begin{bmatrix}  V_\mathrm{sa}^\top C^\top C(V_\mathrm{s} a_\mathrm{s} - q_\mathrm{d}) \\ D_\mathrm{sa}^\top W_\mathrm{sa}^\top C^\top C(V_\mathrm{s} a_\mathrm{s} - q_\mathrm{d}) \end{bmatrix} \label{eq_def:OC_FOTR_sPODG_3},  \\
    a_{\mathrm{sa}}(t_f) = V_\mathrm{sa}(z_\mathrm{sa}(t_f))^\top p(t_f)\Rightarrow a_\mathrm{sa}(t_f) = 0, \label{eq_def:OC_FOTR_sPODG_4}\\[1em]
    \mu u + B^\top V_\mathrm{sa} a_\mathrm{sa} = 0. \label{eq_def:OC_FOTR_sPODG_5}
 \end{subnumcases}
As the state equations \eqref{eq_def:OC_FOTR_sPODG_1} and \eqref{eq_def:OC_FOTR_sPODG_2} are identical, in both cases the cost functional to be minimized is \eqref{eq_def:sPOD_costFunc_discrete}. 
All the shift-dependent matrices and the system matrices are constructed as given in \eqref{eq_def:sPOD_galerkin_basis} and \eqref{eq_def:sPOD_galerkin_sysmat} respectively. 
However, we have a minor difference while assembling $A_{1, \mathrm{sa}}$ and $A_{2, \mathrm{sa}}$ for FOTR where we consider $A^\top$ instead of $A$ following the adjoint equation \eqref{eq_def:OC_FOM_3}. 
\begin{rmrk}\label{rem:FOTR_sPOD_adjoint_terminal}
    We note that from \eqref{eq_def:OC_FOM_4} we have $p(t_f)=0$ which yields $a_\mathrm{sa}(t_f)=0$ in \eqref{eq_def:OC_FOTR_sPODG_4}. This implies that we are free to choose $V_\mathrm{sa}(z_\mathrm{sa}(t_f))$ and in turn $z_\mathrm{sa}(t_f)$. Moreover, we conclude that the adjoint shifts $z_\mathrm{sa}$ at the terminal time point $t_f$ do not influence the adjoint solution. However, if we set $a_\mathrm{sa}(t_f)=0$, it follows $D_\mathrm{sa}(t_f)=0$ which in turn makes the mass matrix on the left-hand side of the adjoint equation \eqref{eq_def:OC_FOTR_sPODG_3} singular. 
    One would then need to regularize the mass matrix in order avoid dealing with a constrained differential (algebraic) equation. 
\end{rmrk}

\begin{rmrk}    
    One key observation from the derivations for the sPOD-G method is that the FOTR and FRTO frameworks do not appear to commute. 
    A major difference between them is that the FRTO framework results in a linear adjoint equation, whereas the FOTR framework leads to a non-linear approximation of the adjoint equation. 
    Consequently, the necessary optimality conditions for each framework are different. 
    Although the state equation remains the same for both frameworks, notable differences include: 
    (1) the adjoint equations, where the FRTO framework involves 3-tensors and includes the control term, while the FOTR framework does not; (2) different terminal conditions for the adjoint variables; and 
    (3) differing relationships between the control and adjoint variables.  
\end{rmrk}

\section{Backtracking line search based on reduced order surrogates}\label{sec:algorithms}

In view of the above discussion, one may use the sPOD-G-based optimality systems \eqref{eq_def:OC_FRTO_sPODG_1}-\eqref{eq_def:OC_FRTO_sPODG_5} and \eqref{eq_def:OC_FOTR_sPODG_1}-\eqref{eq_def:OC_FOTR_sPODG_5} to compute an optimal control for the full order problem. 
However, as the reduced order models are based on snapshots of solutions for a specific control $u$, one can generally not expect them to be accurate if the control changes, as is naturally the case in an optimal control approach. 
Already in the POD-G case, this is a non-trivial issue that has caused many authors to study suitable adaptation strategies for POD which generate basis updates whenever current reduced order models are no longer to be trusted. 
Among many of those techniques, we refer specifically to the OSPOD \cite{kunisch_proper_2008}, TRPOD \cite{fahl_reduced_2003}, and the works by Gr{\"a}\ss le et al. \cite{grasle_pod_2017, grasle_pod_2019}. 
While we made the first steps into generalizations of the previous methods to the sPOD-G framework, the additional challenges we faced do not yet allow for a full replacement of the optimal control problem by a reduced-order one. 
For FRTO the sPOD-G reduced adjoint equation implementation becomes tedious and complex and the equation itself scales with the FOM dimension because of the presence of $q_\mathrm{d}$ thus providing minimal computational benefits.
This FOM scaling is also true for FOTR where additionally, solving the sPOD-G reduced adjoint equation requires solving the FOM adjoint equation for basis generation. 
For transport-dominated systems, this happens frequently and we thus end up solving the FOM adjoint equation almost always.
Moreover, solving the reduced-order adjoint equation for either FRTO or FOTR with user-defined truncation ranks could result in a possible discrepancy in the gradient information and the optimal control problem then may not converge as desired. 
For all these reasons, here we instead focus on the use of reduced-order models as efficient surrogates within a classical Armijo backtracking line search. As this requires several forward evaluations in each optimization step, we can resort to the sPOD-G reduced-order state equation. Pseudocodes for such an approach, based on both the sPOD-G and the POD-G methods, are given in \Cref{alg:FOTR} and \Cref{alg:FOTR-PODG}, respectively. Regarding the former, the following remarks are in order. 
\begin{algorithm}  
  \caption{Optimal control with sPOD-G reduced two-way backtracking line search}\label{alg:FOTR}
  \begin{algorithmic}[1]
  \Require{$A$, $B$, $q_0$, $z_0$, $q_\mathrm{d}$, $\{\mathrm{p}^1_\mathrm{s}, \ldots, \mathrm{p}^K_\mathrm{s}\}$, $\mu$, $\omega^0$, $n_\mathrm{iter}$, $\delta$, $\beta$, $n_\mathrm{samples}$} 
    \State \textbf{Initialize:} $u=u_0$
      \For{$i = 1, \ldots, n_{\mathrm{iter}}$}
        \State $q^i = \textsc{State}(u^i)$ \Comment{Solve \eqref{eq_def:OC_FOM_1} and \eqref{eq_def:OC_FOM_2}}
        \State $p^i = \textsc{Adjoint}(q^i, u^i)$ \Comment{Solve \eqref{eq_def:OC_FOM_3} and \eqref{eq_def:OC_FOM_4}}
        \State $\tfrac{\mathrm{d}\mathcal{L}}{\mathrm{d}u^i}=\mu u^i + B^\top p^i$
        \State $V^i_\mathrm{s}, W^i_\mathrm{s} = \textsc{Basis} ([q^i(t_1) \: \ldots \: q^i(t_n)], \{\mathrm{p}^1_\mathrm{s}, \ldots, \mathrm{p}^K_\mathrm{s}\}, n_\mathrm{samples})$ \label{alg:FOTR_lbasis1} \Comment{Use \eqref{eq_def:sPOD-decomposition-discrete} and \eqref{eq_def:sPOD_galerkin_basis}}
        \State $a^i_\mathrm{s}, z^i_\mathrm{s} = \textsc{ReducedState}(V^i_\mathrm{s}, W^i_\mathrm{s}, u^i)$ \Comment{Solve \eqref{eq_def:OC_FOTR_sPODG_1} and \eqref{eq_def:OC_FOTR_sPODG_2}}
        \State $\omega^i = \textsc{StepSize}(\omega^{i- 1}, \tfrac{\mathrm{d}\mathcal{L}}{\mathrm{d}u^i}, V^i_\mathrm{s}, a^i_\mathrm{s}, z^i_\mathrm{s}, u^i)$  \label{alg:TWBT} 
        \State $u^{i+1} = u^i - \omega^i\left(\tfrac{\mathrm{d}\mathcal{L}}{\mathrm{d}u^i}\right)$ 
        \If{$i==n_{\mathrm{iter}}$} 
            \State break
        \ElsIf{$\norm{\tfrac{\mathrm{d}\mathcal{L}}{\mathrm{d}u^i}}_\mathrm{2} / \norm{\tfrac{\mathrm{d}\mathcal{L}}{\mathrm{d}u^1}}_\mathrm{2} < \delta$}
            \State set $u=u^{i+1}$ and return
        \EndIf
      \EndFor
  \Ensure{$u$}
  \end{algorithmic}
\end{algorithm}

\begin{rmrk}
    In \Cref{alg:FOTR} \cref{alg:TWBT} we use a slightly modified version of the two-way backtracking algorithm \cite{truong_backtracking_2021} which led to a noticeable speed up when compared to a conventional line-search method for selecting a suitable step size $\omega$.
\end{rmrk} 

\begin{rmrk}\label{rem:FOTR_sPODG_precompute}
    We observe that in the state equation \eqref{eq_def:OC_FOTR_sPODG_1} the coefficient matrices depend on $a_\mathrm{s}$ and $z_\mathrm{s}$. 
    Since these values change at each time step, repeatedly constructing these coefficient matrices could be time-consuming unless they are scaled with the reduced dimension. 
    To address this, we pre-construct the matrices dependent on the shift $z_\mathrm{s}$ (that scales with the FOM dimension) by sampling $n_{\mathrm{samples}}$ values of $z_\mathrm{s}$ from a sufficiently large sample space and then interpolate linearly to obtain their values during the time evolution. 
    However, since the matrix $D_\mathrm{s}$, which depends on $a_\mathrm{s}$, scales with the reduced dimension, it is computed dynamically during the time evolution.
    This approach is similar to the one used in \cite{black_efficient_2021}. 
\end{rmrk}

\begin{rmrk}
    We note here that in our example of the 1D linear advection equation with periodic boundary conditions, we only use a single co-moving frame ($K=1$) for the sPOD-G method. 
    Moreover, the right-hand side of the advection equation is also linear. 
    Now subject to this, we can make significant simplifications to the sPOD-G method presented here.
    To make the analysis plausible we consider the infinite-dimensional analog of the sPOD-ansatz \eqref{eq_def:sPOD-decomposition-discrete} for a single frame:
    \begin{equation}
        l(x,t) = \mathcal{T}(z(t))\phi(x) g(t)
    \end{equation}
    where $z(t)$ is the shift and $g(t)$, $l(x, t)$ are the reduced transported field and the transported field respectively.
    Following \cite{black_projection-based_2020} (Ex. 5.12) we assume that the transport operator $\mathcal{T}(z(t))$ which performs shift translations is globally isometric in Euclidean space such that $\norm{\mathcal{T}(z(t))(x) - \mathcal{T}(z(t))(y)}_2 = \norm{x + z(t) - y - z(t)}_2 = \norm{x - y}_2$ and thus $\mathcal{T}^{-1}\mathcal{T} = \mathcal{I}$ (this similarly holds true for rotations and reflections).
    Moreover, the right-hand side of the linear advection equation is equivariant with respect to this operator which is true for the negative advection and the constant control operator $-\nabla$ and $\mathcal{B}$ respectively.
    Moreover, for our example problem, we also get that $\mathcal{T}(-z)[\mathcal{T}'(z)\phi]$ and $\langle [\mathcal{T}'(z)\phi]g, [\mathcal{T}'(z)\psi]h\rangle_{L^2(0, t_f;\Omega)}$ do not depend upon $z$ for all $\phi, \psi \in L^2(\Omega)$ and $g, h \in \mathbb{R}$. 
    We now reference this result in the finite-dimensional setting and deduce that the shift-dependent matrices $M_{1, s},\: M_{2, s},\: N_s,\: A_{1, s},\: A_{2, s}$ are independent of $z$ and thus are constant.
    For a more detailed analysis of this aspect, we refer to \cite{black_projection-based_2020, black_efficient_2021}.  
\end{rmrk}

\begin{rmrk}\label{rem:adaptive_shift_refine}
    In \Cref{alg:FOTR} \cref{alg:FOTR_lbasis1}, the sPOD ansatz is used to decompose the snapshot matrix and obtain the basis $V_\mathrm{s}^i$. 
    At each optimization step, this requires computing the shifts $z_\mathrm{s}^i$ for $q^i$ and the associated transformation operators $T^{z_\mathrm{s}^i}$, which are sparse matrices that are computationally costly to assemble. 
    To reduce this expense, we adopt an adaptive shift update strategy. 
    Initially, shifts are fixed based on the uncontrolled state equation and reused until optimization stagnates. 
    At that point, shifts are updated based on the current $q^i$, and this process is repeated until convergence. 
    Numerical tests showed this approach reduces computational effort while maintaining satisfactory results.
\end{rmrk}

\begin{rmrk}\label{rem:FOTR_sPODG_precompute_parallel}
    As noted in Rem.~\autoref{rem:FOTR_sPODG_precompute}, achieving computational efficiency in the sPOD-G reduced-order model compared to POD-G requires pre-computing the shift-dependent matrices that scale with the FOM dimension and using linear interpolation during the reduced-order solve. 
    In optimal control settings, however, this pre-computation must be repeated at every optimization step due to changing dynamics, creating a computational bottleneck. 
    Since these matrices are computed independently for each shift value, we employed parallelization and achieved appropriate speed-ups in our numerical tests.
\end{rmrk}

\begin{algorithm}
  \caption{Optimal control with POD-G reduced two-way backtracking line search}\label{alg:FOTR-PODG}
  \begin{algorithmic}[1]
  \Require{$A$, $B$, $q_0$, $q_\mathrm{d}$, $\mathrm{p}$, $\mu$, $\omega^0$, $n_\mathrm{iter}$, $\delta$, $\beta$} 
    \State \textbf{Initialize} : $u=0$
      \For{$i = 1, \ldots, n_{\mathrm{iter}}$}
        \State $q^i = \textsc{State}(u^i)$ \Comment{Solve \eqref{eq_def:OC_FOM_1} and \eqref{eq_def:OC_FOM_2}}
        \State $p^i = \textsc{Adjoint}(q^i, u^i)$ \Comment{Solve \eqref{eq_def:OC_FOM_3} and \eqref{eq_def:OC_FOM_4}}
        \State $\tfrac{\mathrm{d}\mathcal{L}}{\mathrm{d}u^i}=\mu u^i + B^\top p^i$
        \State $U^i_\mathrm{p} = \textsc{Basis} ([q^i(t_1) \: \ldots \: q^i(t_n)], \mathrm{p})$ \Comment{Use \eqref{eq_def:POD}}
        \State $a^i_\mathrm{p} = \textsc{ReducedState}(U^i_\mathrm{p}, u^i)$ \Comment{Solve \eqref{eq_def:OC_FOTR_PODG_1} and \eqref{eq_def:OC_FOTR_PODG_2}}
        \State $\omega^i = \textsc{StepSize}(\omega^{i- 1}, \tfrac{\mathrm{d}\mathcal{L}}{\mathrm{d}u^i}, U^i_\mathrm{p}, a^i_\mathrm{p}, u^i)$
        \State $u^{i+1} = u^i - \omega^i\left(\tfrac{\mathrm{d}\mathcal{L}}{\mathrm{d}u^i}\right)$ 
        \If{$i==n_{\mathrm{iter}}$} 
            \State break
        \ElsIf{$\norm{\tfrac{\mathrm{d}\mathcal{L}}{\mathrm{d}u^i}}_\mathrm{2} / \norm{\tfrac{\mathrm{d}\mathcal{L}}{\mathrm{d}u^1}}_\mathrm{2} < \delta$}
            \State set $u=u^{i+1}$ and return
        \EndIf
      \EndFor
  \Ensure{$u$}
  \end{algorithmic}
\end{algorithm}

\section{Numerical results}\label{sec:results}
\begin{figure}[htp!]
    \centering
    \begin{subfigure}{0.45\textwidth}
        \centering
        \setlength\figureheight{0.7\linewidth}
        \setlength\figurewidth{0.7\linewidth}
        \input{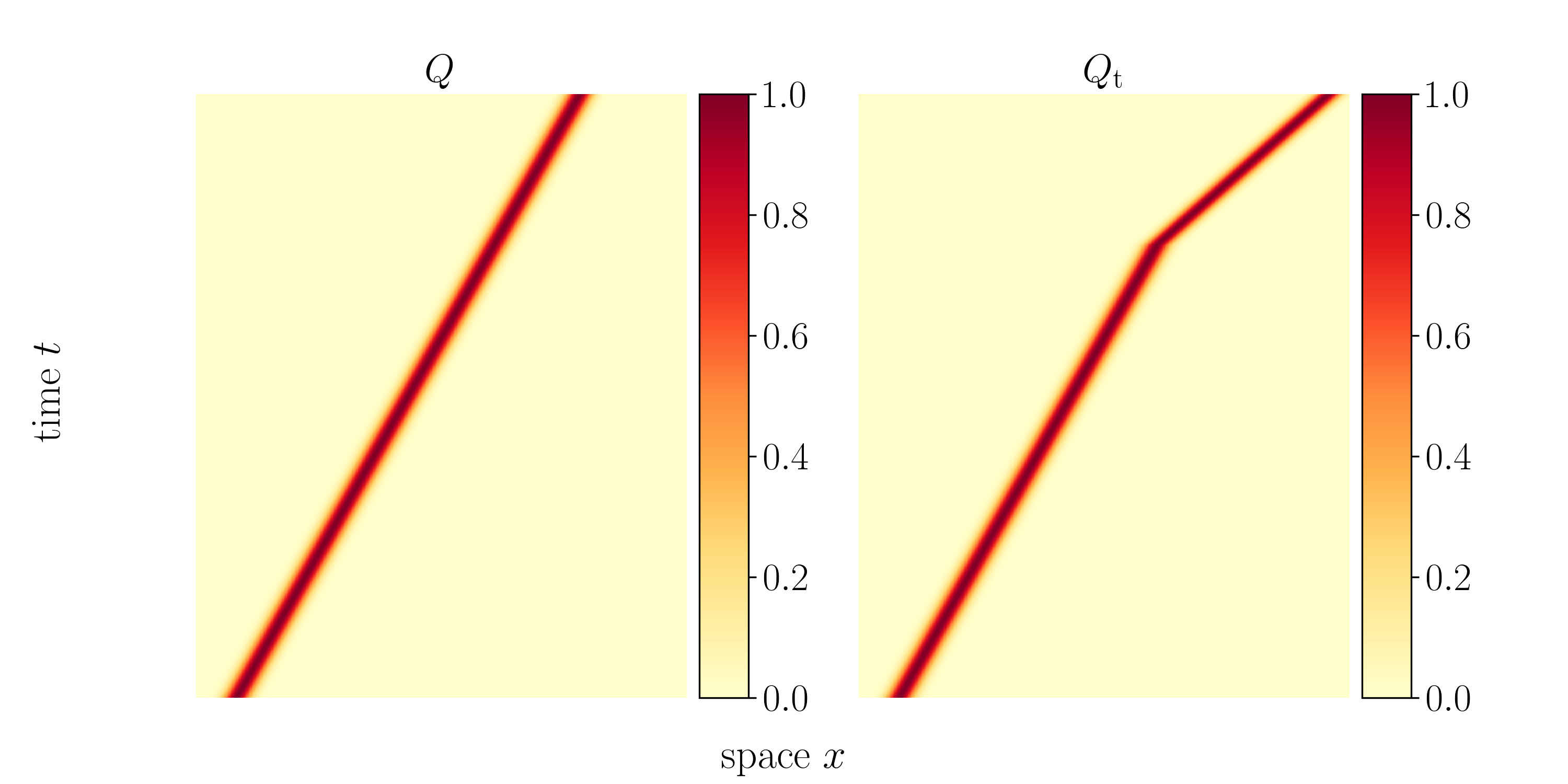} % Input the first plot
        \caption{Example 1: Broad wave with a kink at $\frac{3}{4}th$ of time domain.}
        \label{fig:Example_1}
    \end{subfigure}
    \hspace{0.05\textwidth}
    \begin{subfigure}{0.45\textwidth}
        \centering
        \setlength\figureheight{0.7\linewidth}
        \setlength\figurewidth{0.7\linewidth}
        \input{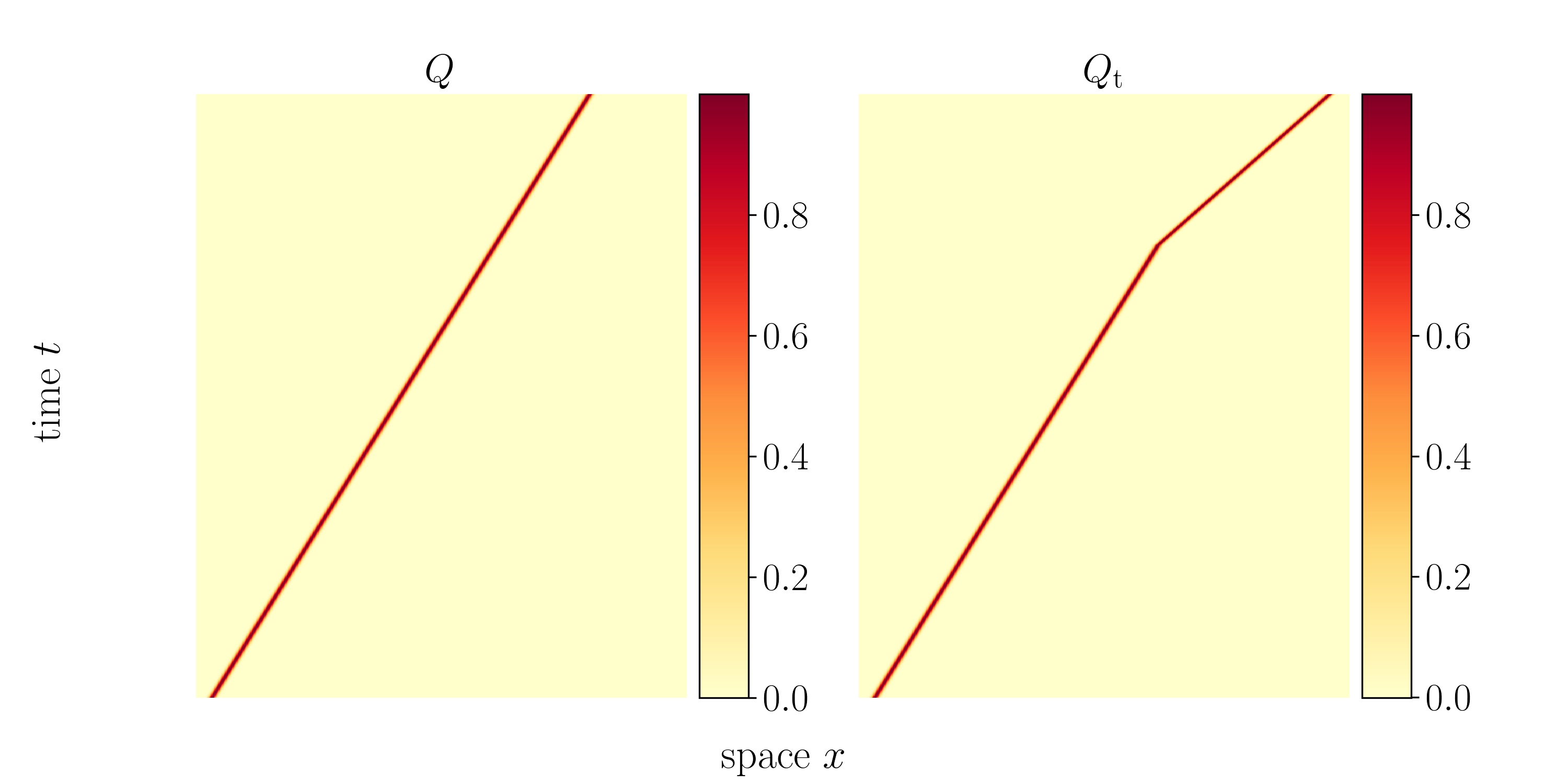} % Input the second plot
        \caption{Example 2: Sharp wave with a kink at $\frac{3}{4}th$ of time domain.}
        \label{fig:Example_2}
    \end{subfigure}
    \vspace{0.5cm} % Add space between rows
    \begin{subfigure}{0.45\textwidth}
        \centering
        \setlength\figureheight{0.7\linewidth}
        \setlength\figurewidth{0.7\linewidth}
        \input{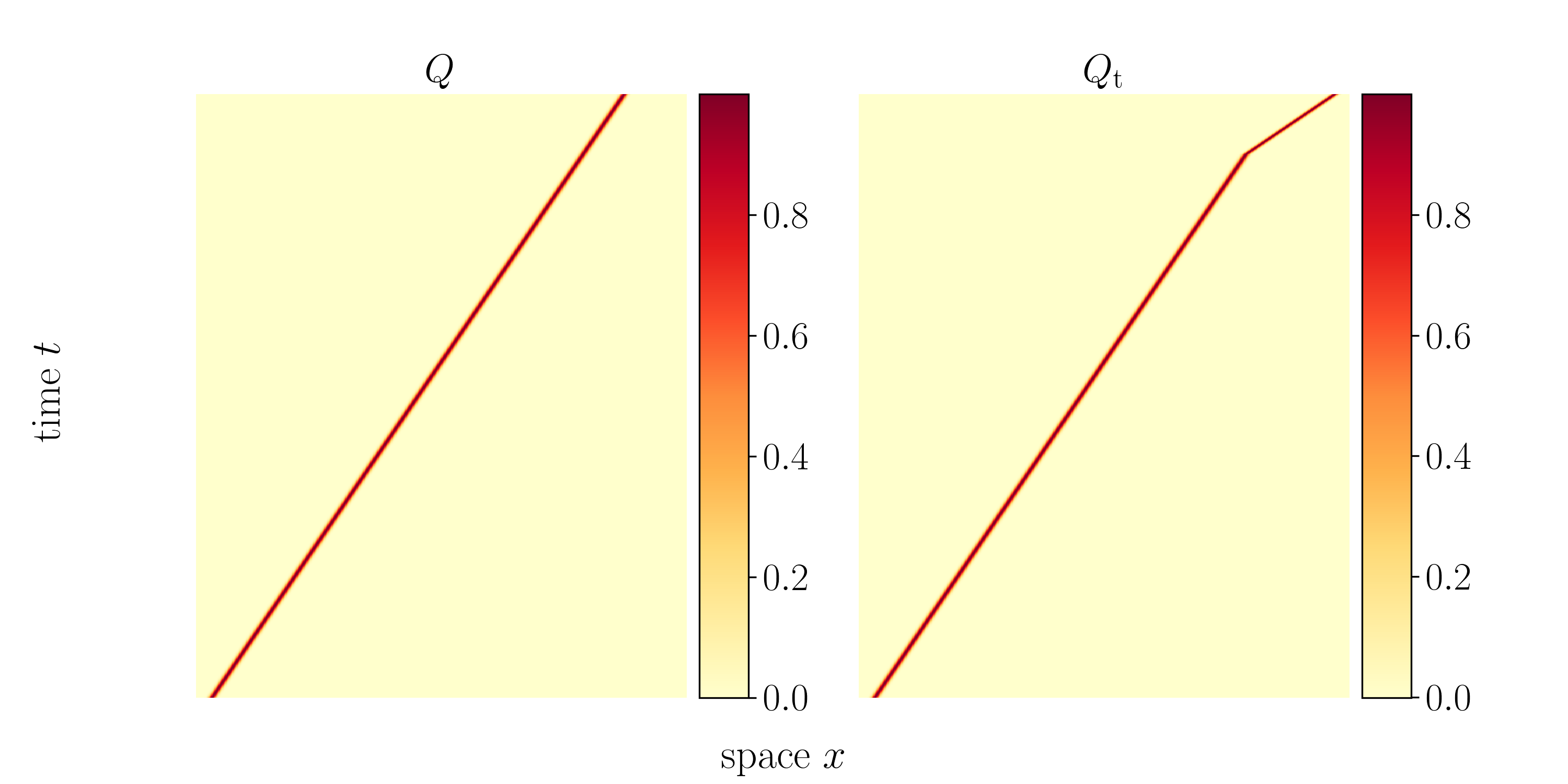} % Input the third plot
        \caption{Example 3: Sharp wave with a kink at $\frac{9}{10}th$ of time domain.}
        \label{fig:Example_3}
    \end{subfigure}
    \caption{Plots for the state and the target for all the three example problems}
    \label{fig:advection_and_target}
\end{figure}

All the numerical tests were run using Python 3.12 on a Macbook Air M1(2020) with an 8-core CPU and 16GB of RAM. 
In this section, we test the proposed methodologies on three different variations of the 1D linear advection example problem.
We consider a one-dimensional strip of length $l=100\:\mathrm{m}$ with $x\in ]0, l]$ discretized into $m=3200$ grid points with $\mathrm{d}x=0.03125 \:\mathrm{m}$ with periodic boundary conditions. 
The system \eqref{eq_def:FOM_discrete} is then simulated with $u=0$ for all $t\in [0, t_f[$ with  $t_f=140\:\mathrm{s}$ considering $n = 3360$ time steps and $\mathrm{d}t= 0.041 \: \mathrm{s}$ calculated with the help of the CFL criterion $\mathrm{d}t = \mathrm{cfl} \:\frac{\mathrm{d}x}{c}$, where we prescribe $\mathrm{cfl} = 4 / 3$ and the characteristic propagation speed of the quantity $q$ as $c=1 \: \textrm{m/s}$.
The solution obtained is shown in \Cref{fig:advection_and_target} along with the desired target profile for varying initial conditions. 
For the variations shown in \Cref{fig:advection_and_target} we consider two different initial conditions:
\begin{equation}\label{eq_def:advection_init}
    q_0 :=
    \begin{cases}
    \mathrm{exp}\left(- \frac{(x-l/12)^2}{7}\right) \:\: \text{for example 1 shown in \Cref{fig:Example_1}}, \\
    \\
    \mathrm{exp}\left(- \frac{(x-l/30)^2}{0.5}\right) \:\: \text{for example 2 and 3 shown in \Cref{fig:Example_2} and \Cref{fig:Example_3} resp.}
    \end{cases}
\end{equation}

For example 1, 2, and 3 we consider the advection velocities to be $v=0.5 c$, $0.55c$, and $0.6c$ respectively. 
Subsequently, we now solve the optimal control problem for achieving the desired target considering both the POD-G and sPOD-G methods.
In our tests, we consider the control shape functions $\mathcal{B}_k$ to be Gaussian functions given as: 
\begin{equation}\label{eq_def:mask_gaussian}
    \mathcal{B}_k(x) = \mathrm{exp}\left(- 4\left(x-\frac{L(k + 1)}{n_c}\right)^2\right)
\end{equation}
for $k=1, \ldots, n_c$ and $n_c=40$. 
The optimization parameters for the procedure are given in \Cref{tab:params}.
\begin{table}[h!]
\begin{center}
\begin{minipage}{\textwidth}
\centering
\begin{tabular}{|l||c|c|c|c|c|c|} 
 \hline
Parameters & $\mu$ & $\delta$ & $\omega$ & $\beta$ & $n_{\mathrm{iter}}$ & $n_{\mathrm{samples}}$ \\
\hline 
Values & $10^{-3}$ & $1\times 10^{-4}$ & $1$ & $0.5$ & $100000$ & $800$ \\
\hline  
\end{tabular}
\caption{Constant parameters used in the optimal control experiments}
\label{tab:params}
\end{minipage}
\end{center}
\end{table}

\begin{rmrk}
    To solve both the state and adjoint equations, we employ a 6th-order central finite difference scheme for spatial discretization and an explicit 4th-order Runge-Kutta (RK4) method for time integration. 
    The RK4 scheme requires variable values of $u$ in the state equation and ($q$, $q_\mathrm{d}$) in the adjoint equation at half-time steps, $t + \mathrm{d}t/2$, which we obtain through linear interpolation between $t$ and $t + \mathrm{d}t$. 
    However, this approach reduces the overall convergence order. Alternatively, implicit methods like Crank-Nicholson could be considered, but these also demand significant computation due to the need to solve a linear system at each optimization step. 
    Simpler explicit methods, such as the Euler scheme or RK2, are potential options but tend to be unstable for transport-dominated problems like ours and would require much finer time discretization, further increasing computational demands. 
    Ultimately, we strike a balance between convergence order, stability for transport-dominated systems, and computational efficiency by using RK4 with linear interpolation for the half-step values.
 \end{rmrk}

The selection of the specific value for $n_\mathrm{iter}$ depends on the convergence criteria met for the FOM with the given $\delta$. 
Consequently, we terminate the optimal control loop for the discussed techniques either upon reaching the prescribed $n_\mathrm{iter}$ or when the relative norm of the gradient of the Lagrangian falls below the specified $\delta$, whichever occurs first.

\begin{rmrk}
    We observed that for some test cases involving the POD-G and sPOD-G methods, the optimal control loop began to stagnate after a certain number of iterations. 
    This stagnation occurs when the prescribed number of truncated modes is insufficient to reduce the cost functional any further. 
    In such cases, the Armijo condition will fail to be satisfied and, as a result, the relative norm of the gradient of the Lagrangian does not decrease further.  In such cases, we terminated the optimal control loop early. 
\end{rmrk}

\begin{figure}[htp!]
    \centering
    \begin{subfigure}{0.49\textwidth}
        \centering
        \setlength\figureheight{0.8\linewidth}
        \setlength\figurewidth{1.0\linewidth}
        \input{ModesVsCost_1}   % Input the first plot
        \caption{Example 1}
        \label{fig:J_vs_modes_Example_1}
    \end{subfigure}
    \hspace{0.001\textwidth}
    \begin{subfigure}{0.49\textwidth}
        \centering
        \setlength\figureheight{0.8\linewidth}
        \setlength\figurewidth{1.0\linewidth}
        \input{ModesVsCost_2}   % Input the second plot
        \caption{Example 2}
        \label{fig:J_vs_modes_Example_2}
    \end{subfigure}
    \vspace{0.5cm} % Add space between rows
    \begin{subfigure}{0.49\textwidth}
        \centering
        \setlength\figureheight{0.8\linewidth}
        \setlength\figurewidth{1.0\linewidth}
        \input{ModesVsCost_3}   % Input the third plot
        \caption{Example 3}
        \label{fig:J_vs_modes_Example_3}
    \end{subfigure}
    \caption{Plots for $\mathcal{J}$ vs. ROM dimension for all the three examples}
    \label{fig:J_vs_modes}
\end{figure}
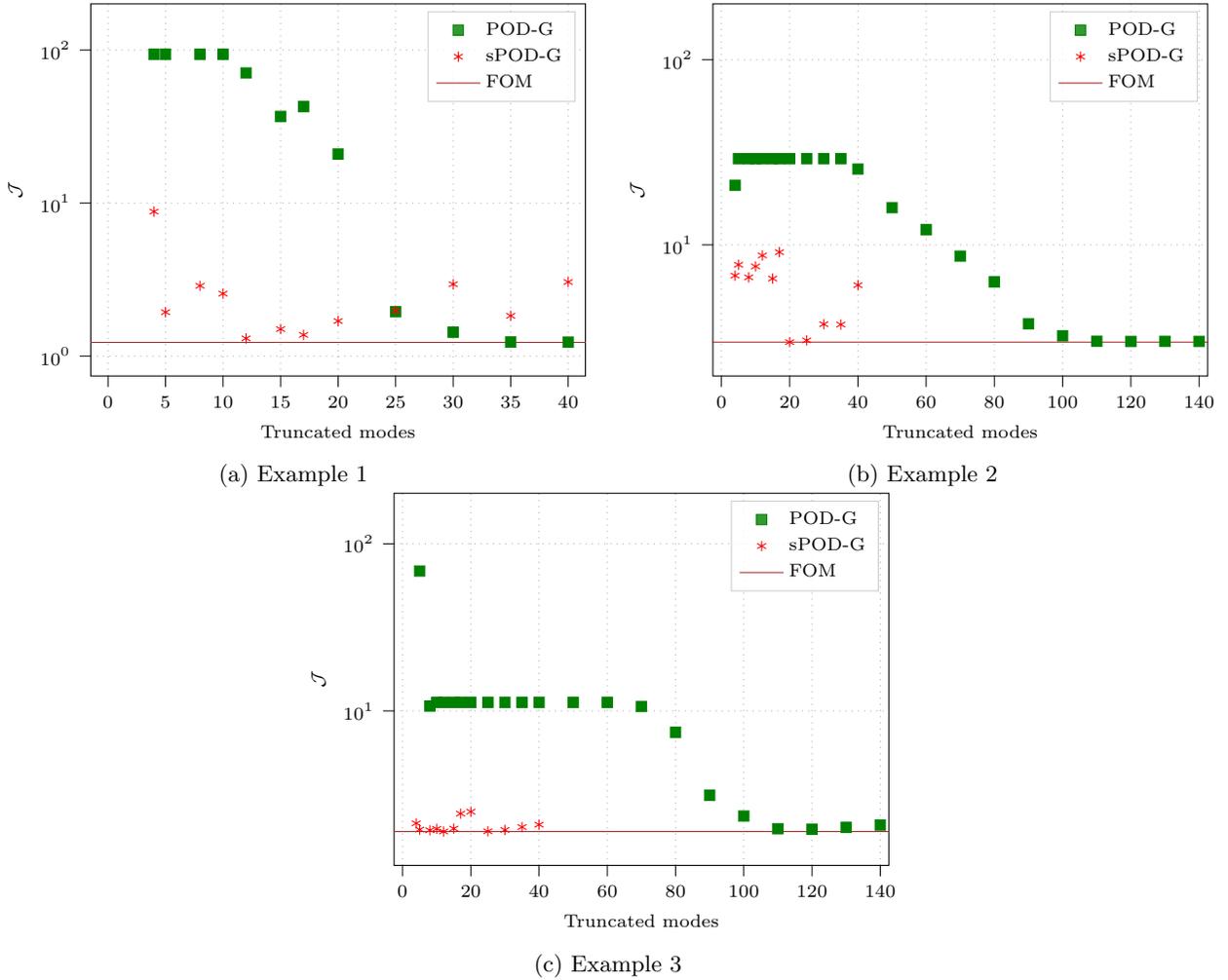

The numerical results for our example problems are shown in \Cref{fig:J_vs_modes} where we plot the value of the cost functional \eqref{eq_def:FOM_costFunc_discrete} with $u$ being the sub-optimal control obtained from either the POD-G or the sPOD-G method and $q$ obtained by solving \eqref{eq_def:FOM_discrete} with the sub-optimal control. 
For example 1, we run the tests for truncated modes $\mathrm{p}\in [1, 40]$ and $r \in [1, 40]$ for the POD-G and sPOD-G methods, respectively. Here, a comparably low number of modes is already sufficient to get the cost functional with sub-optimal control sufficiently close to the FOM cost functional. For examples 2 and 3, we considered truncated modes $\mathrm{p}\in [1, 140]$ and $r \in [1, 40]$ for the POD-G and sPOD-G methods, respectively, in order to come sufficiently close to the costs obtained with the full order optimization.  Furthermore, let us note that for the sPOD-G method, we have just a single co-moving frame ($K=1$) so that $r = \sum^{K}_{k=1}\mathrm{p}^k = \mathrm{p}^1$.

From \Cref{fig:J_vs_modes}, we observe that the sPOD-G method generally requires fewer truncated modes than POD-G to approach the FOM cost functional, significantly reducing the problem's dimensionality. 
However, the difference in required modes varies across the examples.
In example 1, the traveling wave is relatively broad despite being transport-dominated (\Cref{fig:advection_and_target}). 
Here, POD-G requires only about 35 modes to match the FOM cost functional (\Cref{fig:J_vs_modes_Example_1}), while sPOD-G achieves this with just 12 modes—a savings factor of 2.
In example 2, the wave is sharper, with the kink located at $\frac{3}{4}$th of the time domain. 
As shown in \Cref{fig:J_vs_modes_Example_2}, sPOD-G requires approximately $20$ modes, whereas POD-G needs close to $100$, a reduction by a factor of $5$ for sPOD-G.
In example 3, the wave sharpness remains the same, but the kink shifts later to $\frac{9}{10}$th of the time domain. 
From \Cref{fig:J_vs_modes_Example_3}, sPOD-G requires just $5–8$ modes, while POD-G requires nearly 20 times more modes, further demonstrating sPOD-G's efficiency.

\begin{rmrk}\label{rem:oscillations_kink}
    The variations in the advection equation example (\Cref{fig:advection_and_target}) illustrate the increasing difficulty for POD-G in capturing system dynamics with fewer truncated modes, while sPOD-G becomes progressively more effective from example 1 to example 3. 
    Example 3 is particularly favorable for sPOD-G because intermediate states during optimization exhibit diffused behavior, especially above the kink, where sharp fronts are absent. 
    This reduces the sPOD-G method's reliance on additional modes to handle oscillations.
    Since the kink occurs later in time for example 3 compared to example 2, the region susceptible to oscillations is smaller, allowing sPOD-G to perform well with fewer modes. 
    In contrast, POD-G struggles with sharp fronts, as seen in \Cref{fig:J_vs_modes_Example_2} and \Cref{fig:J_vs_modes_Example_3}, performing worse in Example 3 due to the reduced oscillatory regions, which are more favorable for sPOD-G.
\end{rmrk}

From the plots, we observe that the cost functional does not exhibit a monotonous decrease with an increasing number of modes for either method. 
This behavior is not surprising, as the ROM online error depends not only on the number of modes but also on factors like numerical instabilities and model truncation effects. 
Adding more modes can often lead to ill-conditioned systems, particularly in transport-dominated problems. 
Furthermore, higher modes are prone to capturing noise and irrelevant dynamics, which may degrade the ROM's accuracy.
Notably, the oscillations in the cost functional are more pronounced for the sPOD-G method compared to the POD-G method.
In the sPOD-G method, the $V_\mathrm{s}$ and $W_\mathrm{s}$ matrices \eqref{eq_def:sPOD_galerkin_basis} must be reassembled at each optimization step due to the dynamically changing stationary basis $U^k_{\mathrm{p}^k}$. 
As the number of modes increases, the higher modes tend to exhibit oscillatory behavior, which is further amplified when computing $W_\mathrm{s}$ since it involves taking spatial derivatives of these modes. This issue is absent in the POD-G method.
To address this challenge, one could use a very fine spatial grid to accurately capture the derivatives of the highly oscillatory modes 
(in our experiments a much coarser spatial discretization resulted in significantly stronger oscillatory effects).
Thus, the sPOD-G method requires balancing the need for a fine enough spatial resolution with computational feasibility along with keeping the mode numbers as low as possible.

\begin{rmrk}
    In our experiments, we consistently used non-overlapping Gaussian control shape functions 
    $\mathcal{B}_k(x)$ with a fixed variance, as defined in \eqref{eq_def:mask_gaussian}. Switching to other control shapes, such as indicator functions, adversely affected convergence. 
    Specifically, indicator functions introduced significant smearing and oscillatory behavior in the intermediate states during optimization, reducing the efficiency of the sPOD-G method as explained in Remark \autoref{rem:oscillations_kink}. 
    Additionally, the final optimal state computed by the FOM exhibited noticeable disturbances. To mitigate these effects, we chose smooth, non-overlapping Gaussian functions as the control shapes.
\end{rmrk}

Thus far, we have focused on numerical tests where the optimal control is performed using a fixed number of modes throughout the optimization process. 
However, it is not necessary to maintain the same number of modes at each step to accurately represent the dynamics. 
To address this, we implemented a tolerance-based strategy, dynamically selecting the number of truncated modes at each optimization step based on the  criterion:
\begin{equation}\label{eq_def:epsilon_criteria}
    d = \sum^{\min(m, n)}_{i=1} \mathds{1}\left(\frac{s^i}{s^1} > \epsilon\right)\,.
\end{equation}
Here $\{s_i\}$ are the singular values and for POD-G method $s_i \in \Sigma_\mathrm{p}$ and for the sPOD-G method $s_i \in \Sigma^1_{\mathrm{p}^1}$, and $\epsilon$ is a tolerance selected upfront. 
We tested the tolerance-based strategy for all the examples and the results are shown in \Cref{tab:tolerance_0.01} and \Cref{tab:tolerance_0.001}.
\begin{table}[tbp]
  \caption*{Numerical results for all the example problems with tolerance-based mode selection}
  \vspace{-0.3cm}
  \begin{center}
    \begin{minipage}[t]{0.48\linewidth}
      \centering
      \setlength\tabcolsep{6.5pt}
      \begin{tabular}{lllll}
        \toprule
        \multirow{1}{*}{\textbf{Example}} & Method  & $\mathcal{J}$  & $N_{avg}$ \\
        \midrule
        \multirow{3}{*}{Example 1}        & FOM     & 1.22       & --        \\
                                          & POD-G   & 1.24       & $\approx$ 43        \\
                                          & sPOD-G  & 2.53       & $\approx$ 12        \\
        \midrule
        \multirow{3}{*}{Example 2}        & FOM     & 2.98       & --        \\
                                          & POD-G   & 3.00       & $\approx$ 154        \\
                                          & sPOD-G  & 8.44       & $\approx$ 7        \\
        \midrule
        \multirow{3}{*}{Example 3}        & FOM     & 1.89       & --        \\
                                          & POD-G   & 1.96       & $\approx$ 165        \\
                                          & sPOD-G  & 1.97       & $\approx$ 10        \\
        \bottomrule
      \end{tabular}
      \caption{Tolerance $\epsilon = 0.01$}
      \label{tab:tolerance_0.01}
    \end{minipage}
    \hfill
    \begin{minipage}[t]{0.48\linewidth}
      \centering
      \setlength\tabcolsep{6.5pt}
      \begin{tabular}{lllll}
        \toprule
        \multirow{1}{*}{\textbf{Example}} & Method  & $\mathcal{J}$  & $N_{avg}$ \\
        \midrule
        \multirow{3}{*}{Example 1}        & FOM     & 1.22       & --        \\
                                          & POD-G   & 1.25       & $\approx$ 57        \\
                                          & sPOD-G  & 1.30       & $\approx$ 15        \\
        \midrule
        \multirow{3}{*}{Example 2}        & FOM     & 2.98       & --        \\
                                          & POD-G   & 2.98       & $\approx$ 191        \\
                                          & sPOD-G  & 3.85       & $\approx$ 33        \\
        \midrule
        \multirow{3}{*}{Example 3}        & FOM     & 1.89       & --        \\
                                          & POD-G   & 1.94       & $\approx$ 203        \\
                                          & sPOD-G  & 1.90       & $\approx$ 24        \\
        \bottomrule
      \end{tabular}
      \caption{Tolerance $\epsilon = 0.001$}
      \label{tab:tolerance_0.001}
    \end{minipage}
  \end{center}
\end{table}
These tables show the value of the cost functional reached with the obtained sub-optimal control from both methods for two different tolerance values. 
It also shows the average number of modes $N_{avg}$ required by the methods per optimization step.
We observe that in all of the examples, as we decrease the tolerance, the costs obtained based on the reduced-order surrogates get closer to the  FOM costs.
Interestingly, the POD-G method requires significantly more modes for examples 2 and 3, consistent with findings from the fixed-mode study. 
In contrast, the sPOD-G method consistently uses fewer modes. 
Finally, it is important to emphasize that this tolerance-based test is fundamentally different from the fixed-mode study shown in \Cref{fig:J_vs_modes}, and the mode numbers here cannot be directly compared with those in the fixed-mode study.

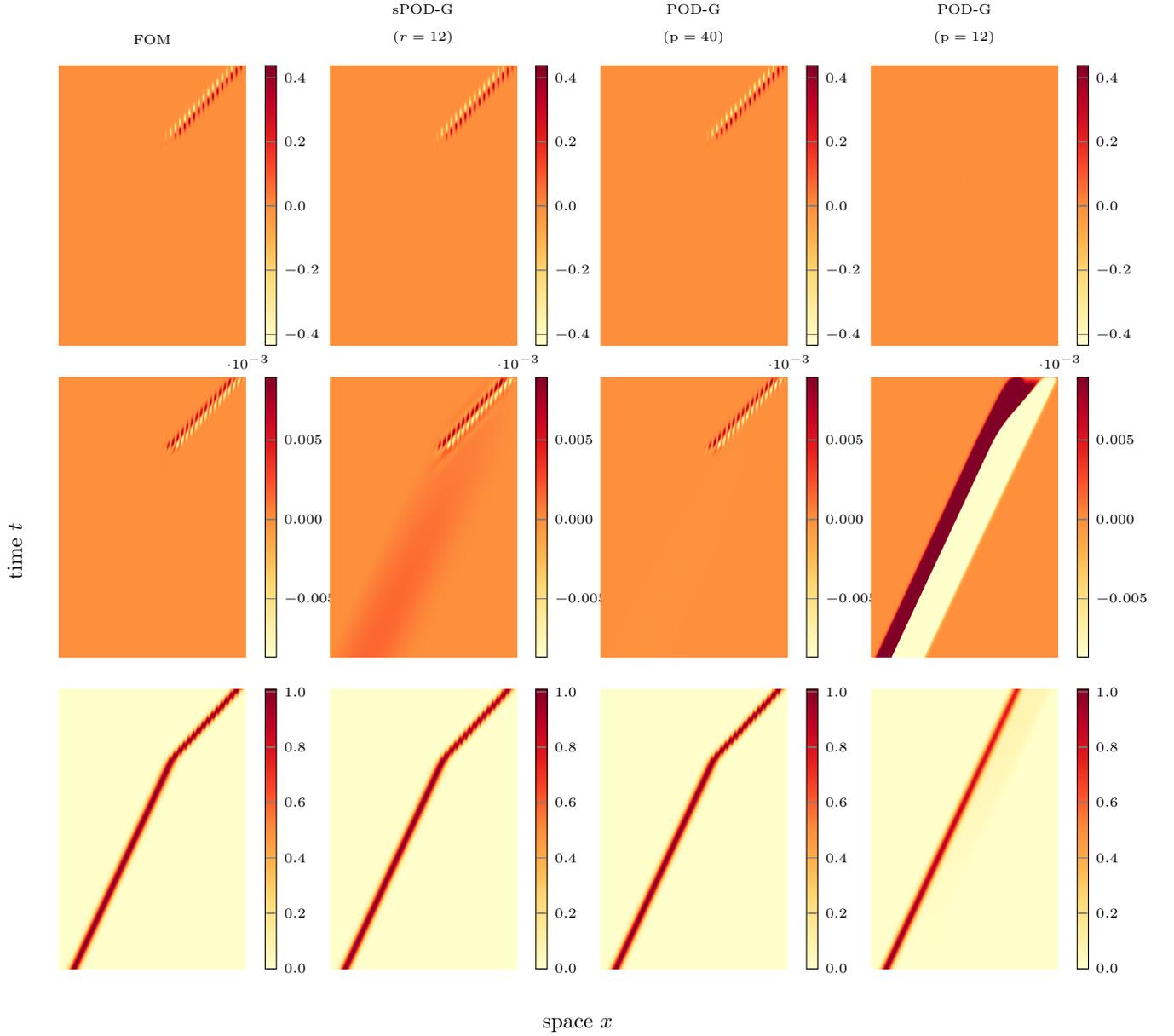
\begin{figure}[htp!]
  \centering
  \setlength\figureheight{0.35\linewidth}
  \setlength\figurewidth{0.265\linewidth}
  \input{Example_1_snapshots}
  \caption{Plots for the best controls, adjoints, and states obtained from the optimization procedure for example 1}
  \label{fig:example_1_snapshot_plots}
\end{figure}

\begin{figure}[htp!]
  \centering
  \setlength\figureheight{0.35\linewidth}
  \setlength\figurewidth{0.265\linewidth}
  \input{Example_2_snapshots}
  \caption{Plots for the best controls, adjoints, and states obtained from the optimization procedure for example 2}
  \label{fig:example_2_snapshot_plots}
\end{figure}
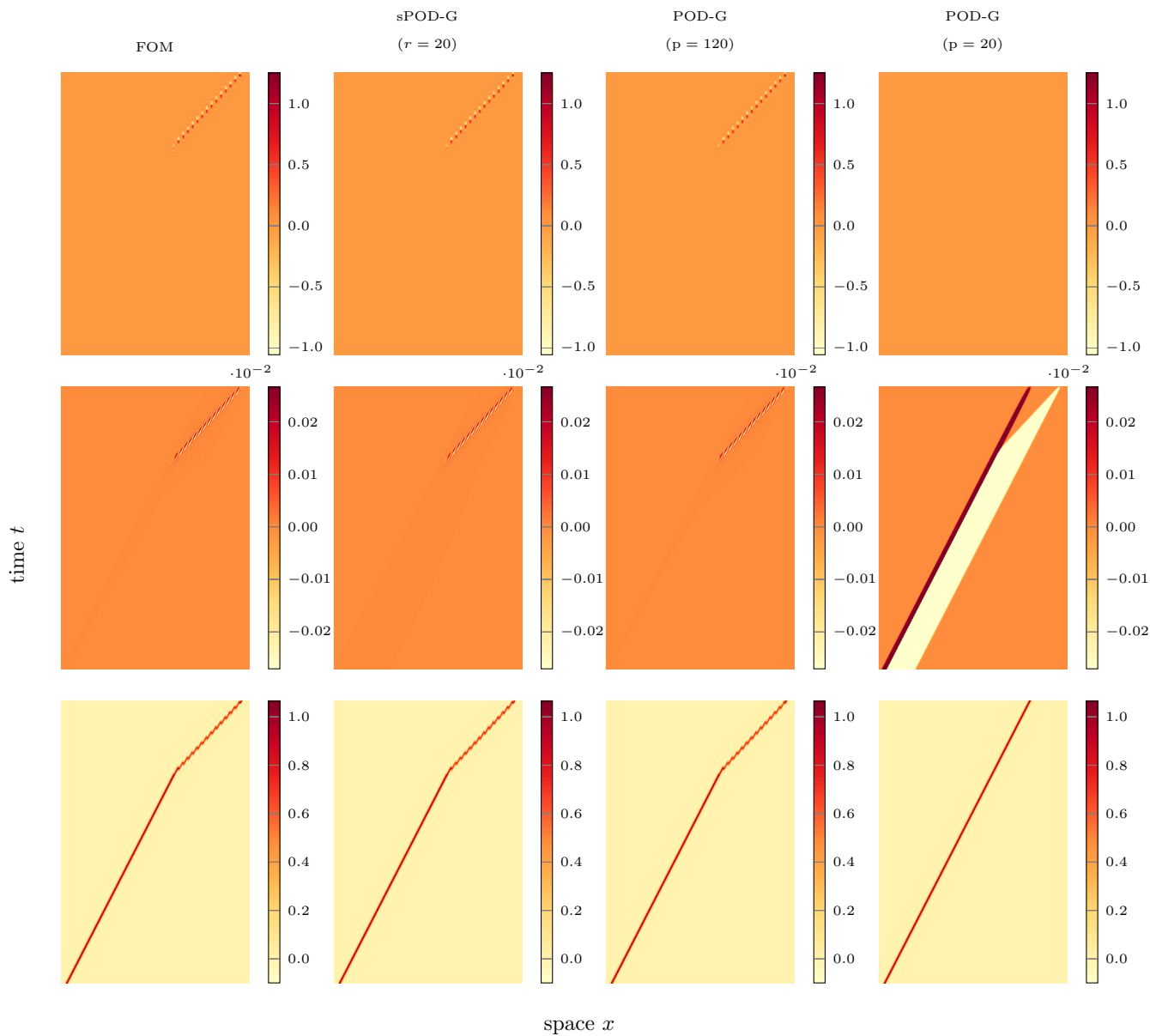

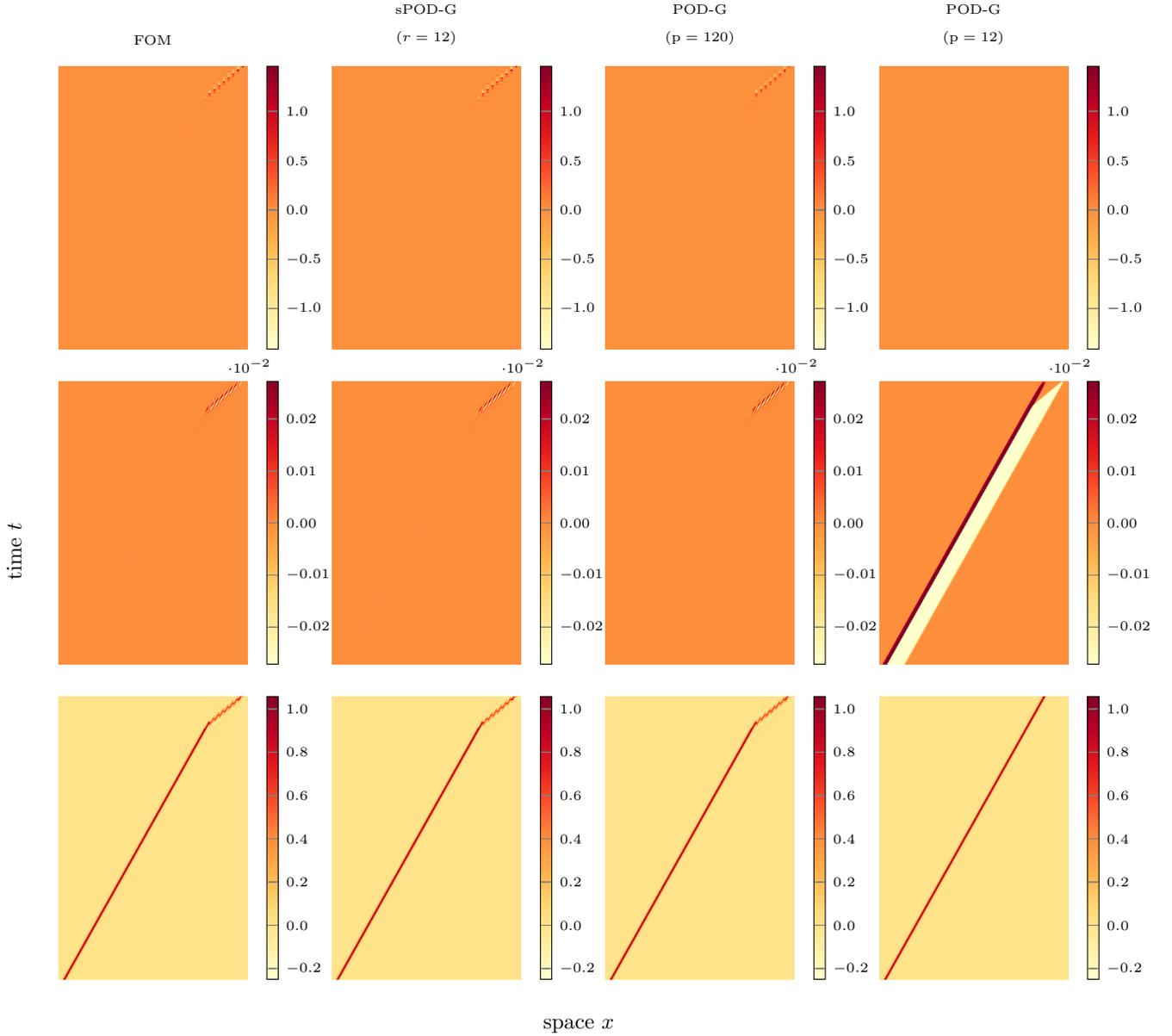
\begin{figure}[htp!]
  \centering
  \setlength\figureheight{0.35\linewidth}
  \setlength\figurewidth{0.265\linewidth}
  \input{Example_3_snapshots}
  \caption{Plots for the best controls, adjoints, and states obtained from the optimization procedure for example 3}
  \label{fig:example_3_snapshot_plots}
\end{figure}

For a visual comparison of the optimal controls, adjoints, and states, we refer to \Cref{fig:example_1_snapshot_plots}, which presents results for example 1.
The top row displays the sub-optimal control snapshot matrix $[Bu(t_1), \cdots, Bu(t_n)]\in \mathbb{R}^{m \times n}$, comparing the FOM results with the best results from both the sPOD-G and POD-G methods. 
Additionally, the last plot in the first row shows the sub-optimal control for POD-G using the same number of modes as the best result from the sPOD-G method. 
The second row compares the adjoints, while the third row illustrates the optimal states derived from the corresponding sub-optimal controls. 
Similar plots for examples 2 and 3 are shown in \Cref{fig:example_2_snapshot_plots} and \Cref{fig:example_3_snapshot_plots}.

\subsection{Timing analysis}
\begin{figure}[htp!]
    \centering
    \begin{subfigure}{0.49\textwidth}
        \centering
        \includegraphics[scale=0.42]{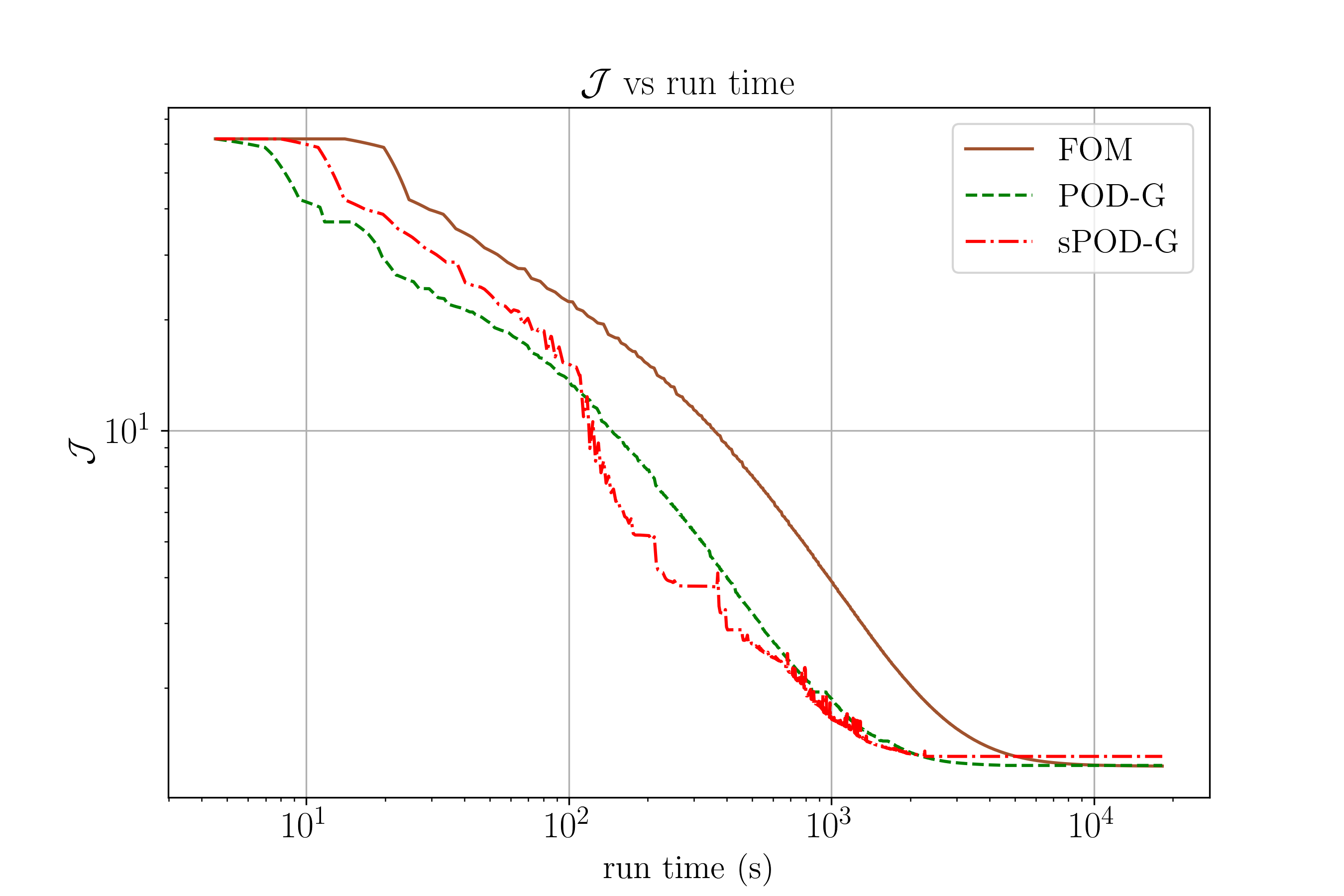}
        \caption{Example 1}
        \label{fig:J_vs_runtime_Example_1}
    \end{subfigure}
    \hspace{0.001\textwidth}
    \begin{subfigure}{0.49\textwidth}
        \centering
        \includegraphics[scale=0.42]{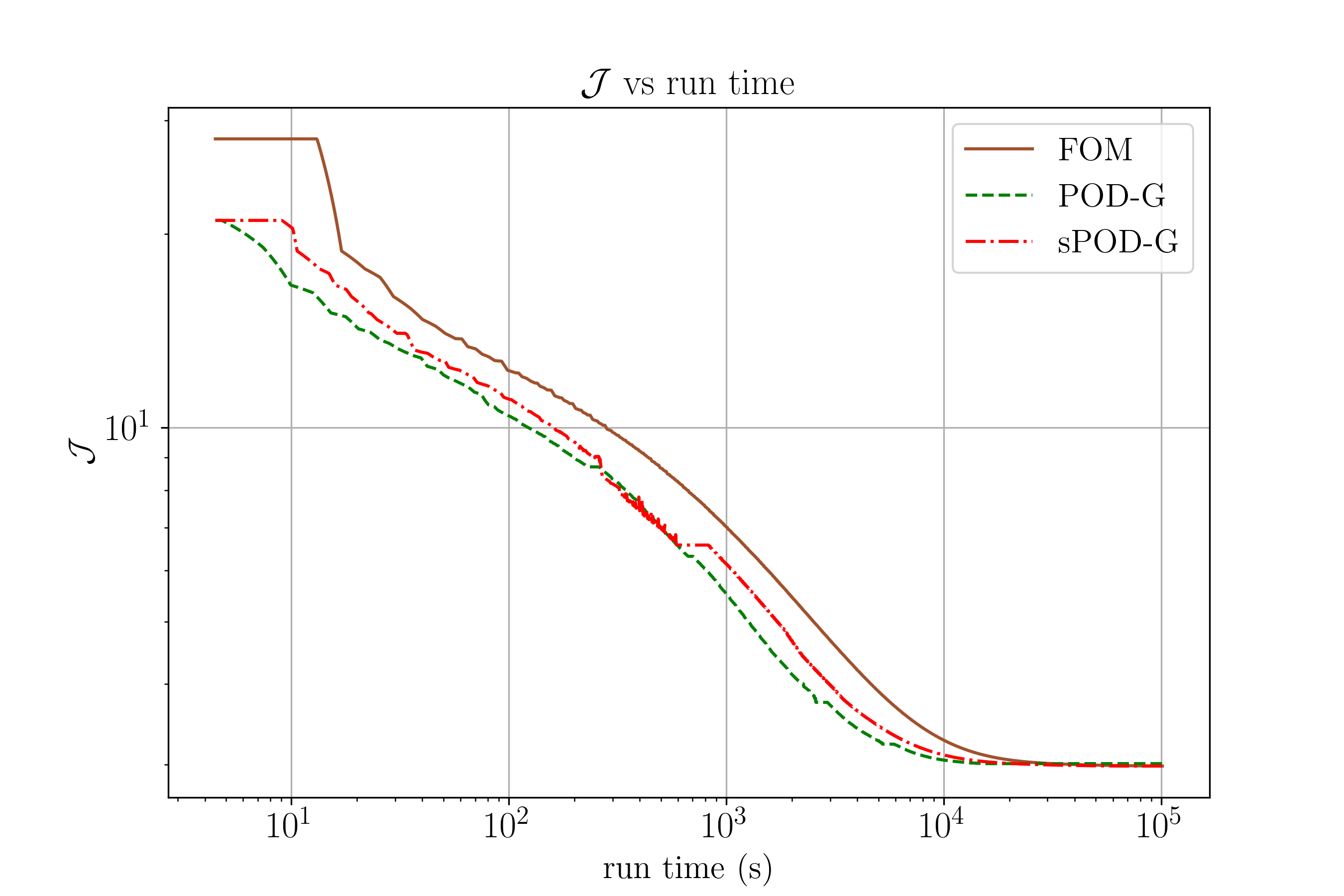}
        \caption{Example 2}
        \label{fig:J_vs_runtime_Example_2}
    \end{subfigure}
    \vspace{0.5cm} % Add space between rows
    \begin{subfigure}{0.49\textwidth}
        \centering
        \includegraphics[scale=0.42]{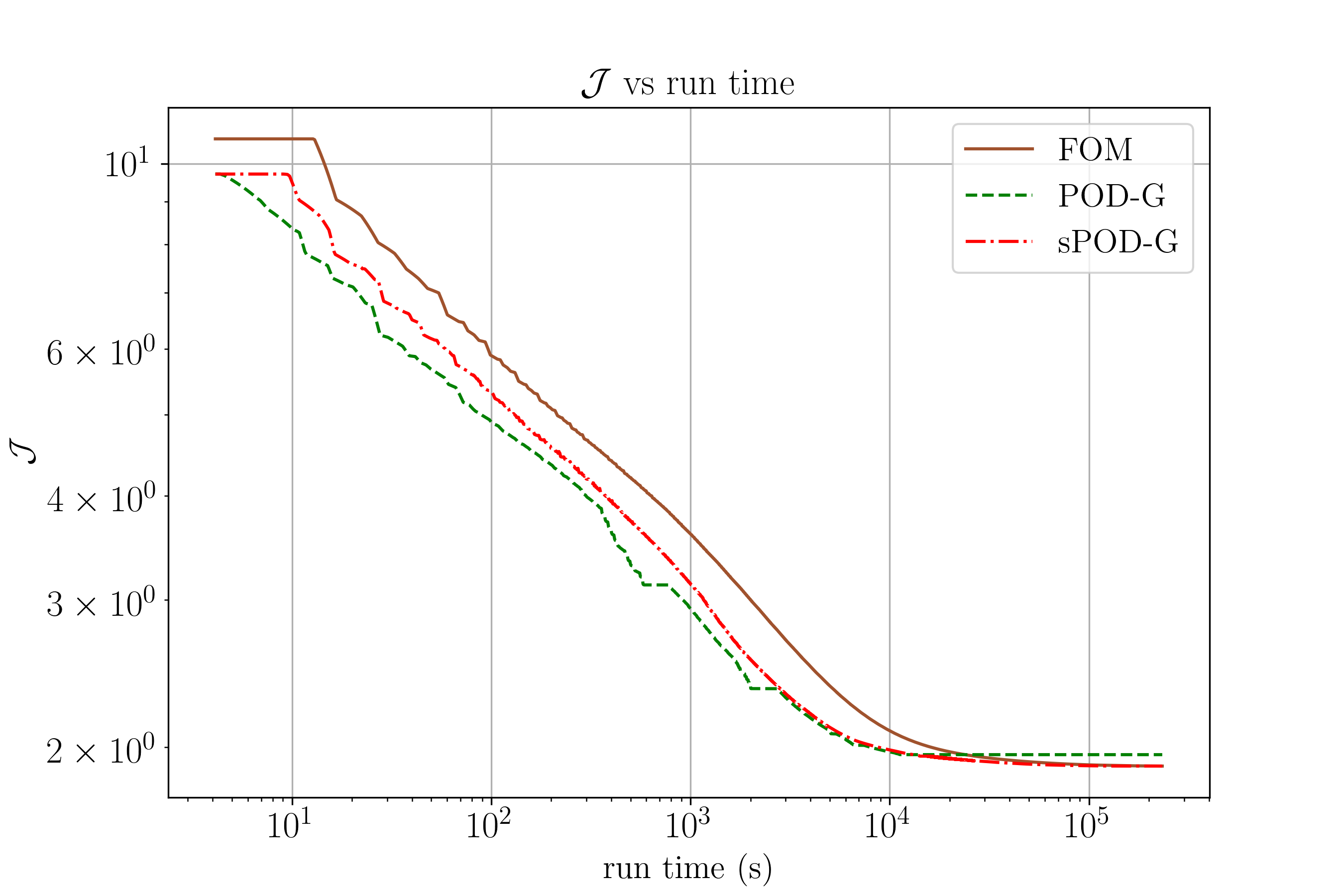}
        \caption{Example 3}
        \label{fig:J_vs_runtime_Example_3}
    \end{subfigure}
    \caption{Plots for $\mathcal{J}$ vs. run time for all the three examples}
    \label{fig:J_vs_runtime}
\end{figure}

So far, we have examined the performance differences between the POD-G and sPOD-G methods in terms of the reduced-order dimension. 
Let us also focus on a comparison between the two techniques based on their computational time.
Since the study for both the POD-G and sPOD-G methods involves varying the number of truncated modes, as shown in \Cref{fig:J_vs_modes}, the computational time for each mode also varies. In \Cref{fig:J_vs_runtime}, we also include a comparison of both methods with respect to the individual run times for all three examples.

We observe that the POD-G and sPOD-G methods are faster than the FOM in achieving a sufficiently low value of the cost functional $\mathcal{J}$. 
This speedup is particularly noteworthy given that both these methods perform basis refinement at every optimization step, which involves solving the FOM at each step. 
Even with this setup, the optimization process using reduced-order models is faster.
Additionally, we find that the runtimes for the POD-G and sPOD-G methods are comparable, with the POD-G method slightly outperforming the sPOD-G method in some cases.
While \Cref{fig:J_vs_modes} shows that sPOD-G requires significantly fewer modes compared to POD-G thus having a faster online stage, the offline computational costs of sPOD-G are higher than those of POD-G. 
This is due to the construction of the $V_\mathrm{s}$ and $W_\mathrm{s}$ matrices for all sampled shift values and their subsequent use in assembling the coefficient matrices for these samples at every optimization step, as noted in Remark \autoref{rem:FOTR_sPODG_precompute_parallel}. 
The construction of these coefficient matrices scales with the full-order dimension, contributing to the offline computational time, which is reflected in the total timing.
Additionally, the intermediate shift refinements further increase computational time since they involve assembling a large number of sparse matrices, as discussed in Remark \autoref{rem:adaptive_shift_refine}.

\section{Conclusion and outlook}\label{sec:conclusion}

In this paper, we formulated an open-loop optimal control problem for a transport-dominated PDE using reduced-order models (ROMs). 
We investigated two ROM approaches: the POD-G and sPOD-G methods. Additionally, we introduced two distinct frameworks for solving the optimal control problem: FOTR and FRTO.
We then assessed the computational performance and accuracy of the POD-G and sPOD-G methods within a line search based optimization loop for a linear advection equation. 
The numerical results presented in \Cref{sec:results} offer a comprehensive comparison between these ROM methods.
Our findings demonstrate that the sPOD-G method, which is particularly effective for transport-dominated systems, requires significantly fewer modes to achieve reasonable accuracy compared to the POD-G method.
Despite needing fewer modes, the sPOD-G method provides similar or slightly smaller speedups than the POD-G method, as shown in the timing analysis. 
This is primarily due to the higher offline computational costs associated with the sPOD-G method. 
As a future direction, we plan to investigate the numerical aspects of both the FOTR and the  FRTO framework in more detail. 
Additionally, we aim to study the effect of adaptive basis refinement on both the accuracy and computational time of the individual methods. 
Finally, we plan to explore the use of the sPOD-G method for optimal control in more complex PDE systems, such as a wildland fire model \cite{burela_parametric_2023}, which involves decomposing the problem into multiple frames using the sPOD ansatz.

\section*{Code Availability}

To enhance the reproducibility and transparency of the research in this paper, the source code for the experiments and analyses has been made publicly available via the following Zenodo repository:\\
\begin{center}
\urlstyle{tt}
\url{https://doi.org/10.5281/zenodo.14355726}
\end{center}

\vspace{1em}
We encourage researchers to utilize and build upon the code for their own research purposes.

\section*{Acknowledgement}
We gratefully acknowledge the support of the Deutsche Forschungsgemeinschaft (DFG) as part of GRK2433 DAEDALUS. 
We also acknowledge the financial support from the SFB TRR154 under the sub-project B03.
The authors were granted access to the HPC resources of the Institute of Mathematics at the TU Berlin.
We also thank Philipp Krah for his insightful comments and feedback and for providing access to the sPOD code base.

\addcontentsline{toc}{section}{Conflict Of Interest (COI)}
\section*{Conflict Of Interest (COI)}
All authors declare that they have no conflicts of interest.

%%-----------------------------
\bibliographystyle{plain}
\bibliography{abbr,refs}
%%-----------------------------
\end{document}

%% file: State_target_p1.tex
% This file was created with tikzplotlib v0.10.1.
\begin{tikzpicture}

\definecolor{darkgray176}{RGB}{176,176,176}

\begin{groupplot}[group style={group size=2 by 1, horizontal sep=1.2cm}]
\nextgroupplot[
colorbar,
colorbar style={ytick={-0.2,0,0.2,0.4,0.6,0.8,1,1.2},yticklabels={
  \(\displaystyle {\ensuremath{-}0.2}\),
  \(\displaystyle {0.0}\),
  \(\displaystyle {0.2}\),
  \(\displaystyle {0.4}\),
  \(\displaystyle {0.6}\),
  \(\displaystyle {0.8}\),
  \(\displaystyle {1.0}\),
  \(\displaystyle {1.2}\)
},ylabel={}, width=0.1*\pgfkeysvalueof{/pgfplots/parent axis width}, ticklabel style={font=\tiny}},
colormap={mymap}{[1pt]
  rgb(0pt)=(1,1,0.8);
  rgb(1pt)=(1,0.929411764705882,0.627450980392157);
  rgb(2pt)=(0.996078431372549,0.850980392156863,0.462745098039216);
  rgb(3pt)=(0.996078431372549,0.698039215686274,0.298039215686275);
  rgb(4pt)=(0.992156862745098,0.552941176470588,0.235294117647059);
  rgb(5pt)=(0.988235294117647,0.305882352941176,0.164705882352941);
  rgb(6pt)=(0.890196078431372,0.101960784313725,0.109803921568627);
  rgb(7pt)=(0.741176470588235,0,0.149019607843137);
  rgb(8pt)=(0.501960784313725,0,0.149019607843137)
},
height=1.2 * \figureheight,
hide x axis,
hide y axis,
point meta max=1.00000201032067,
point meta min=-1.56173253764891e-05,
tick align=outside,
tick pos=left,
title={\(\displaystyle Q\)},
width=0.7 * \figurewidth,
x grid style={darkgray176},
xmin=0, xmax=800,
xtick style={color=black},
xtick={0,200,400,600,800},
xticklabels={
  \(\displaystyle {0}\),
  \(\displaystyle {200}\),
  \(\displaystyle {400}\),
  \(\displaystyle {600}\),
  \(\displaystyle {800}\)
},
y grid style={darkgray176},
ymin=0, ymax=3360,
ytick style={color=black},
ytick={0,500,1000,1500,2000,2500,3000,3500},
yticklabels={
  \(\displaystyle {0}\),
  \(\displaystyle {500}\),
  \(\displaystyle {1000}\),
  \(\displaystyle {1500}\),
  \(\displaystyle {2000}\),
  \(\displaystyle {2500}\),
  \(\displaystyle {3000}\),
  \(\displaystyle {3500}\)
}
]
\addplot graphics [includegraphics cmd=\pgfimage,xmin=0, xmax=800, ymin=0, ymax=3360] {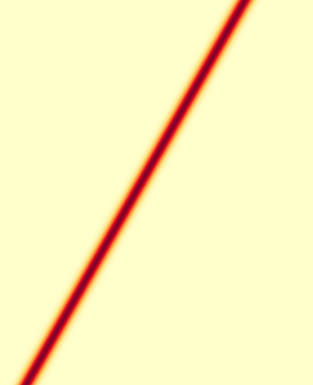};

\nextgroupplot[
colorbar,
colorbar style={ytick={-0.2,0,0.2,0.4,0.6,0.8,1,1.2},yticklabels={
  \(\displaystyle {\ensuremath{-}0.2}\),
  \(\displaystyle {0.0}\),
  \(\displaystyle {0.2}\),
  \(\displaystyle {0.4}\),
  \(\displaystyle {0.6}\),
  \(\displaystyle {0.8}\),
  \(\displaystyle {1.0}\),
  \(\displaystyle {1.2}\)
},ylabel={}, width=0.1*\pgfkeysvalueof{/pgfplots/parent axis width}, ticklabel style={font=\tiny}},
colormap={mymap}{[1pt]
  rgb(0pt)=(1,1,0.8);
  rgb(1pt)=(1,0.929411764705882,0.627450980392157);
  rgb(2pt)=(0.996078431372549,0.850980392156863,0.462745098039216);
  rgb(3pt)=(0.996078431372549,0.698039215686274,0.298039215686275);
  rgb(4pt)=(0.992156862745098,0.552941176470588,0.235294117647059);
  rgb(5pt)=(0.988235294117647,0.305882352941176,0.164705882352941);
  rgb(6pt)=(0.890196078431372,0.101960784313725,0.109803921568627);
  rgb(7pt)=(0.741176470588235,0,0.149019607843137);
  rgb(8pt)=(0.501960784313725,0,0.149019607843137)
},
height=1.2 * \figureheight,
hide x axis,
hide y axis,
point meta max=1.0000020161487,
point meta min=-1.56173253764891e-05,
tick align=outside,
tick pos=left,
title={\(\displaystyle Q_\mathrm{t}\)},
width=0.7 * \figurewidth,
x grid style={darkgray176},
xmin=0, xmax=800,
xtick style={color=black},
xtick={0,200,400,600,800},
xticklabels={
  \(\displaystyle {0}\),
  \(\displaystyle {200}\),
  \(\displaystyle {400}\),
  \(\displaystyle {600}\),
  \(\displaystyle {800}\)
},
y grid style={darkgray176},
ymin=0, ymax=3360,
ytick style={color=black},
ytick={0,500,1000,1500,2000,2500,3000,3500},
yticklabels={
  \(\displaystyle {0}\),
  \(\displaystyle {500}\),
  \(\displaystyle {1000}\),
  \(\displaystyle {1500}\),
  \(\displaystyle {2000}\),
  \(\displaystyle {2500}\),
  \(\displaystyle {3000}\),
  \(\displaystyle {3500}\)
}
]
\addplot graphics [includegraphics cmd=\pgfimage,xmin=0, xmax=800, ymin=0, ymax=3360] {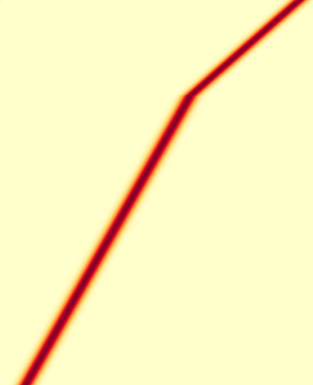};

\end{groupplot}

\draw ({$(current bounding box.south west)!-0.05!(current bounding box.south east)$}|-{$(current bounding box.south west)!0.3!(current bounding box.north west)$}) node[
  scale=0.96,
  anchor=west,
  text=black,
  rotate=90.0
]{time $t$};
\draw ({$(current bounding box.south west)!0.5!(current bounding box.south east)$}|-{$(current bounding box.south west)!-0.2!(current bounding box.north west)$}) node[
  scale=0.96,
  anchor=south,
  text=black,
  rotate=0.0
]{space $x$};
\end{tikzpicture}

%% file: State_target_p2.tex
% This file was created with tikzplotlib v0.10.1.
\begin{tikzpicture}

\definecolor{darkgray176}{RGB}{176,176,176}

\begin{groupplot}[group style={group size=2 by 1, horizontal sep=1.2cm}]
\nextgroupplot[
colorbar,
colorbar style={ytick={-0.2,0,0.2,0.4,0.6,0.8,1},yticklabels={
  \(\displaystyle {\ensuremath{-}0.2}\),
  \(\displaystyle {0.0}\),
  \(\displaystyle {0.2}\),
  \(\displaystyle {0.4}\),
  \(\displaystyle {0.6}\),
  \(\displaystyle {0.8}\),
  \(\displaystyle {1.0}\)
},ylabel={}, width=0.1*\pgfkeysvalueof{/pgfplots/parent axis width}, ticklabel style={font=\tiny}},
colormap={mymap}{[1pt]
  rgb(0pt)=(1,1,0.8);
  rgb(1pt)=(1,0.929411764705882,0.627450980392157);
  rgb(2pt)=(0.996078431372549,0.850980392156863,0.462745098039216);
  rgb(3pt)=(0.996078431372549,0.698039215686274,0.298039215686275);
  rgb(4pt)=(0.992156862745098,0.552941176470588,0.235294117647059);
  rgb(5pt)=(0.988235294117647,0.305882352941176,0.164705882352941);
  rgb(6pt)=(0.890196078431372,0.101960784313725,0.109803921568627);
  rgb(7pt)=(0.741176470588235,0,0.149019607843137);
  rgb(8pt)=(0.501960784313725,0,0.149019607843137)
},
height=1.2 * \figureheight,
hide x axis,
hide y axis,
point meta max=0.999999476581004,
point meta min=-0.000318742791436316,
tick align=outside,
tick pos=left,
title={\(\displaystyle Q\)},
width=0.7 * \figurewidth,
x grid style={darkgray176},
xmin=0, xmax=800,
xtick style={color=black},
xtick={0,200,400,600,800},
xticklabels={
  \(\displaystyle {0}\),
  \(\displaystyle {200}\),
  \(\displaystyle {400}\),
  \(\displaystyle {600}\),
  \(\displaystyle {800}\)
},
y grid style={darkgray176},
ymin=0, ymax=3360,
ytick style={color=black},
ytick={0,500,1000,1500,2000,2500,3000,3500},
yticklabels={
  \(\displaystyle {0}\),
  \(\displaystyle {500}\),
  \(\displaystyle {1000}\),
  \(\displaystyle {1500}\),
  \(\displaystyle {2000}\),
  \(\displaystyle {2500}\),
  \(\displaystyle {3000}\),
  \(\displaystyle {3500}\)
}
]
\addplot graphics [includegraphics cmd=\pgfimage,xmin=0, xmax=800, ymin=0, ymax=3360] {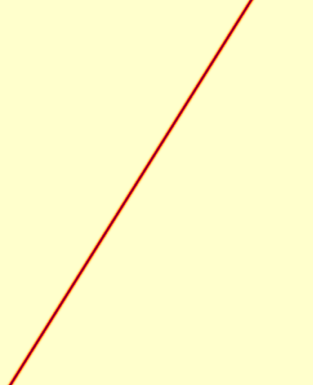};

\nextgroupplot[
colorbar,
colorbar style={ytick={-0.2,0,0.2,0.4,0.6,0.8,1},yticklabels={
  \(\displaystyle {\ensuremath{-}0.2}\),
  \(\displaystyle {0.0}\),
  \(\displaystyle {0.2}\),
  \(\displaystyle {0.4}\),
  \(\displaystyle {0.6}\),
  \(\displaystyle {0.8}\),
  \(\displaystyle {1.0}\)
},ylabel={}, width=0.1*\pgfkeysvalueof{/pgfplots/parent axis width}, ticklabel style={font=\tiny}},
colormap={mymap}{[1pt]
  rgb(0pt)=(1,1,0.8);
  rgb(1pt)=(1,0.929411764705882,0.627450980392157);
  rgb(2pt)=(0.996078431372549,0.850980392156863,0.462745098039216);
  rgb(3pt)=(0.996078431372549,0.698039215686274,0.298039215686275);
  rgb(4pt)=(0.992156862745098,0.552941176470588,0.235294117647059);
  rgb(5pt)=(0.988235294117647,0.305882352941176,0.164705882352941);
  rgb(6pt)=(0.890196078431372,0.101960784313725,0.109803921568627);
  rgb(7pt)=(0.741176470588235,0,0.149019607843137);
  rgb(8pt)=(0.501960784313725,0,0.149019607843137)
},
height=1.2 * \figureheight,
hide x axis,
hide y axis,
point meta max=0.999999476581004,
point meta min=-0.000635634865249002,
tick align=outside,
tick pos=left,
title={\(\displaystyle Q_\mathrm{t}\)},
width=0.7 * \figurewidth,
x grid style={darkgray176},
xmin=0, xmax=800,
xtick style={color=black},
xtick={0,200,400,600,800},
xticklabels={
  \(\displaystyle {0}\),
  \(\displaystyle {200}\),
  \(\displaystyle {400}\),
  \(\displaystyle {600}\),
  \(\displaystyle {800}\)
},
y grid style={darkgray176},
ymin=0, ymax=3360,
ytick style={color=black},
ytick={0,500,1000,1500,2000,2500,3000,3500},
yticklabels={
  \(\displaystyle {0}\),
  \(\displaystyle {500}\),
  \(\displaystyle {1000}\),
  \(\displaystyle {1500}\),
  \(\displaystyle {2000}\),
  \(\displaystyle {2500}\),
  \(\displaystyle {3000}\),
  \(\displaystyle {3500}\)
}
]
\addplot graphics [includegraphics cmd=\pgfimage,xmin=0, xmax=800, ymin=0, ymax=3360] {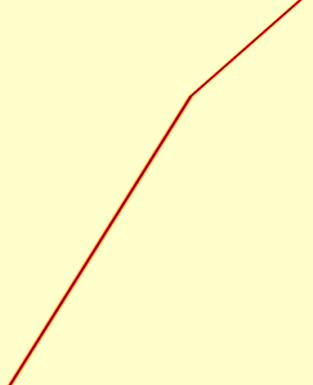};

\end{groupplot}

\draw ({$(current bounding box.south west)!-0.05!(current bounding box.south east)$}|-{$(current bounding box.south west)!0.3!(current bounding box.north west)$}) node[
  scale=0.96,
  anchor=west,
  text=black,
  rotate=90.0
]{time $t$};
\draw ({$(current bounding box.south west)!0.5!(current bounding box.south east)$}|-{$(current bounding box.south west)!-0.2!(current bounding box.north west)$}) node[
  scale=0.96,
  anchor=south,
  text=black,
  rotate=0.0
]{space $x$};
\end{tikzpicture}

%% file: State_target_p3.tex
% This file was created with tikzplotlib v0.10.1.
\begin{tikzpicture}

\definecolor{darkgray176}{RGB}{176,176,176}

\begin{groupplot}[group style={group size=2 by 1, horizontal sep=1.2cm}]
\nextgroupplot[
colorbar,
colorbar style={ytick={-0.2,0,0.2,0.4,0.6,0.8,1},yticklabels={
  \(\displaystyle {\ensuremath{-}0.2}\),
  \(\displaystyle {0.0}\),
  \(\displaystyle {0.2}\),
  \(\displaystyle {0.4}\),
  \(\displaystyle {0.6}\),
  \(\displaystyle {0.8}\),
  \(\displaystyle {1.0}\)
},ylabel={}, width=0.1*\pgfkeysvalueof{/pgfplots/parent axis width}, ticklabel style={font=\tiny}},
colormap={mymap}{[1pt]
  rgb(0pt)=(1,1,0.8);
  rgb(1pt)=(1,0.929411764705882,0.627450980392157);
  rgb(2pt)=(0.996078431372549,0.850980392156863,0.462745098039216);
  rgb(3pt)=(0.996078431372549,0.698039215686274,0.298039215686275);
  rgb(4pt)=(0.992156862745098,0.552941176470588,0.235294117647059);
  rgb(5pt)=(0.988235294117647,0.305882352941176,0.164705882352941);
  rgb(6pt)=(0.890196078431372,0.101960784313725,0.109803921568627);
  rgb(7pt)=(0.741176470588235,0,0.149019607843137);
  rgb(8pt)=(0.501960784313725,0,0.149019607843137)
},
height=1.2 * \figureheight,
hide x axis,
hide y axis,
point meta max=0.999874003918197,
point meta min=-0.000387975143581649,
tick align=outside,
tick pos=left,
title={\(\displaystyle Q\)},
width=0.7 * \figurewidth,
x grid style={darkgray176},
xmin=0, xmax=800,
xtick style={color=black},
xtick={0,200,400,600,800},
xticklabels={
  \(\displaystyle {0}\),
  \(\displaystyle {200}\),
  \(\displaystyle {400}\),
  \(\displaystyle {600}\),
  \(\displaystyle {800}\)
},
y grid style={darkgray176},
ymin=0, ymax=3360,
ytick style={color=black},
ytick={0,500,1000,1500,2000,2500,3000,3500},
yticklabels={
  \(\displaystyle {0}\),
  \(\displaystyle {500}\),
  \(\displaystyle {1000}\),
  \(\displaystyle {1500}\),
  \(\displaystyle {2000}\),
  \(\displaystyle {2500}\),
  \(\displaystyle {3000}\),
  \(\displaystyle {3500}\)
}
]
\addplot graphics [includegraphics cmd=\pgfimage,xmin=0, xmax=800, ymin=0, ymax=3360] {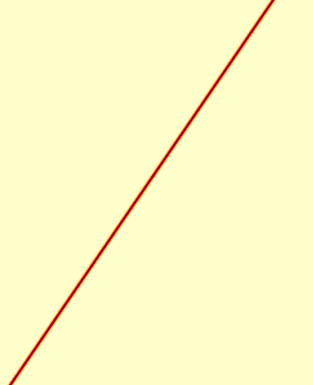};

\nextgroupplot[
colorbar,
colorbar style={ytick={-0.2,0,0.2,0.4,0.6,0.8,1},yticklabels={
  \(\displaystyle {\ensuremath{-}0.2}\),
  \(\displaystyle {0.0}\),
  \(\displaystyle {0.2}\),
  \(\displaystyle {0.4}\),
  \(\displaystyle {0.6}\),
  \(\displaystyle {0.8}\),
  \(\displaystyle {1.0}\)
},ylabel={}, width=0.1*\pgfkeysvalueof{/pgfplots/parent axis width}, ticklabel style={font=\tiny}},
colormap={mymap}{[1pt]
  rgb(0pt)=(1,1,0.8);
  rgb(1pt)=(1,0.929411764705882,0.627450980392157);
  rgb(2pt)=(0.996078431372549,0.850980392156863,0.462745098039216);
  rgb(3pt)=(0.996078431372549,0.698039215686274,0.298039215686275);
  rgb(4pt)=(0.992156862745098,0.552941176470588,0.235294117647059);
  rgb(5pt)=(0.988235294117647,0.305882352941176,0.164705882352941);
  rgb(6pt)=(0.890196078431372,0.101960784313725,0.109803921568627);
  rgb(7pt)=(0.741176470588235,0,0.149019607843137);
  rgb(8pt)=(0.501960784313725,0,0.149019607843137)
},
height=1.2 * \figureheight,
hide x axis,
hide y axis,
point meta max=0.999874003918197,
point meta min=-0.000678053820014014,
tick align=outside,
tick pos=left,
title={\(\displaystyle Q_\mathrm{t}\)},
width=0.7 * \figurewidth,
x grid style={darkgray176},
xmin=0, xmax=800,
xtick style={color=black},
xtick={0,200,400,600,800},
xticklabels={
  \(\displaystyle {0}\),
  \(\displaystyle {200}\),
  \(\displaystyle {400}\),
  \(\displaystyle {600}\),
  \(\displaystyle {800}\)
},
y grid style={darkgray176},
ymin=0, ymax=3360,
ytick style={color=black},
ytick={0,500,1000,1500,2000,2500,3000,3500},
yticklabels={
  \(\displaystyle {0}\),
  \(\displaystyle {500}\),
  \(\displaystyle {1000}\),
  \(\displaystyle {1500}\),
  \(\displaystyle {2000}\),
  \(\displaystyle {2500}\),
  \(\displaystyle {3000}\),
  \(\displaystyle {3500}\)
}
]
\addplot graphics [includegraphics cmd=\pgfimage,xmin=0, xmax=800, ymin=0, ymax=3360] {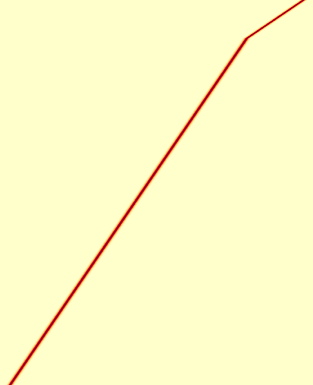};

\end{groupplot}

\draw ({$(current bounding box.south west)!-0.05!(current bounding box.south east)$}|-{$(current bounding box.south west)!0.3!(current bounding box.north west)$}) node[
  scale=0.96,
  anchor=west,
  text=black,
  rotate=90.0
]{time $t$};
\draw ({$(current bounding box.south west)!0.5!(current bounding box.south east)$}|-{$(current bounding box.south west)!-0.2!(current bounding box.north west)$}) node[
  scale=0.96,
  anchor=south,
  text=black,
  rotate=0.0
]{space $x$};
\end{tikzpicture}

%% file: ModesVsCost_1.tex
% This file was created with tikzplotlib v0.10.1.
\begin{tikzpicture}

\definecolor{brown}{RGB}{165,42,42}
\definecolor{darkgray176}{RGB}{176,176,176}
\definecolor{green}{RGB}{0,128,0}
\definecolor{lightgray204}{RGB}{204,204,204}

\begin{axis}[
height=\figureheight,
legend cell align={left},
legend style={fill opacity=0.8, draw opacity=1, text opacity=1, draw=lightgray204},
log basis y={10},
tick align=outside,
tick pos=left,
width=\figurewidth,
x grid style={darkgray176},
xlabel={Truncated modes},
xmajorgrids,
xmin=-1.5, xmax=41.5,
xtick style={color=black},
xtick={0,5,10,15,20,25,30,35,40},
xticklabels={
  \(\displaystyle {0}\),
  \(\displaystyle {5}\),
  \(\displaystyle {10}\),
  \(\displaystyle {15}\),
  \(\displaystyle {20}\),
  \(\displaystyle {25}\),
  \(\displaystyle {30}\),
  \(\displaystyle {35}\),
  \(\displaystyle {40}\),
},
y grid style={darkgray176},
ylabel={\(\displaystyle \mathcal{J}\)},
ymajorgrids,
ymin=0.000, ymax=200.581064323375,
ymode=log,
ytick style={color=black},
ytick={0.1,1,10,100,1000,10000},
yticklabels={
  \(\displaystyle {10^{-1}}\),
  \(\displaystyle {10^{0}}\),
  \(\displaystyle {10^{1}}\),
  \(\displaystyle {10^{2}}\),
  \(\displaystyle {10^{3}}\),
  \(\displaystyle {10^{4}}\)
}
]
\addplot [draw=green, fill=green, mark=square*, only marks]
table{%
x  y
4 93.9293
5 93.9293
8 93.9293
10 93.9293
12 70.8661
15 36.8155
17 42.7626
20 20.9097
25 1.9524
30 1.4339
35 1.2369
40 1.2339
};
\addlegendentry{POD-G}
% \addplot [draw=brown, fill=brown, mark options={rotate=180}, mark=triangle*, only marks]
% table{%
% x  y
% 1 107.93815
% 2 167.68704
% 3 148.0767
% 4 148.50524
% 5 211.13952
% 6 242.84521
% 7 178.25943
% 8 191.88937
% 9 175.48529
% 10 251.02238
% 11 157.31428
% 12 174.71115
% 13 141.60819
% 14 151.94018
% 15 119.15534
% 16 101.19005
% 17 100.45222
% 18 77.96346
% 19 96.89764
% 20 59.88348
% 25 15.59005
% 30 10.03992
% 35 15.73395
% 40 20.47886
% 45 11.84484
% 50 29.11487
% 60 6.76843
% 70 6.31444
% 80 24.92309
% 90 32.56686
% 100 4.19868
% 110 2.77336
% 120 9.17072
% 130 2.63423
% 140 19.44395
% };
% \addlegendentry{FRTO POD-G}
\addplot [draw=red, fill=red, mark=asterisk, only marks]
table{%
x  y
4 8.7990
5 1.9371
8 2.8798
10 2.5581
12 1.3056
15 1.5015
17 1.3766
20 1.6950
25 1.9816
30 2.9485
35 1.8333
40 3.0492
};
\addlegendentry{sPOD-G}
% \addplot [draw=blue, fill=blue, mark=*, only marks]
% table{%
% x  y
% 1 83.4267
% 2 48.8674
% 3 37.09822
% 4 35.23015
% 5 21.28047
% 6 16.48147
% 7 11.17033
% 8 9.72923
% 9 10.67685
% 10 7.92139
% 11 6.42046
% 12 7.82543
% 13 14.55814
% 14 9.15829
% 15 5.12031
% 16 8.31066
% 17 4.74758
% 18 4.7713
% 19 4.65177
% 20 6.65589
% 21 6.54433
% 22 6.62408
% 23 5.32576
% 24 6.57799
% 25 4.88435
% 26 3.5656
% 27 3.74806
% 28 5.10514
% 29 9.37614
% 30 7.02525
% };
% \addlegendentry{FRTO sPOD-G}
\addplot[mark=none, brown, samples=150] coordinates {(-5.95,1.2290)(146.95,1.2290)} ; 
\addlegendentry{FOM}
\end{axis}

\end{tikzpicture}

%% file: ModesVsCost_2.tex
% This file was created with tikzplotlib v0.10.1.
\begin{tikzpicture}

\definecolor{brown}{RGB}{165,42,42}
\definecolor{darkgray176}{RGB}{176,176,176}
\definecolor{green}{RGB}{0,128,0}
\definecolor{lightgray204}{RGB}{204,204,204}

\begin{axis}[
height=\figureheight,
legend cell align={left},
legend style={fill opacity=0.8, draw opacity=1, text opacity=1, draw=lightgray204},
log basis y={10},
tick align=outside,
tick pos=left,
width=\figurewidth,
x grid style={darkgray176},
xlabel={Truncated modes},
xmajorgrids,
xmin=-2.5, xmax=142.5,
xtick style={color=black},
xtick={0,20,40,60,80,100,120,140},
xticklabels={
  \(\displaystyle {0}\),
  % \(\displaystyle {10}\),
  \(\displaystyle {20}\),
  % \(\displaystyle {30}\),
  \(\displaystyle {40}\),
  % \(\displaystyle {50}\),
  \(\displaystyle {60}\),
  % \(\displaystyle {70}\),
  \(\displaystyle {80}\),
  % \(\displaystyle {90}\),
  \(\displaystyle {100}\),
  % \(\displaystyle {110}\),
  \(\displaystyle {120}\),
  % \(\displaystyle {130}\),
  \(\displaystyle {140}\),
},
y grid style={darkgray176},
ylabel={\(\displaystyle \mathcal{J}\)},
ymajorgrids,
ymin=0.000, ymax=200.581064323375,
ymode=log,
ytick style={color=black},
ytick={0.1,1,10,100,1000,10000},
yticklabels={
  \(\displaystyle {10^{-1}}\),
  \(\displaystyle {10^{0}}\),
  \(\displaystyle {10^{1}}\),
  \(\displaystyle {10^{2}}\),
  \(\displaystyle {10^{3}}\),
  \(\displaystyle {10^{4}}\)
}
]
\addplot [draw=green, fill=green, mark=square*, only marks]
table{%
x  y
4 20.9881
5 29.2356
8 29.2356
10 29.2356
12 29.2356
15 29.2356
17 29.2356
20 29.2356
25 29.2356
30 29.2356
35 29.2356
40 25.6912
50 15.8892
60 12.0874
70 8.6945
80 6.3147
90 3.7464
100 3.2270
110 3.0158
120 3.0089
130 3.0125
140 3.0105
};
\addlegendentry{POD-G}

\addplot [draw=red, fill=red, mark=asterisk, only marks]
table{%
x  y
4 6.8227
5 7.8104
8 6.6770
10 7.6506
12 8.7690
15 6.5749
17 9.1331 
20 2.9827
25 3.0494
30 3.7294
35 3.7107
40 6.0582
};
\addlegendentry{sPOD-G}
\addplot[mark=none, brown, samples=150] coordinates {(-5.95,2.9831)(146.95,2.9831)} ; 
\addlegendentry{FOM}
\end{axis}

\end{tikzpicture}

%% file: ModesVsCost_3.tex
% This file was created with tikzplotlib v0.10.1.
\begin{tikzpicture}

\definecolor{brown}{RGB}{165,42,42}
\definecolor{darkgray176}{RGB}{176,176,176}
\definecolor{green}{RGB}{0,128,0}
\definecolor{lightgray204}{RGB}{204,204,204}

\begin{axis}[
height=\figureheight,
legend cell align={left},
legend style={fill opacity=0.8, draw opacity=1, text opacity=1, draw=lightgray204},
log basis y={10},
tick align=outside,
tick pos=left,
width=\figurewidth,
x grid style={darkgray176},
xlabel={Truncated modes},
xmajorgrids,
xmin=-2.5, xmax=142.5,
xtick style={color=black},
xtick={0,20,40,60,80,100,120,140},
xticklabels={
  \(\displaystyle {0}\),
  % \(\displaystyle {10}\),
  \(\displaystyle {20}\),
  % \(\displaystyle {30}\),
  \(\displaystyle {40}\),
  % \(\displaystyle {50}\),
  \(\displaystyle {60}\),
  % \(\displaystyle {70}\),
  \(\displaystyle {80}\),
  % \(\displaystyle {90}\),
  \(\displaystyle {100}\),
  % \(\displaystyle {110}\),
  \(\displaystyle {120}\),
  % \(\displaystyle {130}\),
  \(\displaystyle {140}\),
},
y grid style={darkgray176},
ylabel={\(\displaystyle \mathcal{J}\)},
ymajorgrids,
ymin=0.000, ymax=200.581064323375,
ymode=log,
ytick style={color=black},
ytick={0.1,1,10,100,1000,10000},
yticklabels={
  \(\displaystyle {10^{-1}}\),
  \(\displaystyle {10^{0}}\),
  \(\displaystyle {10^{1}}\),
  \(\displaystyle {10^{2}}\),
  \(\displaystyle {10^{3}}\),
  \(\displaystyle {10^{4}}\)
}
]
\addplot [draw=green, fill=green, mark=square*, only marks]
table{%
x  y
4 265.6590
5 68.6849
8 10.7046
10 11.2481
12 11.2481
15 11.2481
17 11.2481
20 11.2481
25 11.2481
30 11.2481
35 11.2481
40 11.2481
50 11.2481
60 11.2481
70 10.6283
80 7.4471
90 3.1277
100 2.3499
110 1.9720
120 1.9584
130 2.0098
140 2.0740
};
\addlegendentry{POD-G}

\addplot [draw=red, fill=red, mark=asterisk, only marks]
table{%
x  y
4 2.1239
5 1.9504
8 1.9322
10 1.9688
12 1.8976
15 1.9762
17 2.4272 
20 2.4930
25 1.9035
30 1.9405
35 2.0176
40 2.0857
};
\addlegendentry{sPOD-G}
\addplot[mark=none, brown, samples=150] coordinates {(-5.95,1.8977)(146.95,1.8977)} ; 
\addlegendentry{FOM}
\end{axis}

\end{tikzpicture}

%% file: Example_1_snapshots.tex
% This file was created with tikzplotlib v0.10.1.
\begin{tikzpicture}

\definecolor{darkgray176}{RGB}{176,176,176}

\begin{groupplot}[group style={group size=4 by 3, horizontal sep=1.30cm, vertical sep=0.5cm}]
\nextgroupplot[
colorbar,
colorbar style={ytick={-0.6,-0.4,-0.2,0,0.2,0.4,0.6},yticklabels={
  \(\displaystyle {\ensuremath{-}0.6}\),
  \(\displaystyle {\ensuremath{-}0.4}\),
  \(\displaystyle {\ensuremath{-}0.2}\),
  \(\displaystyle {0.0}\),
  \(\displaystyle {0.2}\),
  \(\displaystyle {0.4}\),
  \(\displaystyle {0.6}\)
},ylabel={}, width=0.06*\pgfkeysvalueof{/pgfplots/parent axis width}, ticklabel style={font=\tiny}},
colormap={mymap}{[1pt]
  rgb(0pt)=(1,1,0.8);
  rgb(1pt)=(1,0.929411764705882,0.627450980392157);
  rgb(2pt)=(0.996078431372549,0.850980392156863,0.462745098039216);
  rgb(3pt)=(0.996078431372549,0.698039215686274,0.298039215686275);
  rgb(4pt)=(0.992156862745098,0.552941176470588,0.235294117647059);
  rgb(5pt)=(0.988235294117647,0.305882352941176,0.164705882352941);
  rgb(6pt)=(0.890196078431372,0.101960784313725,0.109803921568627);
  rgb(7pt)=(0.741176470588235,0,0.149019607843137);
  rgb(8pt)=(0.501960784313725,0,0.149019607843137)
},
height=\figureheight,
hide x axis,
hide y axis,
point meta max=0.437558488313281,
point meta min=-0.434043701126128,
tick align=outside,
tick pos=left,
align=center,
title={\tiny{FOM}},
width=\figurewidth,
x grid style={darkgray176},
xmin=0, xmax=3200,
xtick style={color=black},
xtick={0,1000,2000,3000,4000},
xticklabels={
  \(\displaystyle {0}\),
  \(\displaystyle {1000}\),
  \(\displaystyle {2000}\),
  \(\displaystyle {3000}\),
  \(\displaystyle {4000}\)
},
y grid style={darkgray176},
ymin=0, ymax=3360,
ytick style={color=black},
ytick={0,500,1000,1500,2000,2500,3000,3500},
yticklabels={
  \(\displaystyle {0}\),
  \(\displaystyle {500}\),
  \(\displaystyle {1000}\),
  \(\displaystyle {1500}\),
  \(\displaystyle {2000}\),
  \(\displaystyle {2500}\),
  \(\displaystyle {3000}\),
  \(\displaystyle {3500}\)
}
]
\addplot graphics [includegraphics cmd=\pgfimage,xmin=0, xmax=3200, ymin=0, ymax=3360] {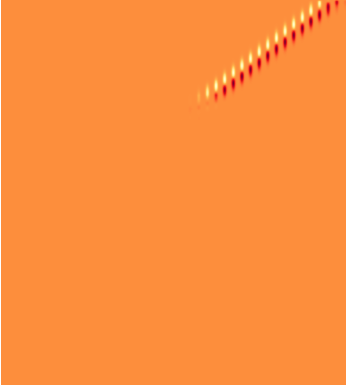};

\nextgroupplot[
colorbar,
colorbar style={ytick={-0.6,-0.4,-0.2,0,0.2,0.4,0.6},yticklabels={
  \(\displaystyle {\ensuremath{-}0.6}\),
  \(\displaystyle {\ensuremath{-}0.4}\),
  \(\displaystyle {\ensuremath{-}0.2}\),
  \(\displaystyle {0.0}\),
  \(\displaystyle {0.2}\),
  \(\displaystyle {0.4}\),
  \(\displaystyle {0.6}\)
},ylabel={}, width=0.06*\pgfkeysvalueof{/pgfplots/parent axis width}, ticklabel style={font=\tiny}},
colormap={mymap}{[1pt]
  rgb(0pt)=(1,1,0.8);
  rgb(1pt)=(1,0.929411764705882,0.627450980392157);
  rgb(2pt)=(0.996078431372549,0.850980392156863,0.462745098039216);
  rgb(3pt)=(0.996078431372549,0.698039215686274,0.298039215686275);
  rgb(4pt)=(0.992156862745098,0.552941176470588,0.235294117647059);
  rgb(5pt)=(0.988235294117647,0.305882352941176,0.164705882352941);
  rgb(6pt)=(0.890196078431372,0.101960784313725,0.109803921568627);
  rgb(7pt)=(0.741176470588235,0,0.149019607843137);
  rgb(8pt)=(0.501960784313725,0,0.149019607843137)
},
height=\figureheight,
hide x axis,
hide y axis,
point meta max=0.437558488313281,
point meta min=-0.434043701126128,
tick align=outside,
tick pos=left,
align=center,
title={\tiny{sPOD-G}\\ \tiny{($r=12$)}},
width=\figurewidth,
x grid style={darkgray176},
xmin=0, xmax=3200,
xtick style={color=black},
xtick={0,1000,2000,3000,4000},
xticklabels={
  \(\displaystyle {0}\),
  \(\displaystyle {1000}\),
  \(\displaystyle {2000}\),
  \(\displaystyle {3000}\),
  \(\displaystyle {4000}\)
},
y grid style={darkgray176},
ymin=0, ymax=3360,
ytick style={color=black},
ytick={0,500,1000,1500,2000,2500,3000,3500},
yticklabels={
  \(\displaystyle {0}\),
  \(\displaystyle {500}\),
  \(\displaystyle {1000}\),
  \(\displaystyle {1500}\),
  \(\displaystyle {2000}\),
  \(\displaystyle {2500}\),
  \(\displaystyle {3000}\),
  \(\displaystyle {3500}\)
}
]
\addplot graphics [includegraphics cmd=\pgfimage,xmin=0, xmax=3200, ymin=0, ymax=3360] {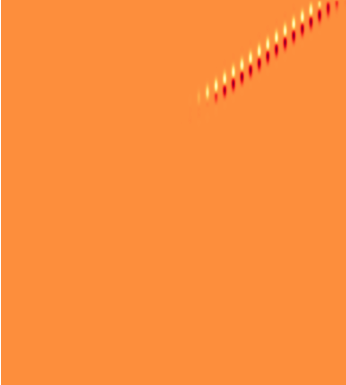};

\nextgroupplot[
colorbar,
colorbar style={ytick={-0.6,-0.4,-0.2,0,0.2,0.4,0.6},yticklabels={
  \(\displaystyle {\ensuremath{-}0.6}\),
  \(\displaystyle {\ensuremath{-}0.4}\),
  \(\displaystyle {\ensuremath{-}0.2}\),
  \(\displaystyle {0.0}\),
  \(\displaystyle {0.2}\),
  \(\displaystyle {0.4}\),
  \(\displaystyle {0.6}\)
},ylabel={}, width=0.06*\pgfkeysvalueof{/pgfplots/parent axis width}, ticklabel style={font=\tiny}},
colormap={mymap}{[1pt]
  rgb(0pt)=(1,1,0.8);
  rgb(1pt)=(1,0.929411764705882,0.627450980392157);
  rgb(2pt)=(0.996078431372549,0.850980392156863,0.462745098039216);
  rgb(3pt)=(0.996078431372549,0.698039215686274,0.298039215686275);
  rgb(4pt)=(0.992156862745098,0.552941176470588,0.235294117647059);
  rgb(5pt)=(0.988235294117647,0.305882352941176,0.164705882352941);
  rgb(6pt)=(0.890196078431372,0.101960784313725,0.109803921568627);
  rgb(7pt)=(0.741176470588235,0,0.149019607843137);
  rgb(8pt)=(0.501960784313725,0,0.149019607843137)
},
height=\figureheight,
hide x axis,
hide y axis,
point meta max=0.437558488313281,
point meta min=-0.434043701126128,
tick align=outside,
tick pos=left,
align=center,
title={\tiny{POD-G}\\ \tiny{($\mathrm{p}=40$)}},
width=\figurewidth,
x grid style={darkgray176},
xmin=0, xmax=3200,
xtick style={color=black},
xtick={0,1000,2000,3000,4000},
xticklabels={
  \(\displaystyle {0}\),
  \(\displaystyle {1000}\),
  \(\displaystyle {2000}\),
  \(\displaystyle {3000}\),
  \(\displaystyle {4000}\)
},
y grid style={darkgray176},
ymin=0, ymax=3360,
ytick style={color=black},
ytick={0,500,1000,1500,2000,2500,3000,3500},
yticklabels={
  \(\displaystyle {0}\),
  \(\displaystyle {500}\),
  \(\displaystyle {1000}\),
  \(\displaystyle {1500}\),
  \(\displaystyle {2000}\),
  \(\displaystyle {2500}\),
  \(\displaystyle {3000}\),
  \(\displaystyle {3500}\)
}
]
\addplot graphics [includegraphics cmd=\pgfimage,xmin=0, xmax=3200, ymin=0, ymax=3360] {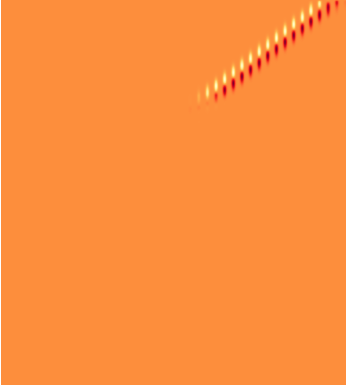};

\nextgroupplot[
colorbar,
colorbar style={ytick={-0.6,-0.4,-0.2,0,0.2,0.4,0.6},yticklabels={
  \(\displaystyle {\ensuremath{-}0.6}\),
  \(\displaystyle {\ensuremath{-}0.4}\),
  \(\displaystyle {\ensuremath{-}0.2}\),
  \(\displaystyle {0.0}\),
  \(\displaystyle {0.2}\),
  \(\displaystyle {0.4}\),
  \(\displaystyle {0.6}\)
},ylabel={}, width=0.06*\pgfkeysvalueof{/pgfplots/parent axis width}, ticklabel style={font=\tiny}},
colormap={mymap}{[1pt]
  rgb(0pt)=(1,1,0.8);
  rgb(1pt)=(1,0.929411764705882,0.627450980392157);
  rgb(2pt)=(0.996078431372549,0.850980392156863,0.462745098039216);
  rgb(3pt)=(0.996078431372549,0.698039215686274,0.298039215686275);
  rgb(4pt)=(0.992156862745098,0.552941176470588,0.235294117647059);
  rgb(5pt)=(0.988235294117647,0.305882352941176,0.164705882352941);
  rgb(6pt)=(0.890196078431372,0.101960784313725,0.109803921568627);
  rgb(7pt)=(0.741176470588235,0,0.149019607843137);
  rgb(8pt)=(0.501960784313725,0,0.149019607843137)
},
height=\figureheight,
hide x axis,
hide y axis,
point meta max=0.437558488313281,
point meta min=-0.434043701126128,
tick align=outside,
tick pos=left,
align=center,
title={\tiny{POD-G} \\ \tiny{($\mathrm{p}=12$)}},
width=\figurewidth,
x grid style={darkgray176},
xmin=0, xmax=3200,
xtick style={color=black},
xtick={0,1000,2000,3000,4000},
xticklabels={
  \(\displaystyle {0}\),
  \(\displaystyle {1000}\),
  \(\displaystyle {2000}\),
  \(\displaystyle {3000}\),
  \(\displaystyle {4000}\)
},
y grid style={darkgray176},
ymin=0, ymax=3360,
ytick style={color=black},
ytick={0,500,1000,1500,2000,2500,3000,3500},
yticklabels={
  \(\displaystyle {0}\),
  \(\displaystyle {500}\),
  \(\displaystyle {1000}\),
  \(\displaystyle {1500}\),
  \(\displaystyle {2000}\),
  \(\displaystyle {2500}\),
  \(\displaystyle {3000}\),
  \(\displaystyle {3500}\)
}
]
\addplot graphics [includegraphics cmd=\pgfimage,xmin=0, xmax=3200, ymin=0, ymax=3360] {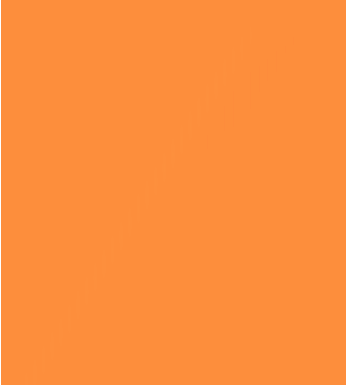};

%%%%%%%%%%%%%%%%%%%%%%%%%%%%%%%%%%%%%%%%%%%%%%%%%%%%%%%%%%%%%%%%%%%%%%%%%%%%%

\nextgroupplot[
colorbar,
colorbar style={ytick={-0.01,-0.005,0,0.005,0.01},yticklabels={
  \(\displaystyle {\ensuremath{-}0.010}\),
  \(\displaystyle {\ensuremath{-}0.005}\),
  \(\displaystyle {0.000}\),
  \(\displaystyle {0.005}\),
  \(\displaystyle {0.010}\)
},ylabel={}, width=0.06*\pgfkeysvalueof{/pgfplots/parent axis width}, ticklabel style={font=\tiny}},
colormap={mymap}{[1pt]
  rgb(0pt)=(1,1,0.8);
  rgb(1pt)=(1,0.929411764705882,0.627450980392157);
  rgb(2pt)=(0.996078431372549,0.850980392156863,0.462745098039216);
  rgb(3pt)=(0.996078431372549,0.698039215686274,0.298039215686275);
  rgb(4pt)=(0.992156862745098,0.552941176470588,0.235294117647059);
  rgb(5pt)=(0.988235294117647,0.305882352941176,0.164705882352941);
  rgb(6pt)=(0.890196078431372,0.101960784313725,0.109803921568627);
  rgb(7pt)=(0.741176470588235,0,0.149019607843137);
  rgb(8pt)=(0.501960784313725,0,0.149019607843137)
},
height=\figureheight,
hide x axis,
hide y axis,
point meta max=0.00894688826680106,
point meta min=-0.00867484563333541,
tick align=outside,
tick pos=left,
width=\figurewidth,
x grid style={darkgray176},
xmin=0, xmax=3200,
xtick style={color=black},
xtick={0,1000,2000,3000,4000},
xticklabels={
  \(\displaystyle {0}\),
  \(\displaystyle {1000}\),
  \(\displaystyle {2000}\),
  \(\displaystyle {3000}\),
  \(\displaystyle {4000}\)
},
y grid style={darkgray176},
ymin=0, ymax=3360,
ytick style={color=black},
ytick={0,500,1000,1500,2000,2500,3000,3500},
yticklabels={
  \(\displaystyle {0}\),
  \(\displaystyle {500}\),
  \(\displaystyle {1000}\),
  \(\displaystyle {1500}\),
  \(\displaystyle {2000}\),
  \(\displaystyle {2500}\),
  \(\displaystyle {3000}\),
  \(\displaystyle {3500}\)
}
]
\addplot graphics [includegraphics cmd=\pgfimage,xmin=0, xmax=3200, ymin=0, ymax=3360] {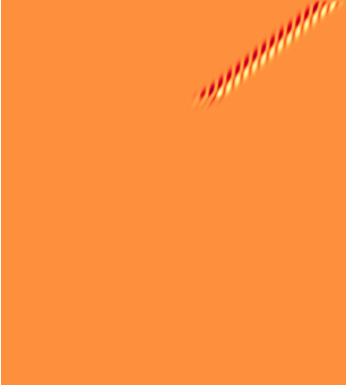};

\nextgroupplot[
colorbar,
colorbar style={ytick={-0.01,-0.005,0,0.005,0.01},yticklabels={
  \(\displaystyle {\ensuremath{-}0.010}\),
  \(\displaystyle {\ensuremath{-}0.005}\),
  \(\displaystyle {0.000}\),
  \(\displaystyle {0.005}\),
  \(\displaystyle {0.010}\)
},ylabel={}, width=0.06*\pgfkeysvalueof{/pgfplots/parent axis width}, ticklabel style={font=\tiny}},
colormap={mymap}{[1pt]
  rgb(0pt)=(1,1,0.8);
  rgb(1pt)=(1,0.929411764705882,0.627450980392157);
  rgb(2pt)=(0.996078431372549,0.850980392156863,0.462745098039216);
  rgb(3pt)=(0.996078431372549,0.698039215686274,0.298039215686275);
  rgb(4pt)=(0.992156862745098,0.552941176470588,0.235294117647059);
  rgb(5pt)=(0.988235294117647,0.305882352941176,0.164705882352941);
  rgb(6pt)=(0.890196078431372,0.101960784313725,0.109803921568627);
  rgb(7pt)=(0.741176470588235,0,0.149019607843137);
  rgb(8pt)=(0.501960784313725,0,0.149019607843137)
},
height=\figureheight,
hide x axis,
hide y axis,
point meta max=0.00894688826680106,
point meta min=-0.00867484563333541,
tick align=outside,
tick pos=left,
width=\figurewidth,
x grid style={darkgray176},
xmin=0, xmax=3200,
xtick style={color=black},
xtick={0,1000,2000,3000,4000},
xticklabels={
  \(\displaystyle {0}\),
  \(\displaystyle {1000}\),
  \(\displaystyle {2000}\),
  \(\displaystyle {3000}\),
  \(\displaystyle {4000}\)
},
y grid style={darkgray176},
ymin=0, ymax=3360,
ytick style={color=black},
ytick={0,500,1000,1500,2000,2500,3000,3500},
yticklabels={
  \(\displaystyle {0}\),
  \(\displaystyle {500}\),
  \(\displaystyle {1000}\),
  \(\displaystyle {1500}\),
  \(\displaystyle {2000}\),
  \(\displaystyle {2500}\),
  \(\displaystyle {3000}\),
  \(\displaystyle {3500}\)
}
]
\addplot graphics [includegraphics cmd=\pgfimage,xmin=0, xmax=3200, ymin=0, ymax=3360] {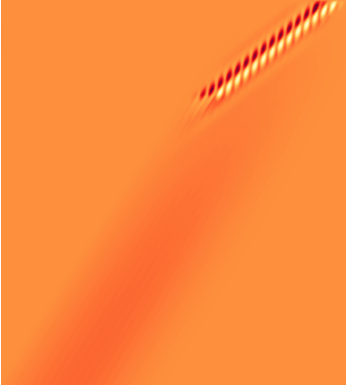};

\nextgroupplot[
colorbar,
colorbar style={ytick={-0.01,-0.005,0,0.005,0.01},yticklabels={
  \(\displaystyle {\ensuremath{-}0.010}\),
  \(\displaystyle {\ensuremath{-}0.005}\),
  \(\displaystyle {0.000}\),
  \(\displaystyle {0.005}\),
  \(\displaystyle {0.010}\)
},ylabel={}, width=0.06*\pgfkeysvalueof{/pgfplots/parent axis width}, ticklabel style={font=\tiny}},
colormap={mymap}{[1pt]
  rgb(0pt)=(1,1,0.8);
  rgb(1pt)=(1,0.929411764705882,0.627450980392157);
  rgb(2pt)=(0.996078431372549,0.850980392156863,0.462745098039216);
  rgb(3pt)=(0.996078431372549,0.698039215686274,0.298039215686275);
  rgb(4pt)=(0.992156862745098,0.552941176470588,0.235294117647059);
  rgb(5pt)=(0.988235294117647,0.305882352941176,0.164705882352941);
  rgb(6pt)=(0.890196078431372,0.101960784313725,0.109803921568627);
  rgb(7pt)=(0.741176470588235,0,0.149019607843137);
  rgb(8pt)=(0.501960784313725,0,0.149019607843137)
},
height=\figureheight,
hide x axis,
hide y axis,
point meta max=0.00894688826680106,
point meta min=-0.00867484563333541,
tick align=outside,
tick pos=left,
width=\figurewidth,
x grid style={darkgray176},
xmin=0, xmax=3200,
xtick style={color=black},
xtick={0,1000,2000,3000,4000},
xticklabels={
  \(\displaystyle {0}\),
  \(\displaystyle {1000}\),
  \(\displaystyle {2000}\),
  \(\displaystyle {3000}\),
  \(\displaystyle {4000}\)
},
y grid style={darkgray176},
ymin=0, ymax=3360,
ytick style={color=black},
ytick={0,500,1000,1500,2000,2500,3000,3500},
yticklabels={
  \(\displaystyle {0}\),
  \(\displaystyle {500}\),
  \(\displaystyle {1000}\),
  \(\displaystyle {1500}\),
  \(\displaystyle {2000}\),
  \(\displaystyle {2500}\),
  \(\displaystyle {3000}\),
  \(\displaystyle {3500}\)
}
]
\addplot graphics [includegraphics cmd=\pgfimage,xmin=0, xmax=3200, ymin=0, ymax=3360] {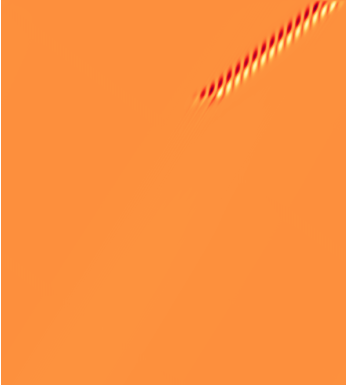};

\nextgroupplot[
colorbar,
colorbar style={ytick={-0.01,-0.005,0,0.005,0.01},yticklabels={
  \(\displaystyle {\ensuremath{-}0.010}\),
  \(\displaystyle {\ensuremath{-}0.005}\),
  \(\displaystyle {0.000}\),
  \(\displaystyle {0.005}\),
  \(\displaystyle {0.010}\)
},ylabel={}, width=0.06*\pgfkeysvalueof{/pgfplots/parent axis width}, ticklabel style={font=\tiny}},
colormap={mymap}{[1pt]
  rgb(0pt)=(1,1,0.8);
  rgb(1pt)=(1,0.929411764705882,0.627450980392157);
  rgb(2pt)=(0.996078431372549,0.850980392156863,0.462745098039216);
  rgb(3pt)=(0.996078431372549,0.698039215686274,0.298039215686275);
  rgb(4pt)=(0.992156862745098,0.552941176470588,0.235294117647059);
  rgb(5pt)=(0.988235294117647,0.305882352941176,0.164705882352941);
  rgb(6pt)=(0.890196078431372,0.101960784313725,0.109803921568627);
  rgb(7pt)=(0.741176470588235,0,0.149019607843137);
  rgb(8pt)=(0.501960784313725,0,0.149019607843137)
},
height=\figureheight,
hide x axis,
hide y axis,
point meta max=0.00894688826680106,
point meta min=-0.00867484563333541,
tick align=outside,
tick pos=left,
width=\figurewidth,
x grid style={darkgray176},
xmin=0, xmax=3200,
xtick style={color=black},
xtick={0,1000,2000,3000,4000},
xticklabels={
  \(\displaystyle {0}\),
  \(\displaystyle {1000}\),
  \(\displaystyle {2000}\),
  \(\displaystyle {3000}\),
  \(\displaystyle {4000}\)
},
y grid style={darkgray176},
ymin=0, ymax=3360,
ytick style={color=black},
ytick={0,500,1000,1500,2000,2500,3000,3500},
yticklabels={
  \(\displaystyle {0}\),
  \(\displaystyle {500}\),
  \(\displaystyle {1000}\),
  \(\displaystyle {1500}\),
  \(\displaystyle {2000}\),
  \(\displaystyle {2500}\),
  \(\displaystyle {3000}\),
  \(\displaystyle {3500}\)
}
]
\addplot graphics [includegraphics cmd=\pgfimage,xmin=0, xmax=3200, ymin=0, ymax=3360] {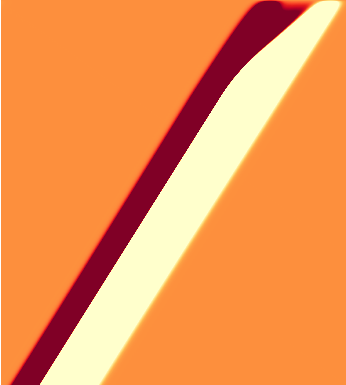};

%%%%%%%%%%%%%%%%%%%%%%%%%%%%%%%%%%%%%%%%%%%%%%%%%%%%%%%%%%%%%%%%%%%%%%%%%%%%%%%%%%%%%%

\nextgroupplot[
colorbar,
colorbar style={ytick={-0.2,0,0.2,0.4,0.6,0.8,1,1.2},yticklabels={
  \(\displaystyle {\ensuremath{-}0.2}\),
  \(\displaystyle {0.0}\),
  \(\displaystyle {0.2}\),
  \(\displaystyle {0.4}\),
  \(\displaystyle {0.6}\),
  \(\displaystyle {0.8}\),
  \(\displaystyle {1.0}\),
  \(\displaystyle {1.2}\)
},ylabel={}, width=0.06*\pgfkeysvalueof{/pgfplots/parent axis width}, ticklabel style={font=\tiny}},
colormap={mymap}{[1pt]
  rgb(0pt)=(1,1,0.8);
  rgb(1pt)=(1,0.929411764705882,0.627450980392157);
  rgb(2pt)=(0.996078431372549,0.850980392156863,0.462745098039216);
  rgb(3pt)=(0.996078431372549,0.698039215686274,0.298039215686275);
  rgb(4pt)=(0.992156862745098,0.552941176470588,0.235294117647059);
  rgb(5pt)=(0.988235294117647,0.305882352941176,0.164705882352941);
  rgb(6pt)=(0.890196078431372,0.101960784313725,0.109803921568627);
  rgb(7pt)=(0.741176470588235,0,0.149019607843137);
  rgb(8pt)=(0.501960784313725,0,0.149019607843137)
},
height=\figureheight,
hide x axis,
hide y axis,
point meta max=1.00981603847325,
point meta min=-0.000458797682120604,
tick align=outside,
tick pos=left,
width=\figurewidth,
x grid style={darkgray176},
xmin=0, xmax=3200,
xtick style={color=black},
xtick={0,1000,2000,3000,4000},
xticklabels={
  \(\displaystyle {0}\),
  \(\displaystyle {1000}\),
  \(\displaystyle {2000}\),
  \(\displaystyle {3000}\),
  \(\displaystyle {4000}\)
},
y grid style={darkgray176},
ymin=0, ymax=3360,
ytick style={color=black},
ytick={0,500,1000,1500,2000,2500,3000,3500},
yticklabels={
  \(\displaystyle {0}\),
  \(\displaystyle {500}\),
  \(\displaystyle {1000}\),
  \(\displaystyle {1500}\),
  \(\displaystyle {2000}\),
  \(\displaystyle {2500}\),
  \(\displaystyle {3000}\),
  \(\displaystyle {3500}\)
}
]
\addplot graphics [includegraphics cmd=\pgfimage,xmin=0, xmax=3200, ymin=0, ymax=3360] {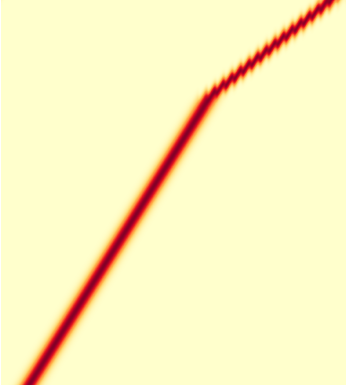};

\nextgroupplot[
colorbar,
colorbar style={ytick={-0.2,0,0.2,0.4,0.6,0.8,1,1.2},yticklabels={
  \(\displaystyle {\ensuremath{-}0.2}\),
  \(\displaystyle {0.0}\),
  \(\displaystyle {0.2}\),
  \(\displaystyle {0.4}\),
  \(\displaystyle {0.6}\),
  \(\displaystyle {0.8}\),
  \(\displaystyle {1.0}\),
  \(\displaystyle {1.2}\)
},ylabel={}, width=0.06*\pgfkeysvalueof{/pgfplots/parent axis width}, ticklabel style={font=\tiny}},
colormap={mymap}{[1pt]
  rgb(0pt)=(1,1,0.8);
  rgb(1pt)=(1,0.929411764705882,0.627450980392157);
  rgb(2pt)=(0.996078431372549,0.850980392156863,0.462745098039216);
  rgb(3pt)=(0.996078431372549,0.698039215686274,0.298039215686275);
  rgb(4pt)=(0.992156862745098,0.552941176470588,0.235294117647059);
  rgb(5pt)=(0.988235294117647,0.305882352941176,0.164705882352941);
  rgb(6pt)=(0.890196078431372,0.101960784313725,0.109803921568627);
  rgb(7pt)=(0.741176470588235,0,0.149019607843137);
  rgb(8pt)=(0.501960784313725,0,0.149019607843137)
},
height=\figureheight,
hide x axis,
hide y axis,
point meta max=1.00981603847325,
point meta min=-0.000458797682120604,
tick align=outside,
tick pos=left,
width=\figurewidth,
x grid style={darkgray176},
xmin=0, xmax=3200,
xtick style={color=black},
xtick={0,1000,2000,3000,4000},
xticklabels={
  \(\displaystyle {0}\),
  \(\displaystyle {1000}\),
  \(\displaystyle {2000}\),
  \(\displaystyle {3000}\),
  \(\displaystyle {4000}\)
},
y grid style={darkgray176},
ymin=0, ymax=3360,
ytick style={color=black},
ytick={0,500,1000,1500,2000,2500,3000,3500},
yticklabels={
  \(\displaystyle {0}\),
  \(\displaystyle {500}\),
  \(\displaystyle {1000}\),
  \(\displaystyle {1500}\),
  \(\displaystyle {2000}\),
  \(\displaystyle {2500}\),
  \(\displaystyle {3000}\),
  \(\displaystyle {3500}\)
}
]
\addplot graphics [includegraphics cmd=\pgfimage,xmin=0, xmax=3200, ymin=0, ymax=3360] {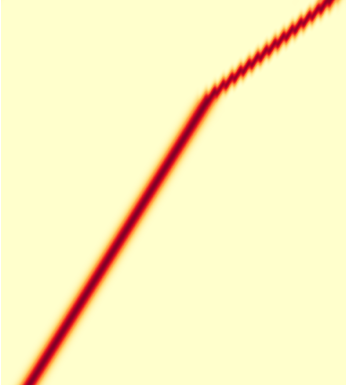};

\nextgroupplot[
colorbar,
colorbar style={ytick={-0.2,0,0.2,0.4,0.6,0.8,1,1.2},yticklabels={
  \(\displaystyle {\ensuremath{-}0.2}\),
  \(\displaystyle {0.0}\),
  \(\displaystyle {0.2}\),
  \(\displaystyle {0.4}\),
  \(\displaystyle {0.6}\),
  \(\displaystyle {0.8}\),
  \(\displaystyle {1.0}\),
  \(\displaystyle {1.2}\)
},ylabel={}, width=0.06*\pgfkeysvalueof{/pgfplots/parent axis width}, ticklabel style={font=\tiny}},
colormap={mymap}{[1pt]
  rgb(0pt)=(1,1,0.8);
  rgb(1pt)=(1,0.929411764705882,0.627450980392157);
  rgb(2pt)=(0.996078431372549,0.850980392156863,0.462745098039216);
  rgb(3pt)=(0.996078431372549,0.698039215686274,0.298039215686275);
  rgb(4pt)=(0.992156862745098,0.552941176470588,0.235294117647059);
  rgb(5pt)=(0.988235294117647,0.305882352941176,0.164705882352941);
  rgb(6pt)=(0.890196078431372,0.101960784313725,0.109803921568627);
  rgb(7pt)=(0.741176470588235,0,0.149019607843137);
  rgb(8pt)=(0.501960784313725,0,0.149019607843137)
},
height=\figureheight,
hide x axis,
hide y axis,
point meta max=1.00981603847325,
point meta min=-0.000458797682120604,
tick align=outside,
tick pos=left,
width=\figurewidth,
x grid style={darkgray176},
xmin=0, xmax=3200,
xtick style={color=black},
xtick={0,1000,2000,3000,4000},
xticklabels={
  \(\displaystyle {0}\),
  \(\displaystyle {1000}\),
  \(\displaystyle {2000}\),
  \(\displaystyle {3000}\),
  \(\displaystyle {4000}\)
},
y grid style={darkgray176},
ymin=0, ymax=3360,
ytick style={color=black},
ytick={0,500,1000,1500,2000,2500,3000,3500},
yticklabels={
  \(\displaystyle {0}\),
  \(\displaystyle {500}\),
  \(\displaystyle {1000}\),
  \(\displaystyle {1500}\),
  \(\displaystyle {2000}\),
  \(\displaystyle {2500}\),
  \(\displaystyle {3000}\),
  \(\displaystyle {3500}\)
}
]
\addplot graphics [includegraphics cmd=\pgfimage,xmin=0, xmax=3200, ymin=0, ymax=3360] {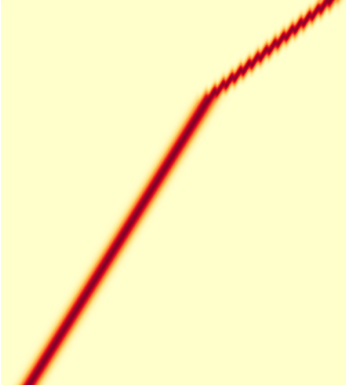};

\nextgroupplot[
colorbar,
colorbar style={ytick={-0.2,0,0.2,0.4,0.6,0.8,1,1.2},yticklabels={
  \(\displaystyle {\ensuremath{-}0.2}\),
  \(\displaystyle {0.0}\),
  \(\displaystyle {0.2}\),
  \(\displaystyle {0.4}\),
  \(\displaystyle {0.6}\),
  \(\displaystyle {0.8}\),
  \(\displaystyle {1.0}\),
  \(\displaystyle {1.2}\)
},ylabel={}, width=0.06*\pgfkeysvalueof{/pgfplots/parent axis width}, ticklabel style={font=\tiny}},
colormap={mymap}{[1pt]
  rgb(0pt)=(1,1,0.8);
  rgb(1pt)=(1,0.929411764705882,0.627450980392157);
  rgb(2pt)=(0.996078431372549,0.850980392156863,0.462745098039216);
  rgb(3pt)=(0.996078431372549,0.698039215686274,0.298039215686275);
  rgb(4pt)=(0.992156862745098,0.552941176470588,0.235294117647059);
  rgb(5pt)=(0.988235294117647,0.305882352941176,0.164705882352941);
  rgb(6pt)=(0.890196078431372,0.101960784313725,0.109803921568627);
  rgb(7pt)=(0.741176470588235,0,0.149019607843137);
  rgb(8pt)=(0.501960784313725,0,0.149019607843137)
},
height=\figureheight,
hide x axis,
hide y axis,
point meta max=1.00981603847325,
point meta min=-0.000458797682120604,
tick align=outside,
tick pos=left,
width=\figurewidth,
x grid style={darkgray176},
xmin=0, xmax=3200,
xtick style={color=black},
xtick={0,1000,2000,3000,4000},
xticklabels={
  \(\displaystyle {0}\),
  \(\displaystyle {1000}\),
  \(\displaystyle {2000}\),
  \(\displaystyle {3000}\),
  \(\displaystyle {4000}\)
},
y grid style={darkgray176},
ymin=0, ymax=3360,
ytick style={color=black},
ytick={0,500,1000,1500,2000,2500,3000,3500},
yticklabels={
  \(\displaystyle {0}\),
  \(\displaystyle {500}\),
  \(\displaystyle {1000}\),
  \(\displaystyle {1500}\),
  \(\displaystyle {2000}\),
  \(\displaystyle {2500}\),
  \(\displaystyle {3000}\),
  \(\displaystyle {3500}\)
}
]
\addplot graphics [includegraphics cmd=\pgfimage,xmin=0, xmax=3200, ymin=0, ymax=3360] {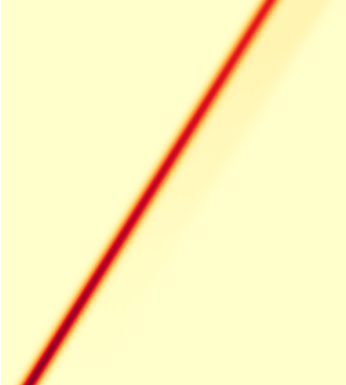};

\end{groupplot}

\draw ({$(current bounding box.south west)!-0.04!(current bounding box.south east)$}|-{$(current bounding box.south west)!0.4!(current bounding box.north west)$}) node[
  scale=0.96,
  anchor=west,
  text=black,
  rotate=90.0
]{time $t$};
\draw ({$(current bounding box.south west)!0.5!(current bounding box.south east)$}|-{$(current bounding box.south west)!-0.06!(current bounding box.north west)$}) node[
  scale=0.96,
  anchor=south,
  text=black,
  rotate=0.0
]{space $x$};
\end{tikzpicture}

%% file: Example_2_snapshots.tex
% This file was created with tikzplotlib v0.10.1.
\begin{tikzpicture}

\definecolor{darkgray176}{RGB}{176,176,176}

\begin{groupplot}[group style={group size=4 by 3, horizontal sep=1.30cm, vertical sep=0.5cm}]
\nextgroupplot[
colorbar,
colorbar style={ytick={-1.5,-1,-0.5,0,0.5,1,1.5},yticklabels={
  \(\displaystyle {\ensuremath{-}1.5}\),
  \(\displaystyle {\ensuremath{-}1.0}\),
  \(\displaystyle {\ensuremath{-}0.5}\),
  \(\displaystyle {0.0}\),
  \(\displaystyle {0.5}\),
  \(\displaystyle {1.0}\),
  \(\displaystyle {1.5}\)
},ylabel={}, width=0.06*\pgfkeysvalueof{/pgfplots/parent axis width}, ticklabel style={font=\tiny}},
colormap={mymap}{[1pt]
  rgb(0pt)=(1,1,0.8);
  rgb(1pt)=(1,0.929411764705882,0.627450980392157);
  rgb(2pt)=(0.996078431372549,0.850980392156863,0.462745098039216);
  rgb(3pt)=(0.996078431372549,0.698039215686274,0.298039215686275);
  rgb(4pt)=(0.992156862745098,0.552941176470588,0.235294117647059);
  rgb(5pt)=(0.988235294117647,0.305882352941176,0.164705882352941);
  rgb(6pt)=(0.890196078431372,0.101960784313725,0.109803921568627);
  rgb(7pt)=(0.741176470588235,0,0.149019607843137);
  rgb(8pt)=(0.501960784313725,0,0.149019607843137)
},
height=\figureheight,
hide x axis,
hide y axis,
point meta max=1.25078919063742,
point meta min=-1.05481896614539,
tick align=outside,
tick pos=left,
align=center,
title={\tiny{FOM}},
width=\figurewidth,
x grid style={darkgray176},
xmin=0, xmax=3200,
xtick style={color=black},
xtick={0,1000,2000,3000,4000},
xticklabels={
  \(\displaystyle {0}\),
  \(\displaystyle {1000}\),
  \(\displaystyle {2000}\),
  \(\displaystyle {3000}\),
  \(\displaystyle {4000}\)
},
y grid style={darkgray176},
ymin=0, ymax=3360,
ytick style={color=black},
ytick={0,500,1000,1500,2000,2500,3000,3500},
yticklabels={
  \(\displaystyle {0}\),
  \(\displaystyle {500}\),
  \(\displaystyle {1000}\),
  \(\displaystyle {1500}\),
  \(\displaystyle {2000}\),
  \(\displaystyle {2500}\),
  \(\displaystyle {3000}\),
  \(\displaystyle {3500}\)
}
]
\addplot graphics [includegraphics cmd=\pgfimage,xmin=0, xmax=3200, ymin=0, ymax=3360] {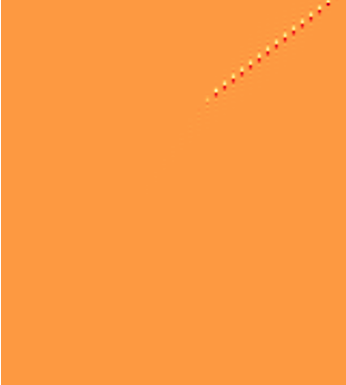};

\nextgroupplot[
colorbar,
colorbar style={ytick={-1.5,-1,-0.5,0,0.5,1,1.5},yticklabels={
  \(\displaystyle {\ensuremath{-}1.5}\),
  \(\displaystyle {\ensuremath{-}1.0}\),
  \(\displaystyle {\ensuremath{-}0.5}\),
  \(\displaystyle {0.0}\),
  \(\displaystyle {0.5}\),
  \(\displaystyle {1.0}\),
  \(\displaystyle {1.5}\)
},ylabel={}, width=0.06*\pgfkeysvalueof{/pgfplots/parent axis width}, ticklabel style={font=\tiny}},
colormap={mymap}{[1pt]
  rgb(0pt)=(1,1,0.8);
  rgb(1pt)=(1,0.929411764705882,0.627450980392157);
  rgb(2pt)=(0.996078431372549,0.850980392156863,0.462745098039216);
  rgb(3pt)=(0.996078431372549,0.698039215686274,0.298039215686275);
  rgb(4pt)=(0.992156862745098,0.552941176470588,0.235294117647059);
  rgb(5pt)=(0.988235294117647,0.305882352941176,0.164705882352941);
  rgb(6pt)=(0.890196078431372,0.101960784313725,0.109803921568627);
  rgb(7pt)=(0.741176470588235,0,0.149019607843137);
  rgb(8pt)=(0.501960784313725,0,0.149019607843137)
},
height=\figureheight,
hide x axis,
hide y axis,
point meta max=1.25078919063742,
point meta min=-1.05481896614539,
tick align=outside,
tick pos=left,
align=center,
title={\tiny{sPOD-G}\\ \tiny{($r=20$)}},
width=\figurewidth,
x grid style={darkgray176},
xmin=0, xmax=3200,
xtick style={color=black},
xtick={0,1000,2000,3000,4000},
xticklabels={
  \(\displaystyle {0}\),
  \(\displaystyle {1000}\),
  \(\displaystyle {2000}\),
  \(\displaystyle {3000}\),
  \(\displaystyle {4000}\)
},
y grid style={darkgray176},
ymin=0, ymax=3360,
ytick style={color=black},
ytick={0,500,1000,1500,2000,2500,3000,3500},
yticklabels={
  \(\displaystyle {0}\),
  \(\displaystyle {500}\),
  \(\displaystyle {1000}\),
  \(\displaystyle {1500}\),
  \(\displaystyle {2000}\),
  \(\displaystyle {2500}\),
  \(\displaystyle {3000}\),
  \(\displaystyle {3500}\)
}
]
\addplot graphics [includegraphics cmd=\pgfimage,xmin=0, xmax=3200, ymin=0, ymax=3360] {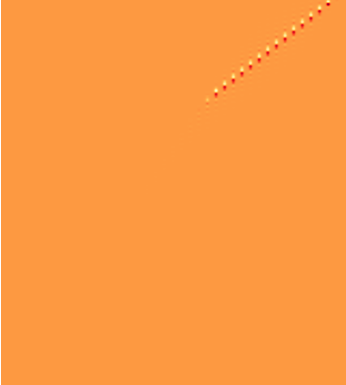};

\nextgroupplot[
colorbar,
colorbar style={ytick={-1.5,-1,-0.5,0,0.5,1,1.5},yticklabels={
  \(\displaystyle {\ensuremath{-}1.5}\),
  \(\displaystyle {\ensuremath{-}1.0}\),
  \(\displaystyle {\ensuremath{-}0.5}\),
  \(\displaystyle {0.0}\),
  \(\displaystyle {0.5}\),
  \(\displaystyle {1.0}\),
  \(\displaystyle {1.5}\)
},ylabel={}, width=0.06*\pgfkeysvalueof{/pgfplots/parent axis width}, ticklabel style={font=\tiny}},
colormap={mymap}{[1pt]
  rgb(0pt)=(1,1,0.8);
  rgb(1pt)=(1,0.929411764705882,0.627450980392157);
  rgb(2pt)=(0.996078431372549,0.850980392156863,0.462745098039216);
  rgb(3pt)=(0.996078431372549,0.698039215686274,0.298039215686275);
  rgb(4pt)=(0.992156862745098,0.552941176470588,0.235294117647059);
  rgb(5pt)=(0.988235294117647,0.305882352941176,0.164705882352941);
  rgb(6pt)=(0.890196078431372,0.101960784313725,0.109803921568627);
  rgb(7pt)=(0.741176470588235,0,0.149019607843137);
  rgb(8pt)=(0.501960784313725,0,0.149019607843137)
},
height=\figureheight,
hide x axis,
hide y axis,
point meta max=1.25078919063742,
point meta min=-1.05481896614539,
tick align=outside,
tick pos=left,
align=center,
title={\tiny{POD-G}\\ \tiny{($\mathrm{p}=120$)}},
width=\figurewidth,
x grid style={darkgray176},
xmin=0, xmax=3200,
xtick style={color=black},
xtick={0,1000,2000,3000,4000},
xticklabels={
  \(\displaystyle {0}\),
  \(\displaystyle {1000}\),
  \(\displaystyle {2000}\),
  \(\displaystyle {3000}\),
  \(\displaystyle {4000}\)
},
y grid style={darkgray176},
ymin=0, ymax=3360,
ytick style={color=black},
ytick={0,500,1000,1500,2000,2500,3000,3500},
yticklabels={
  \(\displaystyle {0}\),
  \(\displaystyle {500}\),
  \(\displaystyle {1000}\),
  \(\displaystyle {1500}\),
  \(\displaystyle {2000}\),
  \(\displaystyle {2500}\),
  \(\displaystyle {3000}\),
  \(\displaystyle {3500}\)
}
]
\addplot graphics [includegraphics cmd=\pgfimage,xmin=0, xmax=3200, ymin=0, ymax=3360] {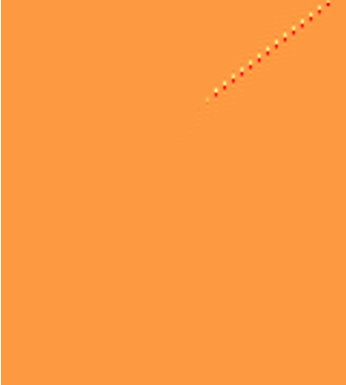};

\nextgroupplot[
colorbar,
colorbar style={ytick={-1.5,-1,-0.5,0,0.5,1,1.5},yticklabels={
  \(\displaystyle {\ensuremath{-}1.5}\),
  \(\displaystyle {\ensuremath{-}1.0}\),
  \(\displaystyle {\ensuremath{-}0.5}\),
  \(\displaystyle {0.0}\),
  \(\displaystyle {0.5}\),
  \(\displaystyle {1.0}\),
  \(\displaystyle {1.5}\)
},ylabel={}, width=0.06*\pgfkeysvalueof{/pgfplots/parent axis width}, ticklabel style={font=\tiny}},
colormap={mymap}{[1pt]
  rgb(0pt)=(1,1,0.8);
  rgb(1pt)=(1,0.929411764705882,0.627450980392157);
  rgb(2pt)=(0.996078431372549,0.850980392156863,0.462745098039216);
  rgb(3pt)=(0.996078431372549,0.698039215686274,0.298039215686275);
  rgb(4pt)=(0.992156862745098,0.552941176470588,0.235294117647059);
  rgb(5pt)=(0.988235294117647,0.305882352941176,0.164705882352941);
  rgb(6pt)=(0.890196078431372,0.101960784313725,0.109803921568627);
  rgb(7pt)=(0.741176470588235,0,0.149019607843137);
  rgb(8pt)=(0.501960784313725,0,0.149019607843137)
},
height=\figureheight,
hide x axis,
hide y axis,
point meta max=1.25078919063742,
point meta min=-1.05481896614539,
tick align=outside,
tick pos=left,
align=center,
title={\tiny{POD-G} \\ \tiny{($\mathrm{p}=20$)}},
width=\figurewidth,
x grid style={darkgray176},
xmin=0, xmax=3200,
xtick style={color=black},
xtick={0,1000,2000,3000,4000},
xticklabels={
  \(\displaystyle {0}\),
  \(\displaystyle {1000}\),
  \(\displaystyle {2000}\),
  \(\displaystyle {3000}\),
  \(\displaystyle {4000}\)
},
y grid style={darkgray176},
ymin=0, ymax=3360,
ytick style={color=black},
ytick={0,500,1000,1500,2000,2500,3000,3500},
yticklabels={
  \(\displaystyle {0}\),
  \(\displaystyle {500}\),
  \(\displaystyle {1000}\),
  \(\displaystyle {1500}\),
  \(\displaystyle {2000}\),
  \(\displaystyle {2500}\),
  \(\displaystyle {3000}\),
  \(\displaystyle {3500}\)
}
]
\addplot graphics [includegraphics cmd=\pgfimage,xmin=0, xmax=3200, ymin=0, ymax=3360] {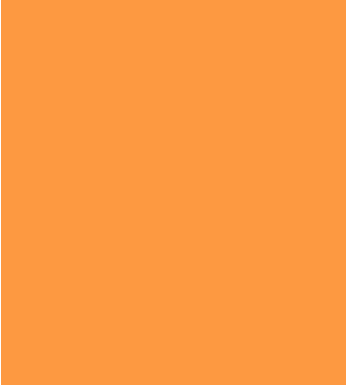};

%%%%%%%%%%%%%%%%%%%%%%%%%%%%%%%%%%%%%%%%%%%%%%%%%%%%%%%%%%%%%%%%%%%%%%%%%%%%%

\nextgroupplot[
colorbar,
colorbar style={ytick={-0.03,-0.02,-0.01,0,0.01,0.02,0.03},yticklabels={
  \(\displaystyle {\ensuremath{-}0.03}\),
  \(\displaystyle {\ensuremath{-}0.02}\),
  \(\displaystyle {\ensuremath{-}0.01}\),
  \(\displaystyle {0.00}\),
  \(\displaystyle {0.01}\),
  \(\displaystyle {0.02}\),
  \(\displaystyle {0.03}\)
},ylabel={}, width=0.06*\pgfkeysvalueof{/pgfplots/parent axis width}, ticklabel style={font=\tiny}},
colormap={mymap}{[1pt]
  rgb(0pt)=(1,1,0.8);
  rgb(1pt)=(1,0.929411764705882,0.627450980392157);
  rgb(2pt)=(0.996078431372549,0.850980392156863,0.462745098039216);
  rgb(3pt)=(0.996078431372549,0.698039215686274,0.298039215686275);
  rgb(4pt)=(0.992156862745098,0.552941176470588,0.235294117647059);
  rgb(5pt)=(0.988235294117647,0.305882352941176,0.164705882352941);
  rgb(6pt)=(0.890196078431372,0.101960784313725,0.109803921568627);
  rgb(7pt)=(0.741176470588235,0,0.149019607843137);
  rgb(8pt)=(0.501960784313725,0,0.149019607843137)
},
height=\figureheight,
hide x axis,
hide y axis,
point meta max=0.0266212879162842,
point meta min=-0.0270272308742578,
tick align=outside,
tick pos=left,
width=\figurewidth,
x grid style={darkgray176},
xmin=0, xmax=3200,
xtick style={color=black},
xtick={0,1000,2000,3000,4000},
xticklabels={
  \(\displaystyle {0}\),
  \(\displaystyle {1000}\),
  \(\displaystyle {2000}\),
  \(\displaystyle {3000}\),
  \(\displaystyle {4000}\)
},
y grid style={darkgray176},
ymin=0, ymax=3360,
ytick style={color=black},
ytick={0,500,1000,1500,2000,2500,3000,3500},
yticklabels={
  \(\displaystyle {0}\),
  \(\displaystyle {500}\),
  \(\displaystyle {1000}\),
  \(\displaystyle {1500}\),
  \(\displaystyle {2000}\),
  \(\displaystyle {2500}\),
  \(\displaystyle {3000}\),
  \(\displaystyle {3500}\)
}
]
\addplot graphics [includegraphics cmd=\pgfimage,xmin=0, xmax=3200, ymin=0, ymax=3360] {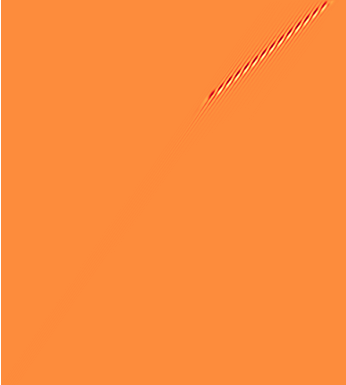};

\nextgroupplot[
colorbar,
colorbar style={ytick={-0.03,-0.02,-0.01,0,0.01,0.02,0.03},yticklabels={
  \(\displaystyle {\ensuremath{-}0.03}\),
  \(\displaystyle {\ensuremath{-}0.02}\),
  \(\displaystyle {\ensuremath{-}0.01}\),
  \(\displaystyle {0.00}\),
  \(\displaystyle {0.01}\),
  \(\displaystyle {0.02}\),
  \(\displaystyle {0.03}\)
},ylabel={}, width=0.06*\pgfkeysvalueof{/pgfplots/parent axis width}, ticklabel style={font=\tiny}},
colormap={mymap}{[1pt]
  rgb(0pt)=(1,1,0.8);
  rgb(1pt)=(1,0.929411764705882,0.627450980392157);
  rgb(2pt)=(0.996078431372549,0.850980392156863,0.462745098039216);
  rgb(3pt)=(0.996078431372549,0.698039215686274,0.298039215686275);
  rgb(4pt)=(0.992156862745098,0.552941176470588,0.235294117647059);
  rgb(5pt)=(0.988235294117647,0.305882352941176,0.164705882352941);
  rgb(6pt)=(0.890196078431372,0.101960784313725,0.109803921568627);
  rgb(7pt)=(0.741176470588235,0,0.149019607843137);
  rgb(8pt)=(0.501960784313725,0,0.149019607843137)
},
height=\figureheight,
hide x axis,
hide y axis,
point meta max=0.0266212879162842,
point meta min=-0.0270272308742578,
tick align=outside,
tick pos=left,
width=\figurewidth,
x grid style={darkgray176},
xmin=0, xmax=3200,
xtick style={color=black},
xtick={0,1000,2000,3000,4000},
xticklabels={
  \(\displaystyle {0}\),
  \(\displaystyle {1000}\),
  \(\displaystyle {2000}\),
  \(\displaystyle {3000}\),
  \(\displaystyle {4000}\)
},
y grid style={darkgray176},
ymin=0, ymax=3360,
ytick style={color=black},
ytick={0,500,1000,1500,2000,2500,3000,3500},
yticklabels={
  \(\displaystyle {0}\),
  \(\displaystyle {500}\),
  \(\displaystyle {1000}\),
  \(\displaystyle {1500}\),
  \(\displaystyle {2000}\),
  \(\displaystyle {2500}\),
  \(\displaystyle {3000}\),
  \(\displaystyle {3500}\)
}
]
\addplot graphics [includegraphics cmd=\pgfimage,xmin=0, xmax=3200, ymin=0, ymax=3360] {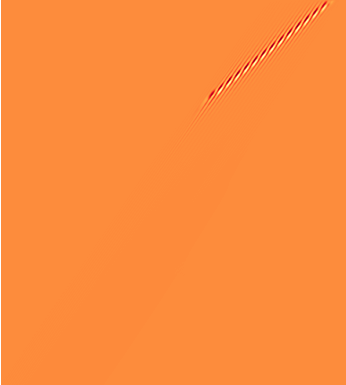};

\nextgroupplot[
colorbar,
colorbar style={ytick={-0.03,-0.02,-0.01,0,0.01,0.02,0.03},yticklabels={
  \(\displaystyle {\ensuremath{-}0.03}\),
  \(\displaystyle {\ensuremath{-}0.02}\),
  \(\displaystyle {\ensuremath{-}0.01}\),
  \(\displaystyle {0.00}\),
  \(\displaystyle {0.01}\),
  \(\displaystyle {0.02}\),
  \(\displaystyle {0.03}\)
},ylabel={}, width=0.06*\pgfkeysvalueof{/pgfplots/parent axis width}, ticklabel style={font=\tiny}},
colormap={mymap}{[1pt]
  rgb(0pt)=(1,1,0.8);
  rgb(1pt)=(1,0.929411764705882,0.627450980392157);
  rgb(2pt)=(0.996078431372549,0.850980392156863,0.462745098039216);
  rgb(3pt)=(0.996078431372549,0.698039215686274,0.298039215686275);
  rgb(4pt)=(0.992156862745098,0.552941176470588,0.235294117647059);
  rgb(5pt)=(0.988235294117647,0.305882352941176,0.164705882352941);
  rgb(6pt)=(0.890196078431372,0.101960784313725,0.109803921568627);
  rgb(7pt)=(0.741176470588235,0,0.149019607843137);
  rgb(8pt)=(0.501960784313725,0,0.149019607843137)
},
height=\figureheight,
hide x axis,
hide y axis,
point meta max=0.0266212879162842,
point meta min=-0.0270272308742578,
tick align=outside,
tick pos=left,
width=\figurewidth,
x grid style={darkgray176},
xmin=0, xmax=3200,
xtick style={color=black},
xtick={0,1000,2000,3000,4000},
xticklabels={
  \(\displaystyle {0}\),
  \(\displaystyle {1000}\),
  \(\displaystyle {2000}\),
  \(\displaystyle {3000}\),
  \(\displaystyle {4000}\)
},
y grid style={darkgray176},
ymin=0, ymax=3360,
ytick style={color=black},
ytick={0,500,1000,1500,2000,2500,3000,3500},
yticklabels={
  \(\displaystyle {0}\),
  \(\displaystyle {500}\),
  \(\displaystyle {1000}\),
  \(\displaystyle {1500}\),
  \(\displaystyle {2000}\),
  \(\displaystyle {2500}\),
  \(\displaystyle {3000}\),
  \(\displaystyle {3500}\)
}
]
\addplot graphics [includegraphics cmd=\pgfimage,xmin=0, xmax=3200, ymin=0, ymax=3360] {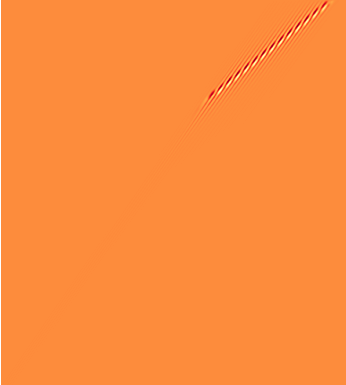};

\nextgroupplot[
colorbar,
colorbar style={ytick={-0.03,-0.02,-0.01,0,0.01,0.02,0.03},yticklabels={
  \(\displaystyle {\ensuremath{-}0.03}\),
  \(\displaystyle {\ensuremath{-}0.02}\),
  \(\displaystyle {\ensuremath{-}0.01}\),
  \(\displaystyle {0.00}\),
  \(\displaystyle {0.01}\),
  \(\displaystyle {0.02}\),
  \(\displaystyle {0.03}\)
},ylabel={}, width=0.06*\pgfkeysvalueof{/pgfplots/parent axis width}, ticklabel style={font=\tiny}},
colormap={mymap}{[1pt]
  rgb(0pt)=(1,1,0.8);
  rgb(1pt)=(1,0.929411764705882,0.627450980392157);
  rgb(2pt)=(0.996078431372549,0.850980392156863,0.462745098039216);
  rgb(3pt)=(0.996078431372549,0.698039215686274,0.298039215686275);
  rgb(4pt)=(0.992156862745098,0.552941176470588,0.235294117647059);
  rgb(5pt)=(0.988235294117647,0.305882352941176,0.164705882352941);
  rgb(6pt)=(0.890196078431372,0.101960784313725,0.109803921568627);
  rgb(7pt)=(0.741176470588235,0,0.149019607843137);
  rgb(8pt)=(0.501960784313725,0,0.149019607843137)
},
height=\figureheight,
hide x axis,
hide y axis,
point meta max=0.0266212879162842,
point meta min=-0.0270272308742578,
tick align=outside,
tick pos=left,
width=\figurewidth,
x grid style={darkgray176},
xmin=0, xmax=3200,
xtick style={color=black},
xtick={0,1000,2000,3000,4000},
xticklabels={
  \(\displaystyle {0}\),
  \(\displaystyle {1000}\),
  \(\displaystyle {2000}\),
  \(\displaystyle {3000}\),
  \(\displaystyle {4000}\)
},
y grid style={darkgray176},
ymin=0, ymax=3360,
ytick style={color=black},
ytick={0,500,1000,1500,2000,2500,3000,3500},
yticklabels={
  \(\displaystyle {0}\),
  \(\displaystyle {500}\),
  \(\displaystyle {1000}\),
  \(\displaystyle {1500}\),
  \(\displaystyle {2000}\),
  \(\displaystyle {2500}\),
  \(\displaystyle {3000}\),
  \(\displaystyle {3500}\)
}
]
\addplot graphics [includegraphics cmd=\pgfimage,xmin=0, xmax=3200, ymin=0, ymax=3360] {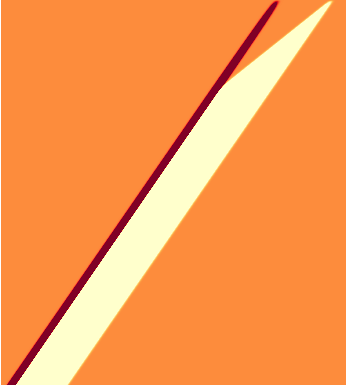};

%%%%%%%%%%%%%%%%%%%%%%%%%%%%%%%%%%%%%%%%%%%%%%%%%%%%%%%%%%%%%%%%%%%%%%%%%%%%%%%%%%%%%%

\nextgroupplot[
colorbar,
colorbar style={ytick={-0.2,0,0.2,0.4,0.6,0.8,1,1.2},yticklabels={
  \(\displaystyle {\ensuremath{-}0.2}\),
  \(\displaystyle {0.0}\),
  \(\displaystyle {0.2}\),
  \(\displaystyle {0.4}\),
  \(\displaystyle {0.6}\),
  \(\displaystyle {0.8}\),
  \(\displaystyle {1.0}\),
  \(\displaystyle {1.2}\)
},ylabel={}, width=0.06*\pgfkeysvalueof{/pgfplots/parent axis width}, ticklabel style={font=\tiny}},
colormap={mymap}{[1pt]
  rgb(0pt)=(1,1,0.8);
  rgb(1pt)=(1,0.929411764705882,0.627450980392157);
  rgb(2pt)=(0.996078431372549,0.850980392156863,0.462745098039216);
  rgb(3pt)=(0.996078431372549,0.698039215686274,0.298039215686275);
  rgb(4pt)=(0.992156862745098,0.552941176470588,0.235294117647059);
  rgb(5pt)=(0.988235294117647,0.305882352941176,0.164705882352941);
  rgb(6pt)=(0.890196078431372,0.101960784313725,0.109803921568627);
  rgb(7pt)=(0.741176470588235,0,0.149019607843137);
  rgb(8pt)=(0.501960784313725,0,0.149019607843137)
},
height=\figureheight,
hide x axis,
hide y axis,
point meta max=1.065146755835,
point meta min=-0.0993750088945431,
tick align=outside,
tick pos=left,
width=\figurewidth,
x grid style={darkgray176},
xmin=0, xmax=3200,
xtick style={color=black},
xtick={0,1000,2000,3000,4000},
xticklabels={
  \(\displaystyle {0}\),
  \(\displaystyle {1000}\),
  \(\displaystyle {2000}\),
  \(\displaystyle {3000}\),
  \(\displaystyle {4000}\)
},
y grid style={darkgray176},
ymin=0, ymax=3360,
ytick style={color=black},
ytick={0,500,1000,1500,2000,2500,3000,3500},
yticklabels={
  \(\displaystyle {0}\),
  \(\displaystyle {500}\),
  \(\displaystyle {1000}\),
  \(\displaystyle {1500}\),
  \(\displaystyle {2000}\),
  \(\displaystyle {2500}\),
  \(\displaystyle {3000}\),
  \(\displaystyle {3500}\)
}
]
\addplot graphics [includegraphics cmd=\pgfimage,xmin=0, xmax=3200, ymin=0, ymax=3360] {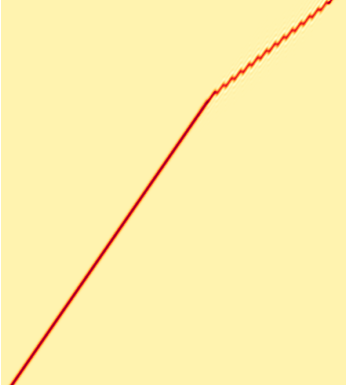};

\nextgroupplot[
colorbar,
colorbar style={ytick={-0.2,0,0.2,0.4,0.6,0.8,1,1.2},yticklabels={
  \(\displaystyle {\ensuremath{-}0.2}\),
  \(\displaystyle {0.0}\),
  \(\displaystyle {0.2}\),
  \(\displaystyle {0.4}\),
  \(\displaystyle {0.6}\),
  \(\displaystyle {0.8}\),
  \(\displaystyle {1.0}\),
  \(\displaystyle {1.2}\)
},ylabel={}, width=0.06*\pgfkeysvalueof{/pgfplots/parent axis width}, ticklabel style={font=\tiny}},
colormap={mymap}{[1pt]
  rgb(0pt)=(1,1,0.8);
  rgb(1pt)=(1,0.929411764705882,0.627450980392157);
  rgb(2pt)=(0.996078431372549,0.850980392156863,0.462745098039216);
  rgb(3pt)=(0.996078431372549,0.698039215686274,0.298039215686275);
  rgb(4pt)=(0.992156862745098,0.552941176470588,0.235294117647059);
  rgb(5pt)=(0.988235294117647,0.305882352941176,0.164705882352941);
  rgb(6pt)=(0.890196078431372,0.101960784313725,0.109803921568627);
  rgb(7pt)=(0.741176470588235,0,0.149019607843137);
  rgb(8pt)=(0.501960784313725,0,0.149019607843137)
},
height=\figureheight,
hide x axis,
hide y axis,
point meta max=1.065146755835,
point meta min=-0.0993750088945431,
tick align=outside,
tick pos=left,
width=\figurewidth,
x grid style={darkgray176},
xmin=0, xmax=3200,
xtick style={color=black},
xtick={0,1000,2000,3000,4000},
xticklabels={
  \(\displaystyle {0}\),
  \(\displaystyle {1000}\),
  \(\displaystyle {2000}\),
  \(\displaystyle {3000}\),
  \(\displaystyle {4000}\)
},
y grid style={darkgray176},
ymin=0, ymax=3360,
ytick style={color=black},
ytick={0,500,1000,1500,2000,2500,3000,3500},
yticklabels={
  \(\displaystyle {0}\),
  \(\displaystyle {500}\),
  \(\displaystyle {1000}\),
  \(\displaystyle {1500}\),
  \(\displaystyle {2000}\),
  \(\displaystyle {2500}\),
  \(\displaystyle {3000}\),
  \(\displaystyle {3500}\)
}
]
\addplot graphics [includegraphics cmd=\pgfimage,xmin=0, xmax=3200, ymin=0, ymax=3360] {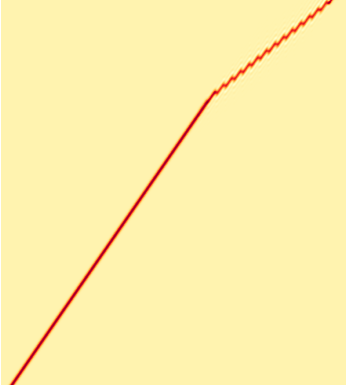};

\nextgroupplot[
colorbar,
colorbar style={ytick={-0.2,0,0.2,0.4,0.6,0.8,1,1.2},yticklabels={
  \(\displaystyle {\ensuremath{-}0.2}\),
  \(\displaystyle {0.0}\),
  \(\displaystyle {0.2}\),
  \(\displaystyle {0.4}\),
  \(\displaystyle {0.6}\),
  \(\displaystyle {0.8}\),
  \(\displaystyle {1.0}\),
  \(\displaystyle {1.2}\)
},ylabel={}, width=0.06*\pgfkeysvalueof{/pgfplots/parent axis width}, ticklabel style={font=\tiny}},
colormap={mymap}{[1pt]
  rgb(0pt)=(1,1,0.8);
  rgb(1pt)=(1,0.929411764705882,0.627450980392157);
  rgb(2pt)=(0.996078431372549,0.850980392156863,0.462745098039216);
  rgb(3pt)=(0.996078431372549,0.698039215686274,0.298039215686275);
  rgb(4pt)=(0.992156862745098,0.552941176470588,0.235294117647059);
  rgb(5pt)=(0.988235294117647,0.305882352941176,0.164705882352941);
  rgb(6pt)=(0.890196078431372,0.101960784313725,0.109803921568627);
  rgb(7pt)=(0.741176470588235,0,0.149019607843137);
  rgb(8pt)=(0.501960784313725,0,0.149019607843137)
},
height=\figureheight,
hide x axis,
hide y axis,
point meta max=1.065146755835,
point meta min=-0.0993750088945431,
tick align=outside,
tick pos=left,
width=\figurewidth,
x grid style={darkgray176},
xmin=0, xmax=3200,
xtick style={color=black},
xtick={0,1000,2000,3000,4000},
xticklabels={
  \(\displaystyle {0}\),
  \(\displaystyle {1000}\),
  \(\displaystyle {2000}\),
  \(\displaystyle {3000}\),
  \(\displaystyle {4000}\)
},
y grid style={darkgray176},
ymin=0, ymax=3360,
ytick style={color=black},
ytick={0,500,1000,1500,2000,2500,3000,3500},
yticklabels={
  \(\displaystyle {0}\),
  \(\displaystyle {500}\),
  \(\displaystyle {1000}\),
  \(\displaystyle {1500}\),
  \(\displaystyle {2000}\),
  \(\displaystyle {2500}\),
  \(\displaystyle {3000}\),
  \(\displaystyle {3500}\)
}
]
\addplot graphics [includegraphics cmd=\pgfimage,xmin=0, xmax=3200, ymin=0, ymax=3360] {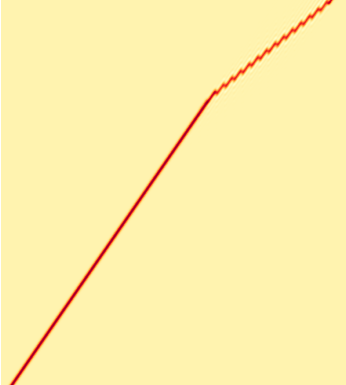};

\nextgroupplot[
colorbar,
colorbar style={ytick={-0.2,0,0.2,0.4,0.6,0.8,1,1.2},yticklabels={
  \(\displaystyle {\ensuremath{-}0.2}\),
  \(\displaystyle {0.0}\),
  \(\displaystyle {0.2}\),
  \(\displaystyle {0.4}\),
  \(\displaystyle {0.6}\),
  \(\displaystyle {0.8}\),
  \(\displaystyle {1.0}\),
  \(\displaystyle {1.2}\)
},ylabel={}, width=0.06*\pgfkeysvalueof{/pgfplots/parent axis width}, ticklabel style={font=\tiny}},
colormap={mymap}{[1pt]
  rgb(0pt)=(1,1,0.8);
  rgb(1pt)=(1,0.929411764705882,0.627450980392157);
  rgb(2pt)=(0.996078431372549,0.850980392156863,0.462745098039216);
  rgb(3pt)=(0.996078431372549,0.698039215686274,0.298039215686275);
  rgb(4pt)=(0.992156862745098,0.552941176470588,0.235294117647059);
  rgb(5pt)=(0.988235294117647,0.305882352941176,0.164705882352941);
  rgb(6pt)=(0.890196078431372,0.101960784313725,0.109803921568627);
  rgb(7pt)=(0.741176470588235,0,0.149019607843137);
  rgb(8pt)=(0.501960784313725,0,0.149019607843137)
},
height=\figureheight,
hide x axis,
hide y axis,
point meta max=1.065146755835,
point meta min=-0.0993750088945431,
tick align=outside,
tick pos=left,
width=\figurewidth,
x grid style={darkgray176},
xmin=0, xmax=3200,
xtick style={color=black},
xtick={0,1000,2000,3000,4000},
xticklabels={
  \(\displaystyle {0}\),
  \(\displaystyle {1000}\),
  \(\displaystyle {2000}\),
  \(\displaystyle {3000}\),
  \(\displaystyle {4000}\)
},
y grid style={darkgray176},
ymin=0, ymax=3360,
ytick style={color=black},
ytick={0,500,1000,1500,2000,2500,3000,3500},
yticklabels={
  \(\displaystyle {0}\),
  \(\displaystyle {500}\),
  \(\displaystyle {1000}\),
  \(\displaystyle {1500}\),
  \(\displaystyle {2000}\),
  \(\displaystyle {2500}\),
  \(\displaystyle {3000}\),
  \(\displaystyle {3500}\)
}
]
\addplot graphics [includegraphics cmd=\pgfimage,xmin=0, xmax=3200, ymin=0, ymax=3360] {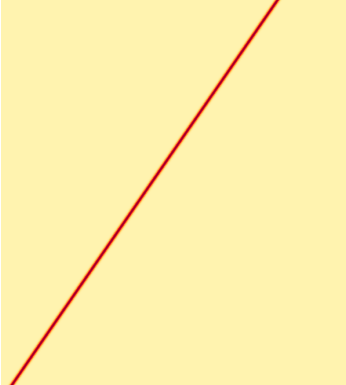};

\end{groupplot}

\draw ({$(current bounding box.south west)!-0.04!(current bounding box.south east)$}|-{$(current bounding box.south west)!0.4!(current bounding box.north west)$}) node[
  scale=0.96,
  anchor=west,
  text=black,
  rotate=90.0
]{time $t$};
\draw ({$(current bounding box.south west)!0.5!(current bounding box.south east)$}|-{$(current bounding box.south west)!-0.06!(current bounding box.north west)$}) node[
  scale=0.96,
  anchor=south,
  text=black,
  rotate=0.0
]{space $x$};
\end{tikzpicture}

%% file: Example_3_snapshots.tex
% This file was created with tikzplotlib v0.10.1.
\begin{tikzpicture}

\definecolor{darkgray176}{RGB}{176,176,176}

\begin{groupplot}[group style={group size=4 by 3, horizontal sep=1.30cm, vertical sep=0.5cm}]
\nextgroupplot[
colorbar,
colorbar style={ytick={-1.5,-1,-0.5,0,0.5,1,1.5},yticklabels={
  \(\displaystyle {\ensuremath{-}1.5}\),
  \(\displaystyle {\ensuremath{-}1.0}\),
  \(\displaystyle {\ensuremath{-}0.5}\),
  \(\displaystyle {0.0}\),
  \(\displaystyle {0.5}\),
  \(\displaystyle {1.0}\),
  \(\displaystyle {1.5}\)
},ylabel={}, width=0.06*\pgfkeysvalueof{/pgfplots/parent axis width}, ticklabel style={font=\tiny}},
colormap={mymap}{[1pt]
  rgb(0pt)=(1,1,0.8);
  rgb(1pt)=(1,0.929411764705882,0.627450980392157);
  rgb(2pt)=(0.996078431372549,0.850980392156863,0.462745098039216);
  rgb(3pt)=(0.996078431372549,0.698039215686274,0.298039215686275);
  rgb(4pt)=(0.992156862745098,0.552941176470588,0.235294117647059);
  rgb(5pt)=(0.988235294117647,0.305882352941176,0.164705882352941);
  rgb(6pt)=(0.890196078431372,0.101960784313725,0.109803921568627);
  rgb(7pt)=(0.741176470588235,0,0.149019607843137);
  rgb(8pt)=(0.501960784313725,0,0.149019607843137)
},
height=\figureheight,
hide x axis,
hide y axis,
point meta max=1.45636822945787,
point meta min=-1.41023053305731,
tick align=outside,
tick pos=left,
align=center,
title={\tiny{FOM}},
width=\figurewidth,
x grid style={darkgray176},
xmin=0, xmax=3200,
xtick style={color=black},
xtick={0,1000,2000,3000,4000},
xticklabels={
  \(\displaystyle {0}\),
  \(\displaystyle {1000}\),
  \(\displaystyle {2000}\),
  \(\displaystyle {3000}\),
  \(\displaystyle {4000}\)
},
y grid style={darkgray176},
ymin=0, ymax=3360,
ytick style={color=black},
ytick={0,500,1000,1500,2000,2500,3000,3500},
yticklabels={
  \(\displaystyle {0}\),
  \(\displaystyle {500}\),
  \(\displaystyle {1000}\),
  \(\displaystyle {1500}\),
  \(\displaystyle {2000}\),
  \(\displaystyle {2500}\),
  \(\displaystyle {3000}\),
  \(\displaystyle {3500}\)
}
]
\addplot graphics [includegraphics cmd=\pgfimage,xmin=0, xmax=3200, ymin=0, ymax=3360] {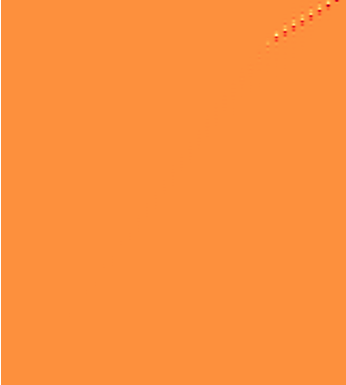};

\nextgroupplot[
colorbar,
colorbar style={ytick={-1.5,-1,-0.5,0,0.5,1,1.5},yticklabels={
  \(\displaystyle {\ensuremath{-}1.5}\),
  \(\displaystyle {\ensuremath{-}1.0}\),
  \(\displaystyle {\ensuremath{-}0.5}\),
  \(\displaystyle {0.0}\),
  \(\displaystyle {0.5}\),
  \(\displaystyle {1.0}\),
  \(\displaystyle {1.5}\)
},ylabel={}, width=0.06*\pgfkeysvalueof{/pgfplots/parent axis width}, ticklabel style={font=\tiny}},
colormap={mymap}{[1pt]
  rgb(0pt)=(1,1,0.8);
  rgb(1pt)=(1,0.929411764705882,0.627450980392157);
  rgb(2pt)=(0.996078431372549,0.850980392156863,0.462745098039216);
  rgb(3pt)=(0.996078431372549,0.698039215686274,0.298039215686275);
  rgb(4pt)=(0.992156862745098,0.552941176470588,0.235294117647059);
  rgb(5pt)=(0.988235294117647,0.305882352941176,0.164705882352941);
  rgb(6pt)=(0.890196078431372,0.101960784313725,0.109803921568627);
  rgb(7pt)=(0.741176470588235,0,0.149019607843137);
  rgb(8pt)=(0.501960784313725,0,0.149019607843137)
},
height=\figureheight,
hide x axis,
hide y axis,
point meta max=1.45636822945787,
point meta min=-1.41023053305731,
tick align=outside,
tick pos=left,
align=center,
title={\tiny{sPOD-G}\\ \tiny{($r=12$)}},
width=\figurewidth,
x grid style={darkgray176},
xmin=0, xmax=3200,
xtick style={color=black},
xtick={0,1000,2000,3000,4000},
xticklabels={
  \(\displaystyle {0}\),
  \(\displaystyle {1000}\),
  \(\displaystyle {2000}\),
  \(\displaystyle {3000}\),
  \(\displaystyle {4000}\)
},
y grid style={darkgray176},
ymin=0, ymax=3360,
ytick style={color=black},
ytick={0,500,1000,1500,2000,2500,3000,3500},
yticklabels={
  \(\displaystyle {0}\),
  \(\displaystyle {500}\),
  \(\displaystyle {1000}\),
  \(\displaystyle {1500}\),
  \(\displaystyle {2000}\),
  \(\displaystyle {2500}\),
  \(\displaystyle {3000}\),
  \(\displaystyle {3500}\)
}
]
\addplot graphics [includegraphics cmd=\pgfimage,xmin=0, xmax=3200, ymin=0, ymax=3360] {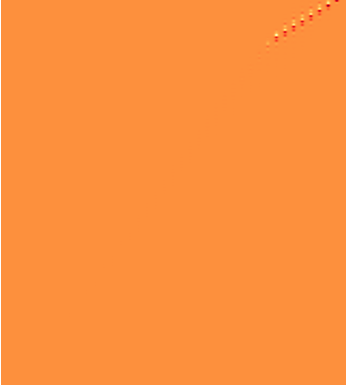};

\nextgroupplot[
colorbar,
colorbar style={ytick={-1.5,-1,-0.5,0,0.5,1,1.5},yticklabels={
  \(\displaystyle {\ensuremath{-}1.5}\),
  \(\displaystyle {\ensuremath{-}1.0}\),
  \(\displaystyle {\ensuremath{-}0.5}\),
  \(\displaystyle {0.0}\),
  \(\displaystyle {0.5}\),
  \(\displaystyle {1.0}\),
  \(\displaystyle {1.5}\)
},ylabel={}, width=0.06*\pgfkeysvalueof{/pgfplots/parent axis width}, ticklabel style={font=\tiny}},
colormap={mymap}{[1pt]
  rgb(0pt)=(1,1,0.8);
  rgb(1pt)=(1,0.929411764705882,0.627450980392157);
  rgb(2pt)=(0.996078431372549,0.850980392156863,0.462745098039216);
  rgb(3pt)=(0.996078431372549,0.698039215686274,0.298039215686275);
  rgb(4pt)=(0.992156862745098,0.552941176470588,0.235294117647059);
  rgb(5pt)=(0.988235294117647,0.305882352941176,0.164705882352941);
  rgb(6pt)=(0.890196078431372,0.101960784313725,0.109803921568627);
  rgb(7pt)=(0.741176470588235,0,0.149019607843137);
  rgb(8pt)=(0.501960784313725,0,0.149019607843137)
},
height=\figureheight,
hide x axis,
hide y axis,
point meta max=1.45636822945787,
point meta min=-1.41023053305731,
tick align=outside,
tick pos=left,
align=center,
title={\tiny{POD-G}\\ \tiny{($\mathrm{p}=120$)}},
width=\figurewidth,
x grid style={darkgray176},
xmin=0, xmax=3200,
xtick style={color=black},
xtick={0,1000,2000,3000,4000},
xticklabels={
  \(\displaystyle {0}\),
  \(\displaystyle {1000}\),
  \(\displaystyle {2000}\),
  \(\displaystyle {3000}\),
  \(\displaystyle {4000}\)
},
y grid style={darkgray176},
ymin=0, ymax=3360,
ytick style={color=black},
ytick={0,500,1000,1500,2000,2500,3000,3500},
yticklabels={
  \(\displaystyle {0}\),
  \(\displaystyle {500}\),
  \(\displaystyle {1000}\),
  \(\displaystyle {1500}\),
  \(\displaystyle {2000}\),
  \(\displaystyle {2500}\),
  \(\displaystyle {3000}\),
  \(\displaystyle {3500}\)
}
]
\addplot graphics [includegraphics cmd=\pgfimage,xmin=0, xmax=3200, ymin=0, ymax=3360] {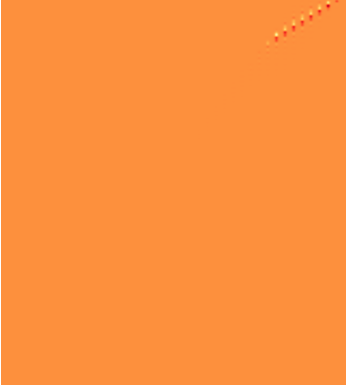};

\nextgroupplot[
colorbar,
colorbar style={ytick={-1.5,-1,-0.5,0,0.5,1,1.5},yticklabels={
  \(\displaystyle {\ensuremath{-}1.5}\),
  \(\displaystyle {\ensuremath{-}1.0}\),
  \(\displaystyle {\ensuremath{-}0.5}\),
  \(\displaystyle {0.0}\),
  \(\displaystyle {0.5}\),
  \(\displaystyle {1.0}\),
  \(\displaystyle {1.5}\)
},ylabel={}, width=0.06*\pgfkeysvalueof{/pgfplots/parent axis width}, ticklabel style={font=\tiny}},
colormap={mymap}{[1pt]
  rgb(0pt)=(1,1,0.8);
  rgb(1pt)=(1,0.929411764705882,0.627450980392157);
  rgb(2pt)=(0.996078431372549,0.850980392156863,0.462745098039216);
  rgb(3pt)=(0.996078431372549,0.698039215686274,0.298039215686275);
  rgb(4pt)=(0.992156862745098,0.552941176470588,0.235294117647059);
  rgb(5pt)=(0.988235294117647,0.305882352941176,0.164705882352941);
  rgb(6pt)=(0.890196078431372,0.101960784313725,0.109803921568627);
  rgb(7pt)=(0.741176470588235,0,0.149019607843137);
  rgb(8pt)=(0.501960784313725,0,0.149019607843137)
},
height=\figureheight,
hide x axis,
hide y axis,
point meta max=1.45636822945787,
point meta min=-1.41023053305731,
tick align=outside,
tick pos=left,
align=center,
title={\tiny{POD-G} \\ \tiny{($\mathrm{p}=12$)}},
width=\figurewidth,
x grid style={darkgray176},
xmin=0, xmax=3200,
xtick style={color=black},
xtick={0,1000,2000,3000,4000},
xticklabels={
  \(\displaystyle {0}\),
  \(\displaystyle {1000}\),
  \(\displaystyle {2000}\),
  \(\displaystyle {3000}\),
  \(\displaystyle {4000}\)
},
y grid style={darkgray176},
ymin=0, ymax=3360,
ytick style={color=black},
ytick={0,500,1000,1500,2000,2500,3000,3500},
yticklabels={
  \(\displaystyle {0}\),
  \(\displaystyle {500}\),
  \(\displaystyle {1000}\),
  \(\displaystyle {1500}\),
  \(\displaystyle {2000}\),
  \(\displaystyle {2500}\),
  \(\displaystyle {3000}\),
  \(\displaystyle {3500}\)
}
]
\addplot graphics [includegraphics cmd=\pgfimage,xmin=0, xmax=3200, ymin=0, ymax=3360] {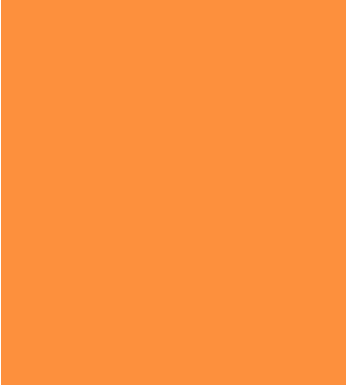};

%%%%%%%%%%%%%%%%%%%%%%%%%%%%%%%%%%%%%%%%%%%%%%%%%%%%%%%%%%%%%%%%%%%%%%%%%%%%%

\nextgroupplot[
colorbar,
colorbar style={ytick={-0.03,-0.02,-0.01,0,0.01,0.02,0.03},yticklabels={
  \(\displaystyle {\ensuremath{-}0.03}\),
  \(\displaystyle {\ensuremath{-}0.02}\),
  \(\displaystyle {\ensuremath{-}0.01}\),
  \(\displaystyle {0.00}\),
  \(\displaystyle {0.01}\),
  \(\displaystyle {0.02}\),
  \(\displaystyle {0.03}\)
},ylabel={}, width=0.06*\pgfkeysvalueof{/pgfplots/parent axis width}, ticklabel style={font=\tiny}},
colormap={mymap}{[1pt]
  rgb(0pt)=(1,1,0.8);
  rgb(1pt)=(1,0.929411764705882,0.627450980392157);
  rgb(2pt)=(0.996078431372549,0.850980392156863,0.462745098039216);
  rgb(3pt)=(0.996078431372549,0.698039215686274,0.298039215686275);
  rgb(4pt)=(0.992156862745098,0.552941176470588,0.235294117647059);
  rgb(5pt)=(0.988235294117647,0.305882352941176,0.164705882352941);
  rgb(6pt)=(0.890196078431372,0.101960784313725,0.109803921568627);
  rgb(7pt)=(0.741176470588235,0,0.149019607843137);
  rgb(8pt)=(0.501960784313725,0,0.149019607843137)
},
height=\figureheight,
hide x axis,
hide y axis,
point meta max=0.0272200884596363,
point meta min=-0.0271155658523269,
tick align=outside,
tick pos=left,
width=\figurewidth,
x grid style={darkgray176},
xmin=0, xmax=3200,
xtick style={color=black},
xtick={0,1000,2000,3000,4000},
xticklabels={
  \(\displaystyle {0}\),
  \(\displaystyle {1000}\),
  \(\displaystyle {2000}\),
  \(\displaystyle {3000}\),
  \(\displaystyle {4000}\)
},
y grid style={darkgray176},
ymin=0, ymax=3360,
ytick style={color=black},
ytick={0,500,1000,1500,2000,2500,3000,3500},
yticklabels={
  \(\displaystyle {0}\),
  \(\displaystyle {500}\),
  \(\displaystyle {1000}\),
  \(\displaystyle {1500}\),
  \(\displaystyle {2000}\),
  \(\displaystyle {2500}\),
  \(\displaystyle {3000}\),
  \(\displaystyle {3500}\)
}
]
\addplot graphics [includegraphics cmd=\pgfimage,xmin=0, xmax=3200, ymin=0, ymax=3360] {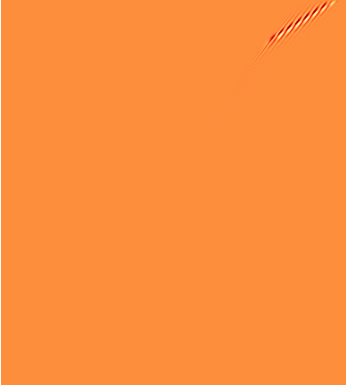};

\nextgroupplot[
colorbar,
colorbar style={ytick={-0.03,-0.02,-0.01,0,0.01,0.02,0.03},yticklabels={
  \(\displaystyle {\ensuremath{-}0.03}\),
  \(\displaystyle {\ensuremath{-}0.02}\),
  \(\displaystyle {\ensuremath{-}0.01}\),
  \(\displaystyle {0.00}\),
  \(\displaystyle {0.01}\),
  \(\displaystyle {0.02}\),
  \(\displaystyle {0.03}\)
},ylabel={}, width=0.06*\pgfkeysvalueof{/pgfplots/parent axis width}, ticklabel style={font=\tiny}},
colormap={mymap}{[1pt]
  rgb(0pt)=(1,1,0.8);
  rgb(1pt)=(1,0.929411764705882,0.627450980392157);
  rgb(2pt)=(0.996078431372549,0.850980392156863,0.462745098039216);
  rgb(3pt)=(0.996078431372549,0.698039215686274,0.298039215686275);
  rgb(4pt)=(0.992156862745098,0.552941176470588,0.235294117647059);
  rgb(5pt)=(0.988235294117647,0.305882352941176,0.164705882352941);
  rgb(6pt)=(0.890196078431372,0.101960784313725,0.109803921568627);
  rgb(7pt)=(0.741176470588235,0,0.149019607843137);
  rgb(8pt)=(0.501960784313725,0,0.149019607843137)
},
height=\figureheight,
hide x axis,
hide y axis,
point meta max=0.0272200884596363,
point meta min=-0.0271155658523269,
tick align=outside,
tick pos=left,
width=\figurewidth,
x grid style={darkgray176},
xmin=0, xmax=3200,
xtick style={color=black},
xtick={0,1000,2000,3000,4000},
xticklabels={
  \(\displaystyle {0}\),
  \(\displaystyle {1000}\),
  \(\displaystyle {2000}\),
  \(\displaystyle {3000}\),
  \(\displaystyle {4000}\)
},
y grid style={darkgray176},
ymin=0, ymax=3360,
ytick style={color=black},
ytick={0,500,1000,1500,2000,2500,3000,3500},
yticklabels={
  \(\displaystyle {0}\),
  \(\displaystyle {500}\),
  \(\displaystyle {1000}\),
  \(\displaystyle {1500}\),
  \(\displaystyle {2000}\),
  \(\displaystyle {2500}\),
  \(\displaystyle {3000}\),
  \(\displaystyle {3500}\)
}
]
\addplot graphics [includegraphics cmd=\pgfimage,xmin=0, xmax=3200, ymin=0, ymax=3360] {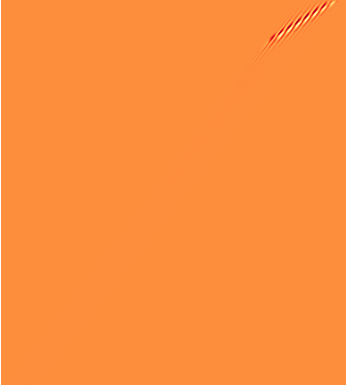};

\nextgroupplot[
colorbar,
colorbar style={ytick={-0.03,-0.02,-0.01,0,0.01,0.02,0.03},yticklabels={
  \(\displaystyle {\ensuremath{-}0.03}\),
  \(\displaystyle {\ensuremath{-}0.02}\),
  \(\displaystyle {\ensuremath{-}0.01}\),
  \(\displaystyle {0.00}\),
  \(\displaystyle {0.01}\),
  \(\displaystyle {0.02}\),
  \(\displaystyle {0.03}\)
},ylabel={}, width=0.06*\pgfkeysvalueof{/pgfplots/parent axis width}, ticklabel style={font=\tiny}},
colormap={mymap}{[1pt]
  rgb(0pt)=(1,1,0.8);
  rgb(1pt)=(1,0.929411764705882,0.627450980392157);
  rgb(2pt)=(0.996078431372549,0.850980392156863,0.462745098039216);
  rgb(3pt)=(0.996078431372549,0.698039215686274,0.298039215686275);
  rgb(4pt)=(0.992156862745098,0.552941176470588,0.235294117647059);
  rgb(5pt)=(0.988235294117647,0.305882352941176,0.164705882352941);
  rgb(6pt)=(0.890196078431372,0.101960784313725,0.109803921568627);
  rgb(7pt)=(0.741176470588235,0,0.149019607843137);
  rgb(8pt)=(0.501960784313725,0,0.149019607843137)
},
height=\figureheight,
hide x axis,
hide y axis,
point meta max=0.0272200884596363,
point meta min=-0.0271155658523269,
tick align=outside,
tick pos=left,
width=\figurewidth,
x grid style={darkgray176},
xmin=0, xmax=3200,
xtick style={color=black},
xtick={0,1000,2000,3000,4000},
xticklabels={
  \(\displaystyle {0}\),
  \(\displaystyle {1000}\),
  \(\displaystyle {2000}\),
  \(\displaystyle {3000}\),
  \(\displaystyle {4000}\)
},
y grid style={darkgray176},
ymin=0, ymax=3360,
ytick style={color=black},
ytick={0,500,1000,1500,2000,2500,3000,3500},
yticklabels={
  \(\displaystyle {0}\),
  \(\displaystyle {500}\),
  \(\displaystyle {1000}\),
  \(\displaystyle {1500}\),
  \(\displaystyle {2000}\),
  \(\displaystyle {2500}\),
  \(\displaystyle {3000}\),
  \(\displaystyle {3500}\)
}
]
\addplot graphics [includegraphics cmd=\pgfimage,xmin=0, xmax=3200, ymin=0, ymax=3360] {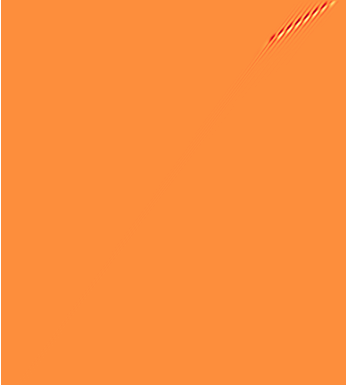};

\nextgroupplot[
colorbar,
colorbar style={ytick={-0.03,-0.02,-0.01,0,0.01,0.02,0.03},yticklabels={
  \(\displaystyle {\ensuremath{-}0.03}\),
  \(\displaystyle {\ensuremath{-}0.02}\),
  \(\displaystyle {\ensuremath{-}0.01}\),
  \(\displaystyle {0.00}\),
  \(\displaystyle {0.01}\),
  \(\displaystyle {0.02}\),
  \(\displaystyle {0.03}\)
},ylabel={}, width=0.06*\pgfkeysvalueof{/pgfplots/parent axis width}, ticklabel style={font=\tiny}},
colormap={mymap}{[1pt]
  rgb(0pt)=(1,1,0.8);
  rgb(1pt)=(1,0.929411764705882,0.627450980392157);
  rgb(2pt)=(0.996078431372549,0.850980392156863,0.462745098039216);
  rgb(3pt)=(0.996078431372549,0.698039215686274,0.298039215686275);
  rgb(4pt)=(0.992156862745098,0.552941176470588,0.235294117647059);
  rgb(5pt)=(0.988235294117647,0.305882352941176,0.164705882352941);
  rgb(6pt)=(0.890196078431372,0.101960784313725,0.109803921568627);
  rgb(7pt)=(0.741176470588235,0,0.149019607843137);
  rgb(8pt)=(0.501960784313725,0,0.149019607843137)
},
height=\figureheight,
hide x axis,
hide y axis,
point meta max=0.0272200884596363,
point meta min=-0.0271155658523269,
tick align=outside,
tick pos=left,
width=\figurewidth,
x grid style={darkgray176},
xmin=0, xmax=3200,
xtick style={color=black},
xtick={0,1000,2000,3000,4000},
xticklabels={
  \(\displaystyle {0}\),
  \(\displaystyle {1000}\),
  \(\displaystyle {2000}\),
  \(\displaystyle {3000}\),
  \(\displaystyle {4000}\)
},
y grid style={darkgray176},
ymin=0, ymax=3360,
ytick style={color=black},
ytick={0,500,1000,1500,2000,2500,3000,3500},
yticklabels={
  \(\displaystyle {0}\),
  \(\displaystyle {500}\),
  \(\displaystyle {1000}\),
  \(\displaystyle {1500}\),
  \(\displaystyle {2000}\),
  \(\displaystyle {2500}\),
  \(\displaystyle {3000}\),
  \(\displaystyle {3500}\)
}
]
\addplot graphics [includegraphics cmd=\pgfimage,xmin=0, xmax=3200, ymin=0, ymax=3360] {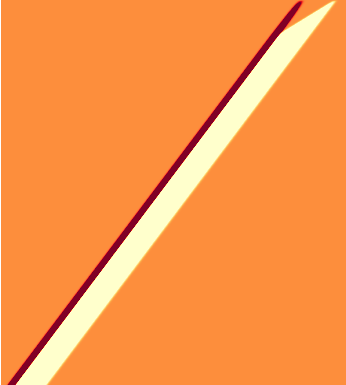};

%%%%%%%%%%%%%%%%%%%%%%%%%%%%%%%%%%%%%%%%%%%%%%%%%%%%%%%%%%%%%%%%%%%%%%%%%%%%%%%%%%%%%%

\nextgroupplot[
colorbar,
colorbar style={ytick={-0.4,-0.2,-2.77555756156289e-17,0.2,0.4,0.6,0.8,1,1.2},yticklabels={
  \(\displaystyle {\ensuremath{-}0.4}\),
  \(\displaystyle {\ensuremath{-}0.2}\),
  \(\displaystyle {0.0}\),
  \(\displaystyle {0.2}\),
  \(\displaystyle {0.4}\),
  \(\displaystyle {0.6}\),
  \(\displaystyle {0.8}\),
  \(\displaystyle {1.0}\),
  \(\displaystyle {1.2}\)
},ylabel={}, width=0.06*\pgfkeysvalueof{/pgfplots/parent axis width}, ticklabel style={font=\tiny}},
colormap={mymap}{[1pt]
  rgb(0pt)=(1,1,0.8);
  rgb(1pt)=(1,0.929411764705882,0.627450980392157);
  rgb(2pt)=(0.996078431372549,0.850980392156863,0.462745098039216);
  rgb(3pt)=(0.996078431372549,0.698039215686274,0.298039215686275);
  rgb(4pt)=(0.992156862745098,0.552941176470588,0.235294117647059);
  rgb(5pt)=(0.988235294117647,0.305882352941176,0.164705882352941);
  rgb(6pt)=(0.890196078431372,0.101960784313725,0.109803921568627);
  rgb(7pt)=(0.741176470588235,0,0.149019607843137);
  rgb(8pt)=(0.501960784313725,0,0.149019607843137)
},
height=\figureheight,
hide x axis,
hide y axis,
point meta max=1.05770831601777,
point meta min=-0.249629751620698,
tick align=outside,
tick pos=left,
width=\figurewidth,
x grid style={darkgray176},
xmin=0, xmax=3200,
xtick style={color=black},
xtick={0,1000,2000,3000,4000},
xticklabels={
  \(\displaystyle {0}\),
  \(\displaystyle {1000}\),
  \(\displaystyle {2000}\),
  \(\displaystyle {3000}\),
  \(\displaystyle {4000}\)
},
y grid style={darkgray176},
ymin=0, ymax=3360,
ytick style={color=black},
ytick={0,500,1000,1500,2000,2500,3000,3500},
yticklabels={
  \(\displaystyle {0}\),
  \(\displaystyle {500}\),
  \(\displaystyle {1000}\),
  \(\displaystyle {1500}\),
  \(\displaystyle {2000}\),
  \(\displaystyle {2500}\),
  \(\displaystyle {3000}\),
  \(\displaystyle {3500}\)
}
]
\addplot graphics [includegraphics cmd=\pgfimage,xmin=0, xmax=3200, ymin=0, ymax=3360] {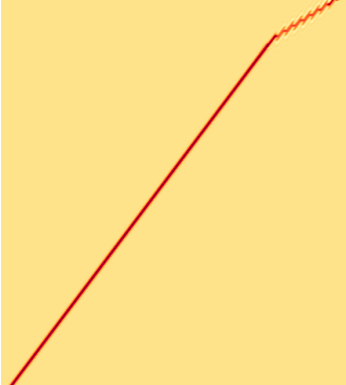};

\nextgroupplot[
colorbar,
colorbar style={ytick={-0.4,-0.2,-2.77555756156289e-17,0.2,0.4,0.6,0.8,1,1.2},yticklabels={
  \(\displaystyle {\ensuremath{-}0.4}\),
  \(\displaystyle {\ensuremath{-}0.2}\),
  \(\displaystyle {0.0}\),
  \(\displaystyle {0.2}\),
  \(\displaystyle {0.4}\),
  \(\displaystyle {0.6}\),
  \(\displaystyle {0.8}\),
  \(\displaystyle {1.0}\),
  \(\displaystyle {1.2}\)
},ylabel={}, width=0.06*\pgfkeysvalueof{/pgfplots/parent axis width}, ticklabel style={font=\tiny}},
colormap={mymap}{[1pt]
  rgb(0pt)=(1,1,0.8);
  rgb(1pt)=(1,0.929411764705882,0.627450980392157);
  rgb(2pt)=(0.996078431372549,0.850980392156863,0.462745098039216);
  rgb(3pt)=(0.996078431372549,0.698039215686274,0.298039215686275);
  rgb(4pt)=(0.992156862745098,0.552941176470588,0.235294117647059);
  rgb(5pt)=(0.988235294117647,0.305882352941176,0.164705882352941);
  rgb(6pt)=(0.890196078431372,0.101960784313725,0.109803921568627);
  rgb(7pt)=(0.741176470588235,0,0.149019607843137);
  rgb(8pt)=(0.501960784313725,0,0.149019607843137)
},
height=\figureheight,
hide x axis,
hide y axis,
point meta max=1.05770831601777,
point meta min=-0.249629751620698,
tick align=outside,
tick pos=left,
width=\figurewidth,
x grid style={darkgray176},
xmin=0, xmax=3200,
xtick style={color=black},
xtick={0,1000,2000,3000,4000},
xticklabels={
  \(\displaystyle {0}\),
  \(\displaystyle {1000}\),
  \(\displaystyle {2000}\),
  \(\displaystyle {3000}\),
  \(\displaystyle {4000}\)
},
y grid style={darkgray176},
ymin=0, ymax=3360,
ytick style={color=black},
ytick={0,500,1000,1500,2000,2500,3000,3500},
yticklabels={
  \(\displaystyle {0}\),
  \(\displaystyle {500}\),
  \(\displaystyle {1000}\),
  \(\displaystyle {1500}\),
  \(\displaystyle {2000}\),
  \(\displaystyle {2500}\),
  \(\displaystyle {3000}\),
  \(\displaystyle {3500}\)
}
]
\addplot graphics [includegraphics cmd=\pgfimage,xmin=0, xmax=3200, ymin=0, ymax=3360] {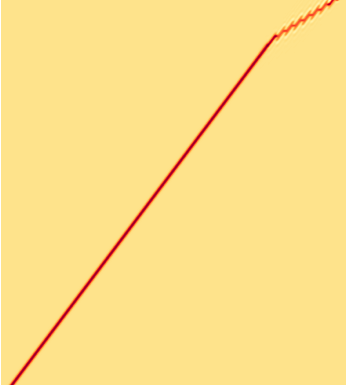};

\nextgroupplot[
colorbar,
colorbar style={ytick={-0.4,-0.2,-2.77555756156289e-17,0.2,0.4,0.6,0.8,1,1.2},yticklabels={
  \(\displaystyle {\ensuremath{-}0.4}\),
  \(\displaystyle {\ensuremath{-}0.2}\),
  \(\displaystyle {0.0}\),
  \(\displaystyle {0.2}\),
  \(\displaystyle {0.4}\),
  \(\displaystyle {0.6}\),
  \(\displaystyle {0.8}\),
  \(\displaystyle {1.0}\),
  \(\displaystyle {1.2}\)
},ylabel={}, width=0.06*\pgfkeysvalueof{/pgfplots/parent axis width}, ticklabel style={font=\tiny}},
colormap={mymap}{[1pt]
  rgb(0pt)=(1,1,0.8);
  rgb(1pt)=(1,0.929411764705882,0.627450980392157);
  rgb(2pt)=(0.996078431372549,0.850980392156863,0.462745098039216);
  rgb(3pt)=(0.996078431372549,0.698039215686274,0.298039215686275);
  rgb(4pt)=(0.992156862745098,0.552941176470588,0.235294117647059);
  rgb(5pt)=(0.988235294117647,0.305882352941176,0.164705882352941);
  rgb(6pt)=(0.890196078431372,0.101960784313725,0.109803921568627);
  rgb(7pt)=(0.741176470588235,0,0.149019607843137);
  rgb(8pt)=(0.501960784313725,0,0.149019607843137)
},
height=\figureheight,
hide x axis,
hide y axis,
point meta max=1.05770831601777,
point meta min=-0.249629751620698,
tick align=outside,
tick pos=left,
width=\figurewidth,
x grid style={darkgray176},
xmin=0, xmax=3200,
xtick style={color=black},
xtick={0,1000,2000,3000,4000},
xticklabels={
  \(\displaystyle {0}\),
  \(\displaystyle {1000}\),
  \(\displaystyle {2000}\),
  \(\displaystyle {3000}\),
  \(\displaystyle {4000}\)
},
y grid style={darkgray176},
ymin=0, ymax=3360,
ytick style={color=black},
ytick={0,500,1000,1500,2000,2500,3000,3500},
yticklabels={
  \(\displaystyle {0}\),
  \(\displaystyle {500}\),
  \(\displaystyle {1000}\),
  \(\displaystyle {1500}\),
  \(\displaystyle {2000}\),
  \(\displaystyle {2500}\),
  \(\displaystyle {3000}\),
  \(\displaystyle {3500}\)
}
]
\addplot graphics [includegraphics cmd=\pgfimage,xmin=0, xmax=3200, ymin=0, ymax=3360] {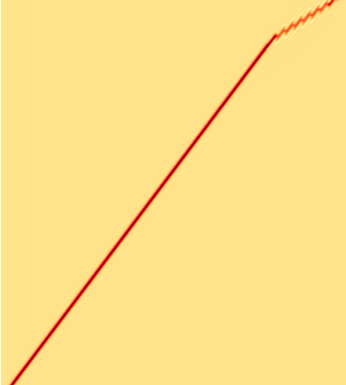};

\nextgroupplot[
colorbar,
colorbar style={ytick={-0.4,-0.2,-2.77555756156289e-17,0.2,0.4,0.6,0.8,1,1.2},yticklabels={
  \(\displaystyle {\ensuremath{-}0.4}\),
  \(\displaystyle {\ensuremath{-}0.2}\),
  \(\displaystyle {0.0}\),
  \(\displaystyle {0.2}\),
  \(\displaystyle {0.4}\),
  \(\displaystyle {0.6}\),
  \(\displaystyle {0.8}\),
  \(\displaystyle {1.0}\),
  \(\displaystyle {1.2}\)
},ylabel={}, width=0.06*\pgfkeysvalueof{/pgfplots/parent axis width}, ticklabel style={font=\tiny}},
colormap={mymap}{[1pt]
  rgb(0pt)=(1,1,0.8);
  rgb(1pt)=(1,0.929411764705882,0.627450980392157);
  rgb(2pt)=(0.996078431372549,0.850980392156863,0.462745098039216);
  rgb(3pt)=(0.996078431372549,0.698039215686274,0.298039215686275);
  rgb(4pt)=(0.992156862745098,0.552941176470588,0.235294117647059);
  rgb(5pt)=(0.988235294117647,0.305882352941176,0.164705882352941);
  rgb(6pt)=(0.890196078431372,0.101960784313725,0.109803921568627);
  rgb(7pt)=(0.741176470588235,0,0.149019607843137);
  rgb(8pt)=(0.501960784313725,0,0.149019607843137)
},
height=\figureheight,
hide x axis,
hide y axis,
point meta max=1.05770831601777,
point meta min=-0.249629751620698,
tick align=outside,
tick pos=left,
width=\figurewidth,
x grid style={darkgray176},
xmin=0, xmax=3200,
xtick style={color=black},
xtick={0,1000,2000,3000,4000},
xticklabels={
  \(\displaystyle {0}\),
  \(\displaystyle {1000}\),
  \(\displaystyle {2000}\),
  \(\displaystyle {3000}\),
  \(\displaystyle {4000}\)
},
y grid style={darkgray176},
ymin=0, ymax=3360,
ytick style={color=black},
ytick={0,500,1000,1500,2000,2500,3000,3500},
yticklabels={
  \(\displaystyle {0}\),
  \(\displaystyle {500}\),
  \(\displaystyle {1000}\),
  \(\displaystyle {1500}\),
  \(\displaystyle {2000}\),
  \(\displaystyle {2500}\),
  \(\displaystyle {3000}\),
  \(\displaystyle {3500}\)
}
]
\addplot graphics [includegraphics cmd=\pgfimage,xmin=0, xmax=3200, ymin=0, ymax=3360] {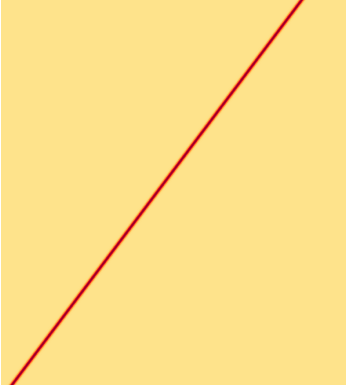};

\end{groupplot}

\draw ({$(current bounding box.south west)!-0.04!(current bounding box.south east)$}|-{$(current bounding box.south west)!0.4!(current bounding box.north west)$}) node[
  scale=0.96,
  anchor=west,
  text=black,
  rotate=90.0
]{time $t$};
\draw ({$(current bounding box.south west)!0.5!(current bounding box.south east)$}|-{$(current bounding box.south west)!-0.06!(current bounding box.north west)$}) node[
  scale=0.96,
  anchor=south,
  text=black,
  rotate=0.0
]{space $x$};
\end{tikzpicture}